\documentclass[a4paper,dvips]{book} 

\usepackage{pstricks,pst-node,pst-plot}

\usepackage{ifthen}

\usepackage{latexsym}
\usepackage{amsfonts}
\usepackage[reqno]{amsmath}
\usepackage{amsthm}
\usepackage{amssymb}
\usepackage{amsopn}	

\usepackage{bbm}			

\usepackage{mathabx}

\usepackage{makeidx}
\usepackage[all]{xy}			
\usepackage{mathrsfs}
\usepackage[nottoc]{tocbibind}		%

\usepackage[refpage]{nomencl}
\renewcommand{\nomgroup}[1]{%
    \ifthenelse{\equal{#1}{G}}{\item[\textbf{Lie groups and algebras}]}{}%
    \ifthenelse{\equal{#1}{F}}{\item[\textbf{Functional analysis}]}{}%
    \ifthenelse{\equal{#1}{D}}{\item[\textbf{Differential geometry}]}{}%
    \ifthenelse{\equal{#1}{A}}{\item[\textbf{Algebra}]}{}%
}

\usepackage[utf8]{inputenc}
\usepackage[T1]{fontenc}

\usepackage{textcomp}
\usepackage{lmodern}
\usepackage[a4paper]{geometry}

\usepackage{textcomp}
\usepackage{lmodern}

\usepackage[ps2pdf]{hyperref} 			
\hypersetup{colorlinks=true,linkcolor=blue}	
\usepackage[numbers,sort&compress]{natbib}

\newcounter{numtho}[chapter]
\renewcommand{\thenumtho}{\arabic{chapter}.\arabic{numtho}}

\newtheoremstyle{mes_tho}%
		{9pt}{9pt}%
		{\itshape}%
		{}%
		{\bfseries}{.}%
		{\newline}%
		{}%

\renewcommand*{\theenumi}{(\roman{enumi})}

\theoremstyle{remark}	\newtheorem{remark}[numtho]{Remark}
			
\theoremstyle{mes_tho}	\newtheorem{lemma}[numtho]{Lemma}
			\newtheorem{theorem}[numtho]{Theorem}
			\newtheorem{corollary}[numtho]{Corollary}
			\newtheorem{proposition}[numtho]{Proposition}
			\newtheorem{definition}[numtho]{Definition}

\newcommand{\arxiv}[2][]{
\newline
\ifthenelse{\equal{#1}{}}{%
	\href{http://xxx.lanl.gov/abs/#2}{{\tt arxiv:#2}}%
			}
			{%
	\href{http://xxx.lanl.gov/abs/#2}{{\tt arxiv:#2}[#1]}%
}
}				%

\newcommand\abstractname{Abstract}
\makeatletter
  \newenvironment{abstract}{%
      \if@twocolumn
        \section*{\abstractname}%
      \else
        \small
        \begin{center}%
          {\bfseries \abstractname\vspace{-.5em}\vspace{\z@}}%
        \end{center}%
        \quotation
      \fi}
      {\if@twocolumn\else\endquotation\fi}
\makeatother

\newcounter{bidon}

\newcommand{\defe}[2]{\textbf{#1}\index{#2}}
\newcommand{\us}[1]{\frac{1}{#1}}
\newcommand{\me}[1]{(-1)^{#1}}
\newcommand{\dpt}[3]{#1\colon #2\to #3}
\newcommand{\dsdd}[3]{\left.\frac{d}{d#2}#1\right|_{#2=#3}}
\newcommand{\Dsddb}[4]{\frac{d}{d#2}\Big[#1\Big]_{#3=#4}}
\newcommand{\Dsdd}[3]{ \Dsddb{#1}{#2}{#2}{#3}   }
\newcommand{\Dsddp}[3]{\frac{d}{d#2}\Big(#1\Big)_{#2=#3}}

\newcommand{\cvec}{\mathfrak{X}}
\newcommand{\LogOu}{\vee}
\newcommand{\LogEt}{\wedge}
		
\newcommand{\subdem}[1]{\par\noindent {\it #1.} }
\newcommand{\BX}{{\bf X}}
\newcommand{\tq}{\mid}
\newcommand{\scal}[2]{ \langle #1|#2 \rangle }		%
\newcommand{\mue}[1]{(-1)^{#1}}

\newcommand{\RM}{\pr_{\sQ}\sR}
\newcommand{\dcr}[1]{[[#1]]}

\DeclareMathOperator{\SU}{SU}		
\DeclareMathOperator{\SL}{SL}
\DeclareMathOperator{\SO}{SO}			
\DeclareMathOperator{\SP}{SP}			
\DeclareMathOperator{\gO}{O}			%
\DeclareMathOperator{\gsl}{\mathfrak{sl}}
\DeclareMathOperator{\gsp}{\mathfrak{sp}}	
\DeclareMathOperator{\so}{\mathfrak{so}}

\DeclareMathOperator{\gpAff}{Aff}
\DeclareMathOperator{\gpSymp}{Symp}
\DeclareMathOperator{\Spin}{Spin}
\DeclareMathOperator{\pr}{\texttt{pr}}
\DeclareMathOperator{\tr}{Tr}
\DeclareMathOperator{\Fun}{\texttt{Fun}}
\DeclareMathOperator{\ad}{ad}
\DeclareMathOperator{\Ad}{Ad}
\DeclareMathOperator{\AD}{\textbf{Ad}}
\DeclareMathOperator{\Der}{\texttt{Der}}
\DeclareMathOperator{\Aut}{Aut}
\DeclareMathOperator{\Span}{Span}
\DeclareMathOperator{\Iso}{\texttt{Iso}}
\DeclareMathOperator{\Int}{Int}
\DeclareMathOperator{\mfsp}{\mathfrak{sp}}
\DeclareMathOperator{\Fr}{Fr}
\DeclareMathOperator{\End}{End}
\DeclareMathOperator{\id}{id}

\newcommand{\mfo}{\vartheta}
\def\bnu{_{\nu}}

\def\mfs{\mathfrak{s}}

\newcommand{\hH}{\mathscr{H}}
\newcommand{\hS}{\mathscr{S}}

\def\pH{\mathscr{H}}

\newcommand{\swS}{\mathscr{S}}	
\newcommand{\swE}{\mathscr{E}}

\def\lA{\mathfrak{a}}
\def\lA{\mathfrak{a}}

\def\lG{\mathfrak{g}}		
\def\lH{\mathfrak{h}}

\def\lK{\mathfrak{k}}

\def\lN{\mathfrak{n}}
\def\lP{\mathfrak{p}}
\def\lQ{\mathfrak{q}}

\def\lZ{\mathcal{Z}}

\def\iA{\mathcal{A}}
\def\iK{\mathcal{K}}
\def\iN{\mathcal{N}}			%
\def\iP{\mathcal{P}}
\def\iR{\mathcal{R}}
\def\iAH{\mathcal{A_H}}
\def\iKH{\mathcal{K_H}}

\def\iNH{\mathcal{N_H}}
\def\iPH{\mathcal{P_H}}
\def\iRH{\mathcal{R_H}}

\newcommand{\SUR}{\mathrm{R}}
\newcommand{\SUA}{\mathrm{A}}		%
\newcommand{\SUN}{\mathrm{N}}

\def\sA{\mathcal{A}}
\def\sG{\mathcal{G}}
\def\sH{\mathcal{H}}			%
\def\sK{\mathcal{K}}			
		
\def\sN{\mathcal{N}}	
\def\sP{\mathcal{P}}
\def\sQ{\mathcal{Q}}
\def\sR{\mathcal{R}}
\newcommand{\sS}{\mathcal{S}}

\def\mA{\mathcal{A}}

\def\mF{\mathcal{F}}
\def\mG{\mathcal{G}}

\def\mI{\mathcal{I}}
\def\mJ{\mathcal{J}}
\def\mK{\mathcal{K}}
\def\mL{\mathcal{L}}

\def\mN{\mathcal{N}}
\def\mO{\mathcal{O}}

\def\mR{\mathcal{R}}
\def\mS{\mathcal{S}}
\def\mT{\mathcal{T}}
\newcommand\mU{\mathcal{U}}

\def\eA{\mathbbm{A}}
\def\eC{\mathbbm{C}}

\def\eN{\mathbbm{N}}

\def\eR{\mathbbm{R}}
\def\eZ{\mathbbm{Z}}

\newcommand{\mtu}{\mathbbm{1}}  			%

\def\sod{\mathfrak{so}(2)}
\def\SOdn{\SO(2,n)}
\newcommand{\sodn}{  {\mathfrak{so}}(2,n)   }
\def\soun{\mathfrak{so}(1,n)}
\def\sldr{\mathfrak{sl}(2,\eR)}
\def\SOun{\SO(1,n)}
\def\Cinf{C^{\infty}}

\def\tf{\tilde{f}}

\def\hu{\hat{u}}
\def\hv{\hat{v}}

\def\stG{\star^{G}}

\def\stX{\star^{X}}
\def\stM{\ast_M}
\def\stt{\star_{\theta}}

\def\fge{''}
\def\oge{``}

\def\AmS{\protect\pAmS}            \def\LaTeX{\protect\pLaTeX}
\def\pAmS{{\the\textfont2
        A\kern-.1667em\lower.5ex\hbox{M}\kern-.125emS}}
\def\pLaTeX{{\rm L\kern-.36em\raise.3ex\hbox{\the\scriptfont0 A}\kern-.15em
    T\kern-.1667em\lower.7ex\hbox{E}\kern-.125emX}}

\def\amstex/{\AmS-\TeX}         \def\amslatex/{\AmS-\LaTeX{}}
\def\latex/{\LaTeX{}}           \def\pictex/{PIC\TeX}
\def\tex/{\TeX}
\def\bibtex/{{\sc Bib\kern-.1em\TeX}}     \def\tugboat/{{\it TUGboat\/}}
\def\amsfonts/{AMSFonts}

   \def\dbar{\leavevmode\hbox to 0pt{\hskip.2ex
    \accent"16\hss}d}

\makeindex

\begin{document}

\pagestyle{empty}

\newcommand{\TitreThese}{Locally anti de Sitter spaces and deformation quantization}

\title{\TitreThese}
\author{Laurent Claessens}
\maketitle

\newpage
\thispagestyle{empty}

\noindent{\huge \TitreThese}\\ 
Laurent Claessens\footnote{Supported by FRIA, Belgium.}\\
September 2007\footnote{This version has some minor modification with respect to the original publication which can be downloaded here: \url{http://hdl.handle.net/2078.1/5354}.}
\vspace{2cm}

\noindent\textsc{Promoteur :} P. Bieliavsky\\

\vspace{1cm}

\noindent\textsc{Jury :}\\
Yoshiaki Maeda\\
Philippe Bonneau\\
Giuseppe  Dito\\
Jean Mawhin (président)\\
Pascal Lambrechts\\
Yves Felix

\chapter*{Remerciements}

Étant donné qu'il est d'usage de terminer par les choses importantes je ne vais pas commencer par Pierre Bieliavsky, l'homme qui voit des fugues de Bach dès qu'un calcul fait plus de deux lignes. Parlons donc de celui qui tapote des impros de jazz à longueur de calculs, celui qui a pu Yanismer nombre des expressions de Pierre. Remercions-le pour avoir écouté des heures durant les fautes que je lui expliquais au tableau\footnote{Je crois qu'il y a encore un problème de signe en bas à droite.}. 

Le p'ti rigolo de cousin kapitalist doit être remercié pour sa collaboration BTZ, son organisation de Modave\footnote{À propos, je t'ai attendu toute la journée hier à la gare de Bastogne; tu m'as oublié ou quoi ?} (merci aussi aux autres) mais aussi et surtout parce que non content d'être un p'ti rigolo, il est aussi un grand marrant avec qui j'ai passé de très bons moments !

Que l'ensemble de la communauté libre reçoive ma gratitude pour m'avoir fourni des logiciels de qualité qui m'ont permis de mener ma thèse à bien  : \LaTeX, ubuntu, vim, kde, maxima et bien d'autres. Dans la même veine, je voudrais remercier toutes celles et ceux qui ont mis leurs articles sur arXiv pour me permettre d'y accéder plus facilement.

Il y a aussi des choses à dire du côté de ma famille pour la bonne éducation prodiguée (mais dont il est facile de ne pas profiter), pour le désencombrement des meubles en trop, son accueil et sa facilité à me reconnaître sur une photo.

Merci aussi à la petite s\oe ur que j'ai connue anti-\TeX tile, plus chaude pour aller au Canada que sur la banquise et qui ne fut pas comique à convertir; je lui dois huit ans\footnote{d'où il faut décompter les heures de labos en première.} d'amitié.

I want to add an acknowledgement for Ping Xu and Abhay Ashtekar who accepted me for a postdoc in Penn States next year. Thank you Ping for the housing and Matthieu for your helpful collaboration.

Beaucoup de professeurs m'ont appris des choses que je devais savoir tout au long de mes études de physique et de mathématique à différents niveaux~: Daniel Noul, Pierre Gaspard, 
Christianne Schomblond, Philippe Boulanger, Marc Henneaux, Glenn Barnich, Laurent Houart, Simone Gutt, Philippe Spindel et bien d'autres qui ont émaillé ma vie de bons conseils.

Cela étant fait, parlons un peu de Pierre qui a su me montrer que dans les math, il n'y a pas que la physique. Il y a aussi la beauté de l'objet mathématique pour lui même, la cuisine asiatique, le calcul sans coordonnées et bien entendu, les fugues de Bach. Merci Pierre pour tout ce que tu m'as apporté durant ces quatre années.

\tableofcontents

\chapter*{Introduction}
\addcontentsline{toc}{chapter}{Introduction}

The question arises when one watches movies such as Star Trek: what is a black hole ? One knows from special relativity that light speed cannot be exceeded. So, as a first attempt to define the notion of black hole, we just say that it is a region of the space from which even light cannot escape. Such an object causes a scientific problem because it is by assumption not observable. This fact allows science-fiction writers to invent anything without any chance of contradiction. That is a Star Trek black hole.

Physical black holes are much more interesting because they are the signal of a general relativity failure.

The Newtonian gravitational field is given by a potential which increases as $1/r$ when you get closer to a massive object. At $r=0$, this potential makes no sense and physics is in trouble. One can avoid the problem by postulating that there exist no pointwise masses and that particles cannot penetrate each other. From these assumptions, the fact that Newtonian mechanics does not impose any limit speed makes the divergence at $r=0$ unimportant.

In general relativity, the divergence at small distances is much more problematic because there is a limit speed; hence a pointwise mass always creates a whole region from which nothing (not even light) can escape. Worse: even a homogeneous ball produces a divergence in the metric when it is too dense, and such objects may exist in the real world. Stated in a more mathematical way: solutions of Einstein's equations for the real world may be singular. From this point of view, black holes are nothing else than a feature in the mathematical framework of relativity which indicates that this is not a final theory. That is the notion of black hole in general relativity and in cosmology.

The transfer of concept from physics to mathematics always consists in taking the key features of the mathematical framework of a physical theory and posing them as definition of a new mathematical object. What are the main mathematical points in the concept of black hole in general relativity ? 
First, we retain the notion of pseudo-Riemannian manifold. The sign of the norm of a vector is the crucial property which allows one to define causality (the light cone).

The second main feature that we extract from the physical situation is the fact that a general relativity black hole has a non empty interior. We saw that this is the key difference between the Newtonian case in which all points are equivalent except the unique point where the mass lies, and the general relativistic case in which a whole region was causally disconnected from the rest of the space.

More precisely, as mathematicians, we ask a black hole to separate the pseudo-Riemannian manifold into two causally disconnected parts in the sense that no light-like geodesics can reach the second region from the first one. Notice that we do not include metric singularity in our mathematical black hole notion. In cosmology, in contrast, black holes always take root in a divergence of some metric invariant such as the curvature.

The anti de Sitter space is a solution of Einstein's equations with constant negative curvature. We consider this space as our framework. First, we define as \emph{singular} the closed orbits of the action of some subgroup of the isometry group $\SO(2,n)$ of anti de Sitter. This is done in such a way to generalize to any dimensions the celebrated BTZ black hole. Then we prove that the resulting structure is a black hole in the sense that it cuts the space into two parts : an interior region from which every light-like geodesic intersects the singularity and an exterior region in which every point accepts at least one light-like geodesics which does not intersect the singularity.  Notice that our black hole does not present any curvature singularity.

The second theme of this thesis is deformation quantization. The key ingredient of quantum mechanics is the noncommutativity of quantum observables. When one tries to measure the velocity and the position of a classical particle (such as a tennis ball or a planet), one can choose the order of measurement. It does not matter which of velocity or position is measured first. Quantum mechanics (the mechanics which governs subatomic particles) is very different. If you measure the position of an electron and then you measure its velocity, you do not get the same result as if you had measured the velocity first and then the position. That noncommutativity in measurements is the very foundation of the quantum mechanics. In the usual mathematical framework, it is implemented by describing each measurable quantity by an operator acting on a Hilbert space. The eigenvalues of these operators correspond to physical measurements. The position and momentum operators for example are respectively $f(x)\mapsto xf(x)$ and $f(x)\mapsto -i\hbar(\partial_xf)(x)$. These two operators obviously do not commute.

 In a more abstract way, we say that noncommutativity of quantum mechanics is implemented by considering some noncommutative algebra of operators acting on a Hilbert space, while the classical mechanics deals with observables that are usual functions that form a commutative algebra. The procedure to pass from commutative function algebras to noncommutative operator algebras is the so-called \emph{quantization} in physics.

In our sense\footnote{Quantization is a very large field of mathematics; as far I know, the idea of noncommutativity is always present, but precise notion of ``to quantize something'' may vary from one subject to another}, \emph{deforming} a manifold is simply putting a one-parameter family of new noncommutative products on the set of functions on this manifold. We impose that it reduces to the usual commutative product when the parameter goes to zero. In order to speak of \emph{quantization}, we ask the first order term in the expansion with respect to the parameter to somehow ``contain'' the symplectic structure given on the original manifold.

Questions that arise in this context are: is it possible to study causality in a noncommutative framework ? does it apply to real physics ?

The main result of the present work is not to directly address these large questions, but to build a concrete example in which one can work. Namely, we consider the anti de Sitter space --- that is the simplest non trivial solution of Einstein's equations with constant negative curvature --- that we endow with a black hole structure defined from the action of a subgroup of the isometry group. Then we select the physical part of the space --- the one which is causally connected to infinity --- and we perform a deformation of that part.

The work is divided into three main parts. In a first  time (chapter \ref{ChapAdS}) we define a ``BTZ'' black hole in anti de Sitter space in any dimension. That will be done by means of group theoretical and symmetric spaces considerations. A physical ``good domain'' is identified as an open orbit of a subgroup of the isometry group of anti de Sitter. 

Then (chapter \ref{ChDefoBH}) we show that the open orbit is in fact isomorphic to a group (we introduce the notion of \emph{globally group type} manifold) for which a quantization exists. The quantization of the black hole is performed and its Dirac operator is computed.

The third part (appendix \ref{ChapDefo} and \ref{ChapTool}) exposes some previously known results. Appendix \ref{ChapDefo} is given in a pedagogical purpose: it exposes generalities about deformation quantization and  careful examples with $\SL(2,\eR)$ and split extensions of Heisenberg algebras. Appendix \ref{ChapTool} is devoted to some classical results about homogeneous spaces and Iwasawa decompositions. Explicit decompositions are given for every algebra that will be used in the thesis. It serves to make the whole text more self contained and to fix notations. Basics of quantization by group action are given in appendix \ref{SecDefAction}. 

One more chapter is inserted (chapter \ref{ChapNoteDev}). It contains two small results which have no true interest by themselves but which raise questions and call for further development. We discuss a product on the half-plane (or, equivalently, on the Iwasawa subgroup of $\SL(2,\eR)$) due to A. Unterberger. We show that the \emph{quantization by group action} machinery can be applied to this product in order to deform the dual of the Lie algebra of that Iwasawa subgroup. Although this result seems promising, we show by two examples that the product is not universal in the sense that even the product of compactly supported functions cannot be defined on $AdS_2$ by the quantization induced by Unterberger's product. 

 Then we show that the Iwasawa subgroup of $\SO(2,n)$ (i.e. the group which defines the singularity) is a symplectic split extension of the Iwasawa subgroup of $SU(1,1)$ by the Iwasawa subgroup of $SU(1,n)$. A quantization of the two latter groups being known, a quantization of $SO(2,n)$ is in principle possible using an extension lemma (subsection \ref{SubSecExtLem}). Properties of this product and the resulting quantization of $AdS_l$ were not investigated because we found a more economical way to quantize $AdS_4$.

\pagestyle{headings}

\chapter{Black holes in anti de Sitter spaces}
\label{ChapAdS}

\begin{abstract}
This chapter deals with black holes in anti de Sitter spaces. The latter are the simplest non flat solutions to Einstein's equations with constant negative cosmological constant; they are in particular pseudo-Riemannian manifolds that carry a causal structure, physically due to the finiteness of speed of light. That physical restriction is mathematically encoded by the existence of three types of geodesics: the space-, time- and light-like ones, existence which is in turn implied by the non positivity of the metric. A causal structure is introduced by defining two points as \emph{causally connected} when there exists a time- or light-like path connecting them.

 The  originality of our approach is that the $l$-dimensional space $AdS_l$ is seen as a quotient of groups $\SO(2,l-1)/\SO(1,l-1)=G/H$, and that the special causal black hole structure is described in terms of orbits of the action of a subgroup of the isometry group of the manifold.

Using symmetric spaces techniques, we show that closed orbits of the Iwasawa subgroup of $\SO(2,l-1)$ naturally define a causal black hole singularity in anti de Sitter spaces in $l \geq 3$ dimensions. In particular, we recover for $l=3$ the non-rotating massive BTZ black hole. The method presented here is very simple and in principle generalizable to any semisimple symmetric space.

\end{abstract}

    \section{Introduction}

\subsection{General ideas of a black hole}

The basic notions needed in order to define a causal structure on a time orientable pseudo-Riemannian manifold are that of time-, light- and space-like tangent vector. A tangent vector is said to be respectively \emph{time-}, \emph{space-} or \emph{light-like} when its norm is positive, negative or null; physically, only time-like vectors are allowed to be the velocity of an observer (this is the fact that light speed cannot be attained by a massive particle), and it is only possible for massless particle (such as photons) to follow trajectories with light-like tangent vectors.

From a geometric point of view, a black hole is the data of a causal manifold $M$ together with a subset $\hS \subset M$ called \emph{singularity} such that the whole manifold is divided into two parts: the \emph{interior} and the \emph{exterior} of the black hole. A point is said to be \emph{interior} if all future light-like geodesics through the point have a non empty intersection with the singularity. A point is \emph{exterior} if it is not interior. An important subset of the space is the \emph{event horizon}: the boundary between these two subsets.

\subsection{BTZ black hole}		\index{BTZ black hole}

The BTZ black hole introduced and developed by Bañados, Teitelbaum, Zannelli and Henneaux in \cite{BTZ_un,BTZ_deux} is an example of a black hole whose singularity is not motivated by metric divergences\footnote{It turns out that general relativity accepts a lot of solutions presenting metric divergences; or more precisely, there are a lot of \emph{physical situations} from which Einstein's equations lead to divergences of some metric invariant such as the curvature.}. The construction is roughly as follows. We consider the anti de Sitter space in which we pick up a Killing vector field whose sign of norm is not constant. Then we perform a \emph{discrete} quotient along the integral curves of this vector field. Of course we obtain a lot of closed geodesics. The point is that, in the region where the Killing vector field is space-like, these closed curves are space-like. That violates the physical principle of causality. For that reason, we decree that this region is singular or, equivalently, that the boundary of this region is singular. The BTZ singularity is then the loci where the chosen Killing vector field has a vanishing norm. Since discrete quotients do not affect local structures, the resulting space remains a solution of the $(2+1)$-dimensional general relativity with negative cosmological constant\footnote{For honesty, we have to warn the reader that the real world's cosmological constant has been measured very small but positive. We also have to point out that the four dimensional anti de Sitter space is a solution of general relativity \emph{without masses}. From a physical point of view, this thesis has to be seen as a toy model.}. In this context one can define pertinent notions of  \emph{mass} and \emph{angular momentum} which depend on the chosen Killing vector field.

 In the case of the \emph{non-rotating massive} BTZ black hole, the structure of the singularity and the horizon are closely related to the action of a minimal parabolic (Iwasawa) subgroup of the isometry group of anti de Sitter, see \cite{BTZB_deux,Keio}. The whole work on the BTZ black hole and the fact that it belongs to the class of causal symmetric spaces (for definitions and some examples, see \cite{HilgertOlaf}) motivate the following definition:

\begin{definition}
A \defe{causal solvable symmetric black hole}{Causal solvable symmetric black hole} is a causal symmetric space where the closed orbits of minimal parabolic subgroups of its isometry group define a black hole singularity. See section \ref{SecCausal} for definitions of causality and singularity.
\label{Def1}
\end{definition}
 
\subsection{Generalization and group setting}

The original BTZ black hole was constructed in dimension three, but we will see in this chapter that, exploiting their group theoretical description, they can easily be generalized to any dimension, as pointed out in \cite{BDRS,lcTNAdS}.  Notice that higher-dimensional generalizations of the BTZ construction have been studied in the physics literature, by classifying the one-parameter isometry subgroups of $\Iso(AdS_l)=\SO(2,l-1)$, see \cite{Figueroa,AdSBH,Madden,Banados:1997df,Aminneborg,HolstPeldan}, but these approaches do not exploit the symmetric space structure of anti de Sitter.

The structure that will be described with full details in next pages may be summarized as follows. Take $G=\SO(2,l-1)$, fix a Cartan involution $\theta$ and a $\theta$-commuting involutive automorphism $\sigma$ of $G$ such that the subgroup $H$ of $G$ of the elements fixed by $\sigma$ is locally isomorphic to $\SO(1,l-1)$. The quotient space $M=G/H$ is a $l$-dimensional Lorentzian symmetric space, the {\sl anti de Sitter space-time}.  We denote by $\sG$ and $\sH$ the Lie algebras of $G$ and $H$. We have the decomposition $\sG=\sH\oplus\sQ$ into the $\pm 1$-eigenspace  of the differential at $e$ of $\sigma$ that we denote again by $\sigma$.  We also consider $\sG=\sK\oplus\sP$, the Cartan decomposition induced by $\theta$; and $\sA$, a $\sigma$-stable maximally abelian subalgebra of $\sP$. A positive system of roots is chosen  and let $\sN$ be the corresponding nilpotent subalgebra (see Iwasawa decomposition, theorem \ref{ThoIwasawaVrai}).  Set  $\overline{\sN}=\theta(\sN)$, $\sR=\sA\oplus\sN$ and $\overline{\sR}=\sA\oplus\overline{\sN}$. Finally denote by $R=AN$ and $\overline{R}=A\overline{N}$ the corresponding analytic subgroups of $G$.  One then has

\begin{theorem}
The $l$-dimensional anti de Sitter space with $l\geq 3$, seen as the symmetric space $\SO(2,l-1)/\SO(1,l-1)$, becomes a causal solvable symmetric black hole, as defined above, when the closed orbits of the Iwasawa subgroup $R$ of $\SO(2,l-1)$ and its Cartan conjugated $\overline{ R }$ are said to be singular. There exists in particular a non empty event horizon. There are finitely many such closed orbits. 
\label{ThoLeBut}
 \end{theorem}

 This chapter intends to prove this theorem, and for the sake of completeness, we also analyze in some detail in section \ref{sec_AdSdeux} the two-dimensional case, for which the construction does not yield a black hole structure.

The black hole causal structure is thus completely determined by the action of a solvable group.  This observation gives prominence to potential embeddings of these spaces in the framework of noncommutative geometry, in defining noncommutative causal black holes (see also \cite{BDRS}) through the existence of universal deformation formulae for solvable groups actions which have been obtained in the context of WKB-quantization of symplectic symmetric spaces \cite{StrictSolvableSym,Biel-Massar-2}. These issues are investigated in chapter \ref{ChDefoBH} and in \cite{articleBVCS}.

\section{Causality, light cone and related topics} \label{SecCausal}

\subsection{Causality in anti de Sitter spaces}

We consider the $l$-dimensional anti de Sitter space (see appendix \ref{SecSymeStructAdS})
\begin{equation}    \label{eq:defAdS}
  AdS_l=\frac{ \SO(2,l-1) }{ \SO(1,l-1) }(\equiv u^2+t^2-x_1^2-\cdots-x_{l-1}^2=1).
 \end{equation}
According to proposition \ref{PropGHconn}, we can only consider the identity component of $\SO(2,l-1)$ and $\SO(1,l-1)$ instead of full groups. The metric that we put on $AdS_l$ is the one induced from the Killing form of $\SO(2,l-1)$ by formula \eqref{EqDefMetrHomo}. This metric has a Minkowskian signature, so that we have  natural notions of time-, space- and light-like vectors. From now we denote by $G$ and $H$ the groups $\SO(2,l-1)$ and $\SO(1,l-1)$. 

The connected group $\SO_0(2,l-1)$ admits an Iwasawa decomposition $ANK$ (see theorem \ref{ThoIwasawaVrai}). Let $A\bar N$ be the $\theta$-conjugate\footnote{Roughly speaking, it corresponds to different choices in the Iwasawa decomposition of $\SO(2,l-1)$.}group of $AN$ where $\theta$ is the Cartan involution of subsection \ref{SubSecCartandeuxN}. We will see that the actions of $AN$ and $A\bar N$ have closed and open orbits. The closed ones are denoted by $\hS_{AN}$ and $\hS_{A\bar N}$. The following definition is motivated all previously existing work about BTZ black hole.
\begin{definition}
The \defe{singularity}{Singularity} in $AdS_l$ is the set
\[
  \hS=\text{singularity}=\hS_{AN}\cup\hS_{A\bar N},
\]
so that a point is \defe{singular}{Singular point} when it belongs to a closed orbit of $AN$ or $A\bar N$. The \defe{black hole}{Black hole} is defined as
\[
  BH=\{ x\in AdS_{l} \text{ st } \forall \text{ time-like vector } k\in T_xAdS_l,\,  l^k_x\cap\mathcal{S}\neq\emptyset \}
\]
where $l^k_x$ is the (future directed) geodesic in the direction $k$ starting at $x$ (see equation \eqref{EqTousVecLumTy} and the discussion above).
\label{Singular}
\end{definition}

The aim of this chapter is to prove that the so-defined black hole is non trivial in the sense that the following inclusions are strict:
\begin{equation}		\label{EqhSssubBH}
 \hS\subset BH\subset AdS_l.
 \end{equation}

In order to get a full definition of the black hole and its structure, we need to define and characterise the notions of light ray and light cone. These notions are of course directly issued from physics of relativity.  
\begin{definition}
A \defe{light ray}{Light!ray} is a geodesic whose tangent vector is everywhere light-like.
\label{lightraycone}
 \end{definition}

\begin{remark}
The \defe{causal structure}{Causal structure} of a general pseudo-Riemannian manifold is the fact that two points are said to be \emph{causally connected} when there exists a light ray which passes by both points.
\end{remark}

A light ray trough $\mfo$ is given by a vector of $\sQ$ with vanishing norm. So let us study these vectors. Let $E_1=q_0+q_1$ and $k$, a general element of  $\SO(n)$ which reads $k= e^{K}$ with $K=a^{ij}(E_{ij}-E_{ji})$, $i,j\geq 3$ and $a^{ij}=-a^{ji}$.  If we pose $A_j=E_{1j}+E_{j1}$, we have $[K,E_1]=(2a)^{j3}A_j$ and $[K,A_k]=a^{jk}A_j$. Hence,
\[
\ad(K)^nE_1=\big((2a)^n\big)^{k3}A_k,
\]
and
\begin{equation} \label{eq:Adkeu} 
\begin{split}
\Ad(k)E_1=e^{\ad K}E_1&=E_1+\sum_{n\geq 1}\big( (2a)^n\big)^{k3}A_k\\
	      &=E_1+\sum_{n=0}^{\infty}\big(  (2a)^n \big)^{k3}A_k-\delta^{j3}A_j\\
		&=E_1-E_{31}-E_{13}+\big( e^{2a}\big)^{j3}A_j\\
	      &=q_0+\sum_{j=1}^{l-1}w_jq_j
\end{split}
\end{equation}
where $w_i=\big(  e^{2a} \big)^{i3}$. Under an explicit form, we have 
 \begin{equation} \label{eq:AdkE} 
   \Ad(k)E_1=
\begin{pmatrix}
0&1&w_1&w_2&\ldots\\
-1\\
w_1\\
w_2\\
\vdots
\end{pmatrix}
\end{equation}
The exponential $ e^{2a}$ being an element of $\SO(n)$, the parameters $w_i$ are restricted by the condition $\sum_{k}w_k=1$.  Remark moreover that \emph{every} matrix of $\SO(2)$ can be written under the form $e^{2a}$ for a good choice of $a\in\sod$. The light cone is therefore given by the set of vectors of the form $(1,w_i)$ with $\|w\|^2=1$. If we consider the metric $diag(+--\cdots)$ on $\sQ$ with respect to the basis $\{q_i\}$, we have
\[
  \|\Ad(k)E_1\|^2=0.
\]
This is coherent with the intuitive notion of light cone. It is on the one hand also true that \emph{every} light-like vector of $\sQ$ reads $\Ad(k)E_1$ for some $k\in\SO(n)$. On the other hand every nilpotent element of $\sQ$ is light-like because trace of nilpotent matrix is zero (using Engel's theorem). In definitive, we proved the following:

\begin{proposition}
When $E$ is any nilpotent element of $\sQ$, the set of light-like vectors of $\sQ$ is parametrized by $\lambda\Ad(k)E$ with $k\in\SO(n)$ and $\lambda\in\eR$.
\label{PropToutVectLumQ}
\end{proposition}

Let us point out the fact that only the first column of the \oge direction\fge{} $k\in \SO(n)$ has an importance in causality issues. So the word \oge directions\fge{} will often be used to refer to the vector $w$. It is not a particular feature of our particular matrix representation choie. Indeed the element $k$ only appears in the combination $\Ad(k)E$ which is a light-like vector in $\sQ$, i.e. $\Ad(k)E=tq'_0+\sum_i x_iq'_i$ with $t^2-\sum_i x_i^2=0$ for any orthonormal basis $\{q'_i\}$ of $\sQ$. As far as causality is concerned, a rescaling $\Ad(k)E$ to $\lambda\Ad(k)E$ has no importance, so one can choice $t=1$ and find back $\sum_i x_i^2=1$. We see that it is a natural feature that the light-like rays are parametrized by  unital vectors of $\eR^n$.

\begin{lemma}
Let $E$ be a nilpotent element in $\sQ$, and $\pi: G \rightarrow G/H$, the canonical projection. A light ray through $[g]\in AdS_l$ has the form
\begin{equation}
   l^k_{[g]}(t)=\pi\big( ge^{-t\Ad(k)E} \big)
\end{equation}
for a certain $k\in K_H=K\cap H=\SO(n)$.
 \label{lem:AdkEcone}
\end{lemma}

\begin{proof}
General theory of symmetric spaces (see \cite{kobayashi2}, pages 230--233, particularly theorem 3.2) proves that a light ray through $\mfo=[e]$ has the form
\[
  s(t)=\pi\big( e^{tX} \big).
\]
In our context, we have the additional request for the tangent vector to be light-like. Proposition \ref{PropToutVectLumQ} thus imposes $X$ to be of the form $\Ad(k)E$. That proves the claim for geodesics trough $\mfo$.

The fact that $d\tau_g$ is an nondegenerate isometry then extends the result to all points.

\end{proof}

\begin{corollary}
If $E$ is nilpotent in $\sQ$, then $\{\Ad(k)E\}_{k\in K_H}$ is the set of light-like vectors in $T_{[\mfo]}AdS_l\simeq\sQ$. Therefore
\begin{equation}
  \exp_{\mfo}( t\Ad(k)E )=\exp(t\Ad(k)E)\cdot\mfo.
\end{equation}
is the light cone of $\mfo$ in $AdS_l$.

\end{corollary}

In order to fix ideas, we will always use the element $E_1$ as choice of nilpotent element in $\sQ$ in order to parametrize light-cone.  Since $\SO(2,l-1)$ acts on $AdS_l$ by isometries, the \defe{light cone}{Light!cone} at $\pi(g)$ is given by a translation of the one at $\mfo$:
\begin{equation}	\label{eq_defcone}
  C^+_{\pi(g)}=g\cdot C_{\mfo}=\{  \pi\big( g e^{t\Ad(k)E_1}  \big)  \}_{\substack{t\in\eR^+\\ k\in K_H}}.
\end{equation} 
The product being taken at left while the quotient being taken at right, one can fear a problem of well definiteness in this expression. The following proposition shows that all is right.

\begin{proposition}
Definition \eqref{eq_defcone} is independent of the representative $g$ in the class $\pi(g)$. In other words,
\begin{equation}  \label{eq_statdefcone}
  \{ \Ad(hk)E_1 \}_{k\in K_H}=\{ \Ad(k)E_1 \}_{k\in K_H}
\end{equation} 
for all $h\in H$. 
\end{proposition}

\begin{proof}
The metric on $\sQ$ is the restriction of the Killing form of $\sG$ (notice that $\sQ$ has no own Killing form for the simple reason that it is not a Lie algebra). From $\Ad$-invariance, we have in particular
\[
  B\big(\Ad(h)X,\Ad(h)Y \big)=B(X,Y)
\]
for all $h\in \SO(1,l-1)$. The point is that reducibility makes $\Ad(h)X\in\sQ$ when $X\in\sQ$. The element $\Ad(hk)E_1$ in the left hand side of equation \eqref{eq_statdefcone} being zero-normed in $\sQ$, it reads $\Ad(k')E_1$ for some $k'\in K_H$. That proves the inclusion in one sense. For the second inclusion, we have to find a $k'\in K_H$ such that $\Ad(hk')E_1=\Ad(k)E_1$. Existence of such a $k'$ follows from the fact that $\Ad(h^{-1}k)E_1$ is a light-like vector of $\sQ$.
\end{proof}

\subsection{Time orientation}

A \defe{time orientation}{Time orientation} on $\sQ$ is the choice of a vector $T$ such that $\scal{T}{T}>0$. When such a choice is made, a vector $v$ is \defe{future directed}{Future directed vector} when $\scal{v}{T}>0$. In our case, the choice is the intuitive one: the vector $q_0$ defines the time orientation on $\sQ$ and $v=(v^0,v^1,v^2,v^3)$ is future directed if and only if $v^0>0$. So a light-like future directed vector is always --up to a positive multiple-- of the form $(1,\overline{v})$ with $\|\overline{v}\|=1$. For this reason, the set
\begin{equation}	\label{EqTousVecLumTy}
  \{t\Ad(k)E_1\}_{%
\begin{subarray}{l}
t>0\\k\in \SO(3)
\end{subarray}
}
\end{equation}
is exactly the set of light-like future-directed vectors of $\sQ$.

We are now able to define causality as follows.  A point $[g]\in AdS_l$ belongs to the \defe{interior region}{Interior region} if for every direction $k\in K_H$, the future light ray $l^k_{[g]}$ intersects the singularity within a \emph{finite} time.  In other words, it is interior when the whole light cone ends up in the singularity.  A point which is not interior is said to be \defe{exterior}{Exterior}. A particularly important set is the \defe{event horizon}{Event horizon}, or simply \emph{horizon}, defined as the boundary of the interior. When a space contains a non trivial causal structure (i.e. when there exists a non empty horizon), we say that the definition of singularity gives rise to a \defe{black hole}{Black hole}.  By extension, the term ``black hole'' often refers to the set of interior points.

\subsection{Some final remarks}

Remember that we decree closed orbits to be \emph{singular}. Now the fact for a point $\pi(g)\in AdS_l$ to be \emph{exterior} is that there exists an non empty set $\mO$ of $K_H$ such that $\forall k\in\mO$.
\[
  \pi\big( g e^{t\Ad(k)E_1}  \big)\cap\mS=\emptyset.
\]

The restriction of the Killing form to $\sQ$ reads
\begin{subequations}
\begin{align}
	B(q_0,q_0)&=\tr(q_0q_0)=-2,\\
	B(q_{i},q_{i})&=\tr(q_{i},q_{i})=2&\textrm{for $i\geq 1$}.
\end{align}
\end{subequations}
So the norm on $\sQ$ is $\| X \|=-\frac{ 1 }{2}B(X,X)$. The bi-invariance of the Killing form and the fact that the decomposition $\sG=\sQ\oplus\sH$ is reductive  imply $\| \Ad(h)X \|=\| X \|$, hence
\begin{equation}  \label{EqInclAdHSOq}
  \Ad(H)|_{\sQ}\subset\SO(\sQ).
\end{equation} 
A question is to know the kernel of this inclusion: which $h\in H$ fulfill $\Ad(h)q_i=q_i$ for all $i$ ? The equation $Aq_iA^{-1}=q_i$ can be simplified (from a computational point of view) using the relation $A^{-1}=\eta A^t\eta$ which defines $\SO(1,n)$. It is a somewhat long but easy computation to prove that $A=\pm\mtu$ are the only two solutions in $SO(1,n)$ to the system $A(q_i\eta)A^t=q_i\eta$.

One can go further than inclusion \eqref{EqInclAdHSOq} and prove the following
\begin{proposition}		
 Let $h\in H_0$ seen as a matrix acting on $\eR^{1,l-1}$ and let see $\Ad(h)$ as a matrix acting on $\sQ$. In this case we have $\Ad(h)_{ij}=h_{ij}$, in particular
\begin{equation}
   \Ad(H_0)=\SO_0(\sQ)
\end{equation} 
where the index zero denotes the identity component.
\label{PropSOADHequal}
\end{proposition}

\begin{proof}
We will prove that for each unital vector $X\in\sQ$, the element $\Ad(h)X$ is a general element of norm $1$ in $\sQ$ when $h$ runs over $H_0$. Explicit matrix computation will show by the way the equality  $\Ad(h)_{ij}=h_{ij}$. The general product to be computed is
\[ 
\Ad(h)X=
  \begin{pmatrix}
1	&	0\\
0	&
\begin{pmatrix}
&&\\
&h^{-1}\\
&&
\end{pmatrix}
\end{pmatrix}
\begin{pmatrix}
0&-w_0&w_1&\cdots\\
w_0\\
w_1\\
\vdots
\end{pmatrix}
  \begin{pmatrix}
1	&	0\\
0	&
\begin{pmatrix}
&&\\
&h\\
&&
\end{pmatrix}
\end{pmatrix}.
\]
But we know that the result is a matrix of $\sQ$, so it is sufficient to compute the first line. If we denote by $c_i$ the columns of $h$, we find
\[ 
  \Ad(h)X=\sum_{i=0}^{l-1}(w\cdot c_i)q_i
\]
where the dot denotes the inner product of $\eR^{1,l-1}$. Since $\{ c_i \}$ is a general orthonormal basis of $\eR^{1,l-1}$, the latter expression is a general vector of norm $1$ in $\sQ$.
\end{proof}

\section{Open and closed orbits}

\subsection{Openness of orbits in homogeneous spaces} \label{subsec:question}

\begin{proposition}
The orbits of $AN$ are submanifolds of $G/H$.
\label{pg:orbit_ssvar}
\end{proposition}

\begin{proof}
 Indeed proposition 4.4 in \cite{Helgason} (page 125) makes $R/(R\cap H)$ the orbit of $\pi(e)$ by $R$ and assure us that it is a submanifold of $G/H$. That proves the proposition for the orbit of $e$. 

For the other orbits, we consider the group $R_z=\AD(z^{-1})R$ which is also a Lie  subgroup of $G$. The space $R_z/(R_z\cap H)$ is isomorphic to the orbit of $\pi(e)$ under the action of $R_z$. Therefore $zR[z^{-1}]$ is a submanifold of $G/H$ and the very definition of a Lie group makes that  $R[z^{-1}]$ is a submanifold too.

\end{proof}

Let us start by computing the closed orbits of the actions of $AN$ and $A\bar{N}$ on $AdS_l$. In order to see if $[g]\in AdS_l$ belongs to a closed orbit of $AN$, we ``compare'' the space spanned by the basis $\{d\pi dL_g q_i\}$ of $T_{[g]}AdS_l$ and the space spanned by the fundamental vectors of the action. If these two spaces are equal, then $[g]$ belongs to an open orbit (because a submanifold is open if and only if it has same dimension as the main manifold). That idea is precisely contained in the following theorem which holds for any homogeneous space $M=G/H$.

\begin{theorem}
If $R$ is a subgroup of $G$ with Lie algebra $\sR$, then the orbit $R\cdot \mfo$ is open in $G/H$ if and only if the projection $\dpt{\pr}{\sR}{\sQ}$ parallel to $\sH$ is surjective.
\label{tho:pr_ouvert}
\end{theorem}

The projection is defined by $\pr(X)=X_{\sQ}$ if $X=X_{\sQ}+X_{\sH}$ is the decomposition of $X\in\sG$ with respect to the decomposition $\sG=\sH\oplus\sQ$. We need two lemmas before to prove the theorem.

\begin{lemma}
The orbit $R\cdot\mfo$ is open if and only if
\[
    \Span\{X^*_{\mfo}|X\in\mR\}=T_{\mfo}M
\]
where $X^*$ is the fundamental field defined by equation \eqref{EqDefChmpFonfOff}.
\label{lem:equiv_1}
\end{lemma}

\begin{proof}
From general theory of fundamental fields we know that
\[
\Span\{X^*_{\mfo}|X\in\sG\}=T_{\mfo}M.
\]
The game is now to prove that one can replace $\sG$ by $\sR$ if and only if $R\cdot \mfo$ is open.

\subdem{Necessary condition}
If $R\cdot\mfo$ is open, we have a neighbourhood of $\mfo$ which is contained in $R\cdot\mfo$. Then for any $X\in\sG$, and for a small enough $t$, the element $e^{-tX}\cdot\mfo$ belongs to $R\cdot\mfo$. Hence we have a path $r_X(t)$ in $R$ such that $e^{-tX}\cdot\mfo=r_X(t)\cdot\mfo$:
\[
      \Dsdd{e^{-tX}\cdot\mfo}{t}{0}=\Dsdd{r_X(t)\cdot\mfo}{t}{0}.
\]
Since $r_X(t)$ is a path in $R$, we can replace it by a $e^{-tY}$ with a $Y\in\mR$ in the derivative. For this $Y$, we have $X^*_{\mfo}=Y^*_{\mfo}$.

\subdem{Sufficient condition} We have $\dim(R\cdot\mfo)=\dim\Span\{ X^*_{\mfo}\tq X\in\sR \}=\dim T_{\mfo}M$,
so $R\cdot\mfo$ has the same dimension as $M$. The conclusion follows from the fact that a submanifold is open if and only if it has maximal dimension.

\end{proof}

\begin{lemma}
The canonical projection is surjective from $\sR$ to the tangent space to identity:
\begin{equation}\label{eq:equiv_2}
    \Span\{X^*_{\mfo}|X\in\mR\}=d\pi_e(\mR).
\end{equation}

\label{XsdpiR}

\end{lemma}

\begin{proof}
 Consider the following computation when $X\in\mR=T_eR$ is given by the path $X(t)=e^{tX}$:
\begin{equation}
  d\pi_e X=\Dsdd{[X(t)]}{t}{0}
	=\Dsdd{e^{tX}\mfo}{t}{0}
	=Y^*_{\mfo}
\end{equation}
with $Y=-X$. Reading these lines from left to right shows that $d\pi_e(\mR)\subseteq\{X^*_{\mfo}:X\in\mR\}$ while reading it from right to left shows the inverse inclusion.
\end{proof}

\begin{proof}[Proof of theorem \ref{tho:pr_ouvert}]
From lemma \ref{lem:equiv_1} and lemma \ref{XsdpiR}, the orbit $R\cdot\mfo$ is open if and only if $\dpt{d\pi_e}{\mR}{T_{\mfo}M}$ is surjective. On the one hand any $X\in\mR$ can uniquely be written as $X=X_{\sH}+X_{\sQ}$ with $X_{\sH}\in\sH$ and $X_{\sQ}\in\sQ$. On the other hand it is clear that $d\pi_e X_{\sH}=0$, thus $R\cdot\mfo$ is open if and only if $\dpt{d\pi_e}{\RM}{T_{\mfo}M}$ is surjective.

Now, recall that $d\pi_e$ is surjective from $\sG$, hence it is surjective from $\sQ$. The first conclusion is that if $\dpt{\pr}{\mR}{\sQ}$ is surjective, then $R\cdot\mfo$ is open. The inverse implication remains to be proved.

We know that openness $R\cdot\mfo$ implies that $\dpt{d\pi_e}{\RM}{T_{\mfo}M}$ is bijective (surjective because $R\cdot\mfo$ is open and injective because $\dpt{d\pi_e}{\sQ}{T_{\mfo}M}$ is injective by proposition \ref{IPropdpiisomMTM}). From all that, one concludes that $\RM=\sQ$. Indeed,  suppose that $X_{\sQ}\in\sQ$ and $X_{\sQ}\notin\RM$. Since $\dpt{d\pi_e}{\RM}{T_{\mfo}M}$ is surjective, there exists a $X_{\sQ}'\in\RM$ such that $d\pi_eX_{\sQ}'=d\pi_eX_{\sQ}'$. This is impossible because $d\pi_e$ is injective from the whole $\sQ$.

\end{proof}

\subsection{Open orbits in anti de Sitter spaces}

Now the strategy is to to check openness of the $R$-orbit of $[g]$ by checking openness of the $\AD(g^{-1})R$-orbit of $\mfo$ using the theorem \ref{tho:pr_ouvert}.

The problem is simplified by the following remark.  We know that matrices of $K$ and $H$ are given by
\begin{equation}	\label{eq:K_H_SO}
  K\leadsto \begin{pmatrix}
                \SO(2)&   \\
		      & \SO(n)
            \end{pmatrix},\quad
  H\leadsto \begin{pmatrix}
                    1 & \\
		     & \SOun
            \end{pmatrix},
\end{equation}
so we obviously have
\[
\bigcup_{s\in \SO(2)} \tau_{AN}([s]) =\bigcup_{\substack { s\in \SO(2)\\ h\in \SO(n)}}[ANsh] =\bigcup_{k\in K} [ANk] =[G].
\]
This is nothing else than the fact that the $AN$-orbits are $AN$-invariant.
So the $K$ part of $[g]=ank$ alone fixes the orbit which contains $[g]$ and we have at most one orbit for each element in $\SO(2)$. Computations using theorem \ref{tho:pr_ouvert} show that the $R$-orbits of $[\mu]$ with
\[
\mu=
\begin{pmatrix}
\cos\mu &\sin\mu\\
-\sin\mu&\cos\mu\\
&&\mtu
\end{pmatrix}
\]
is not open if and only if $\sin \mu=0$. We will see later that they are actually closed (page \pageref{PgTopoOrb}), so that the singularity is described as
\begin{equation}\label{Sing2}
\hS=[AN(\pm\mtu_{\SO(2)})]\bigcup[A \bar{N}(\pm\mtu_{\SO(2)})].
\end{equation}
 Because of $AN$-invariance of the $AN$-orbits, the equation of the $AN$-closed orbits can be expressed as
\begin{equation}
\sin \mu=0.
\end{equation}

Notice that there are some differences between the two choices of Iwasawa decompositions of subsections \ref{TabelPrem} and \ref{TableSeconde} in the determination of open and closed orbits. In the ``new'' Iwasawa decomposition (the one which is always used if not mentioned), up to matrices of $\sH$, a general matrix of $\sR$ is $jJ_1+mM+lL+kJ_2$. If we note $x=m+l$,
\begin{equation} \label{eq:geneR}
\sR\leadsto
\begin{pmatrix}
0&x&k&-x\\
-x\\
k\\
-x
\end{pmatrix}
\end{equation}
and it is obvious that the matrix $q_0$ can't be obtained by combinations of such matrices. So the $R$-orbit of $\mfo$ is not open.

If we use the ``old'' Iwasawa decomposition, the result is completely different. We have
\begin{align}
  q_{0}&=\pr\left( \frac{ N+M }{ 2 } \right),
&q_{1}&=\pr H_{2},
&q_{2}&=\pr\left( N-\frac{ N+M }{ 2 } \right),
\end{align}
and other elements of $\sQ$ are projections of the matrices $V_{i}$'s.  So we see that the map $\dpt{\pr}{\iR}{\sQ}$ is surjective and \label{pg:mfo_ouvert} the orbit $R\cdot\mfo$ is open.

\subsection{Two other characterizations of the singularity} 

In this short section, we first give a coordinatewise characterization of the singularity (which allows some brute force computations), and then we point out that the vector field $J_1^*$ has vanishing norm on the singularity (see also \ref{PropAdSDeuxJannule}). That should make the connection with the quotient construction of the original BTZ black hole.  Notice that we do not classify all vectors from which vanishing of the norm define a singularity. The point is that one can make our black hole ``causally inextensible'' by making a discrete quotient of $AdS_l$ along the integral curves of $J^*_1$.

\begin{proposition}
In term of the embedding of $AdS_l$ in $\eR^{2,l-1}$, the closed orbits of $AN \subset \SO(2,l-1)$ are located at $y-t = 0$.  Similarly, the closed orbits of $A \bar{N}$ correspond to $y+t=0$. In other words, the equation
\begin{equation} \label{tcarrycarr}
t^2-y^2=0
 \end{equation}
describes the singularity $\hS=\hS_{AN}\cup\hS_{A\bar{N}}$.
\label{Proptcarrycarr}
\end{proposition}

\begin{proof}
The different fundamental vector fields of the $AN$ action can be computed with the matricial relation $X^*_{[g]}=-Xg\cdot\mfo$. For example, in $AdS_3$,
\[
\begin{split}
   M^*_{[g]}&=
\begin{pmatrix}
0&-1&0&1\\
1&0&-1&0\\
0&-1&0&1\\
1&0&-1&0
\end{pmatrix}
\begin{pmatrix}
u\\t\\x\\y
\end{pmatrix}
=
\begin{pmatrix}
-t+y\\u-x\\-t+y\\u-x
\end{pmatrix}\\
&=(y-t)\partial_u+(u-x)\partial_t+(y-t)\partial_x+(u-x)\partial_y.
\end{split}
\]
Full results are
\begin{subequations}\label{Gen}
\begin{align}
J_1^*&=-y\partial_t-t\partial_y							\label{EqNormeJun}\\
J_2^*&=-x\partial_u-u\partial_x                                                      \label{eq:Jds}\\
M^*  &=(y-t)\partial_u+(u-x)\partial_t+(y-t)\partial_x+(u-x)\partial_y\\
L^*  &=(y-t)\partial_u+(u+x)\partial_t+(t-y)\partial_x+(u+x)\partial_y\\
W_i^*&=-x_i\partial_t-x_i\partial_y+(y-t)\partial_i\\
V_j^*&=-x_j\partial_u-x_j\partial_x+(x-u)\partial_j,
\label{eq:Vjs}
\end{align}
\end{subequations}
with $i,j=3,\ldots,l-1$.
First consider points satisfying $t-y=0$. It is clear that, at these points, the $l$ vectors $J_1^*$, $M^*$, $L^*$ and $W_i^*$ are linearly dependent. Thus, there are at most $l-1$ linearly independent vectors amongst the $2(l-1)$ vectors \eqref{Gen}. We conclude that a point satisfying $t-y=0$ belongs to a closed orbit of $AN$.

Now we show that a point with $t-y\neq 0$ belongs to an open orbit of $AN$. It is easy to see that $J_1^*$, $L^*$ and $M^*$ are three linearly independent vectors. The vectors $V_i^*$ gives us $l-3$ more. Then they span a $l$-dimensional space.

The same can be done with the closed orbits of $A\bar{N}$.  The result is that a point belongs to a closed orbit of $A\bar{N}$ if and only if $t+y=0$.
\end{proof}
This shows that in the three dimensional case, our black hole reduces to the previously existing one. 

The following corollary shows that a discrete quotient of $AdS_l$ along the orbits of $J_1^*$ gives a direct higher-dimensional generalization of the non-rotating BTZ black hole.
\begin{corollary}
The singularity coincides with the set of points in $AdS_l$ where $\| J_1^* \|^2 = 0$ for the metric induced from the ambient space $\eR^{2,l-1}$.
\label{CorJannsingul}
\end{corollary}

\begin{proof}
The expression \eqref{EqNormeJun} shows that the norm of $J_1^* $ is $y^2-t^2$ which vanishes on the singularity.
\end{proof}

 In the three-dimensional case, it was shown in \cite{BTZ_deux,BTZB_un} that the non-rotating BTZ black hole singularity is precisely given by equation \eqref{tcarrycarr}. Hence, the following is a particular case of theorem \ref{ThoLeBut}:

\begin{corollary}
 The non-rotating BTZ black hole is a causal symmetric solvable black hole.
\end{corollary}

\subsection{Orbits and topology}
\label{PgTopoOrb}

Let  $D^{\pm}=AN\SO(n)\SO(2)^{\pm}$ where $\SO(2)^{\pm}$ are the subgroups of $\SO(2)\subset \SO(2,n)$ with strictly positive (negative). We see $\SO(2)$ and $\SO(n)$ as subgroups of $\SO(2,n)$ in the way indicated by equation \eqref{eq:K_H_SO}. Notice that the parts $\SO(2)$ and $\SO(n)$ are commuting and that $\SO(n)\subset H$. The notation $-\mtu_{\SO(2)}$ refers to the element of $\SO(2,n)$ which the identity as $AN$-component and $-\mtu$ as $\SO(2)$-component.

A continuous path from $[D^+]$ to $[D^-]$ must pass trough an element of the form $[AN\mtu_{\SO(2)}]$. We saw that the $AN$-orbit of such an element is not open while the $AN$-orbit of an element of $[D^+]$ is open. So we deduce that an orbit passing trough $[D^+]$ does not intersects $[D^-]$.

The set $[D^+]$ is connected in $G/H$ and $D^+$  being open in $G$, the set $[D^+]=\pi(D^+)$ is also open in $G/H$ from the definition of the topology 
(see \cite{Helgason}, chapter II, paragraph 4 and particularly the theorem 4.2). Now, the orbits of $AN$ in $[D^+]$ (who are all open) furnish an open partition of $[D^+]$. Such a partition is impossible for an open connected set. We deduce that $[D^+]$ is only one orbit of $AN$ in $G/H$. The same can be done with $[D^-]$.

We are left with the sets $[AN]$ and $[AN(-\mtu_{\SO(2)})]$ whose union is closed because we just saw that the complement is open. Now we prove that these two sets are disjoint, in such a way that they have to be separately closed. Existence of an intersection point between $[AN]$ and $[AN(-\mtu_{\SO(2)})]$ would lead to the existence of a $h\in H$ such that $an\mtu_{\SO(2)}=(-\mtu_{\SO(2)})h$, or
\[ 
  h=(-\mtu_{\SO(2)})an,
\]
that is a non trivial $K$-component to $h$ in the decomposition $KAN$, but the only $K$-component in $H$ is $\SO(n)$. Hence such a $h$ does not exist and $R[\mtu]\cap R[-\mtu_{\SO(2)}]=\emptyset$.

The conclusion is that the Iwasawa group $AN$ has only four orbits :
\begin{align}
[D^+],&&[D^-],&&[AN\mtu_{\SO(2)}],&&[AN(-\mtu_{\SO(2)})].
\end{align}
The two first are open and the other two are closed. Remark\label{PgNoticeKpassung} that an element of $[K]$ does not belongs to a closed orbit of $AN$ or $A\bar N$.

\section{Horizons}

\subsection{Existence}

We are now able to prove that definition \ref{Singular} provides a non empty horizon as expressed by the condition \eqref{EqhSssubBH}.  First we  consider points of the form $\SO(2)\cdot\mfo$, which are parametrized by an angle $\mu$.  Up to the choice of this parametrization, a light-like geodesic trough $\mu$ is given by
 \begin{equation}
   K\cdot \mbox{e}^{-s\Ad(k)E_1}\cdot\mfo
\end{equation}
with $k\in \SO(l-1)$ and  $s\in\eR$. Using the isomorphism $[g]\mapsto g\cdot \mfo$ between $G/H$ and $AdS_l$, we find
\[
\begin{split}
  l^k_{[u]}(s)= \pi\big( u e^{t\Ad(k)E_1} \big)&=
\begin{pmatrix}
\cos\mu&\sin\mu\\
-\sin\mu&\cos\mu\\
&&1\\
&&&1\\
&&&&1\\
&&&&&\ddots
\end{pmatrix}
 e^{s\Ad(k)E_1}
\begin{pmatrix}
1\\0\\0\\0\\0\\\vdots
\end{pmatrix}
\\
&=
\begin{pmatrix}
u_{k}(s)\\t_{k}(s)\\x_{k}(s)\\y_{k}(s)\\z_{k}(s)\\\vdots
\end{pmatrix}
\end{split}
\]
According to proposition \ref{Proptcarrycarr}, this geodesic reaches the singularity if $t_{k}(s)^{2}-y_{k}(s)^{2}=0$ for a certain (positive) $s$. Since $\Ad(k)E_1$ is nilpotent, the computation of $ e^{s\Ad(k)E_1}$ is simple and we only need the first column because it only acts on the first basis vector. A short computation shows that
\begin{equation}  \label{EqGedCompo}
  l_{[\mu]}^{k}(s)=
\begin{pmatrix}
\cos\mu-s\sin\mu\\
-\sin\mu-s\cos\mu\\
sw_{1}\\
sw_{2}\\
\vdots
\end{pmatrix}.
\end{equation}

We conclude that the geodesic reaches $\hS_{AN}$ and $\hS_{A\bar{N}}$ for values $s_{AN}$ and $s_{A\bar{N}}$ of the affine parameter, given by
\begin{align}   \label{eq:tempssingul}
 s_{AN}&= \frac{\sin\mu}{\cos\mu - w_2}&s_{A\bar{N}}&= \frac{\sin\mu}{\cos\mu + w_2}
\end{align}
where $w_{2}$ is the second component of the first column of $k$, see equation \eqref{eq:AdkE}; in particular $-1\leq w_2 \leq 1$.

Since the part $\sin \mu =0$ is precisely  $\hS_{AN}$, we may restrict ourselves to the open connected domain of $AdS_l$ given by $\sin \mu > 0$. More precisely, $\sin\mu=0$ is the equation of $\hS_{AN}$ in the $ANK$ decomposition. In the same way, $\hS_{A\bar{N}}$ is given by $\sin\mu'=0$ in the $A\bar{N}K$ decomposition.  In order to escape the singularity, the point $[\mu]$ needs both $s_{AN}$ and $s_{A\bar{N}}$ to be strictly positive.  It is only possible to find directions (i.e. a parameter $w_2$) which respects this condition when $\cos \mu>0$.  So the point
\begin{equation}  \label{EqUnPtHoriz}
u\equiv \cos\mu=0
\end{equation}
is one point of the horizon. Theorem \ref{ThoLeBut} is now proved. 

The following proposition contains some physical intuition about the nature of the horizon.

\begin{proposition}
A light-like geodesic which escapes the singularity (i.e. which does not intersect $\hS$) and which passes trough a point of the horizon is contained in the horizon.
\end{proposition}

\begin{proof}
Let $x=[g]$ be a point of the horizon and $\pi(ge^{tAd(k)E_1})$, a light-like geodesic escaping the singularity. Near from $x$, there exists a point $y=[g']$ in the black hole. From definition of a black hole, for all $k\in \SO(3)$ and $t_{0}\in\eR^{+}$, points of the form  $\pi(g'e^{t_0Ad(k)E_1})$ also belong to the black hole. From continuity, in each neighbourhood of $\pi(ge^{t_0Ad(k)E_1})$, there is such a $\pi(g'e^{t_0Ad(k)E_1})$. This proves that $\pi(ge^{t_0Ad(k)E_1})$ belongs to the closure of the black hole. But it is not in the interior of the black hole because (by assumption) the given geodesic escapes the singularity, so every point of the form $\pi\big( g e^{t_0\Ad(k)E_1} \big)$ belongs to the horizon.
\end{proof}

Let us consider the point of the horizon that we know (the one given by \eqref{EqUnPtHoriz}), and see how can that point hope to escape the singularity.  Equations \eqref{eq:tempssingul} which give the time needed to fall into the singularity become
\begin{align}
  t_{AN}&=\frac{1}{w_{2}}&t_{A \bar{N}}&=-\frac{1}{w_{2}}.
\end{align}
So for every $w_{2}\neq 0$, this point reaches the singularity within a finite time. Taking the direction $w_{2}=0$ the point is able to reject his fall to infinity. This agrees to physical intuition which is that the horizon corresponds to points that fall into the singularity within an infinite time.

\subsection{Characterization}

Let $D[g]$ be the set of light-like directions (vectors in $\SO(n)$) for which the point $[g]$ falls into $\hS_{AN}$. Similarly, the set $\overline{D}[g]$ is the one of directions which fall into $\hS_{A \bar{N}}$. One can express $\overline{ D }$ in terms of $D$:
\[
\begin{split}
\overline{ D }[g]&=\{ k\in\SO(n)\tq\exists t\text{ for which }\pi\big( g e^{t\Ad(k)E_1} \big)\in\hS_{A\bar N} \}\\
		&=\{ k\in\SO(n)\tq\exists t\text{ for which }\pi\big( \theta(g)\theta( e^{tAd(k)E_1}) \big)\in\hS_{AN}\}\\
		&=\{ k\in\SO(n)\tq \pi(k)\in D\big( \theta[g] \big)\}\\
		&=\{ k\in\SO(n)\tq k\in\big( D(\theta[g]) \big)_{\theta}\},
\end{split}
\]
So
\begin{equation} \label{eq:DbarD}
\overline{D}[g]=(D\theta[g])_{\theta}
\end{equation}
where by definition, $k_{\theta}=Jk$ with $J$ being defined by $\theta=\Ad(J)$ ($\theta$ is the Cartan involution). It is easy to see that $\theta$ changes the sign of the spacial part of $k$, i.e. changes $w_i\to -w_i$.

 A main property of $k_{\theta}$ is
\[
	\theta(\Ad(k)E_1)=\Ad(k_{\theta})E_1.
\]
Since $k_{\theta}$ only appears in the expression $\Ad(k)E_1$, that property is actually a sufficient characterization of $k_{\theta}$ for our purpose. In particular, $k_{\theta\theta}\neq k$, but $\Ad(k_{\theta\theta})E_1=\Ad(k)E_1$.

How to express the condition $g\in\hH$ in terms of $D[g]$ ? The condition to belong to the black hole is $D[g]\cup \overline{D}[g]=\SO(n)$. If the complementary of $D[g]\cup \overline{D}[g]$ has an interior (i.e. if it contains an open subset), then by continuity the complementary $D[g']\cup \overline{D}[g']$ has also an interior for all $[g']$ near from $[g]$. In this case, $[g]$ cannot belong to the horizon. So a characterization of $\hH$ is the fact that the boundary of $D[g]$ and $\overline{D}[g]$ coincide. Equation \eqref{eq:DbarD} expresses this condition under the form
\begin{equation}
  \Fr D[g]=\Fr \big( D(\theta[g])\big)_{\theta},
 \end{equation}
from which one immediately deduces that $\hH$ is $\theta$-invariant.

We have an expression of $D[\mu]$ for $\mu\in \SO(2)$ by examining equations \eqref{eq:tempssingul}. The set $D[\mu]$ is the set of $w_2\in [-1,1]$ such that $\cos \mu+w_2>0$:
\begin{equation}
  D[\mu]=]-\cos \mu,1[.
\end{equation}
So in order for $\mu$ to belong to $\hH$, the point $[\mu]$ must satisfy
\[
\overline{D}[\mu]=D[\theta \mu]_{\theta}=]-1,-\cos \mu[.
\]
Consequently, if $\mu'$ is the $K$-component of $\theta \mu$ in the $ANK$ decomposition, we impose $]-\cos \mu',1[=D[\theta \mu]\stackrel{!}{=}]-\cos \mu',1[$\,, and we can describe the horizon by
\begin{equation} \label{eq:caractcous}
\cos \mu=-\cos \mu'
\end{equation}
where $\mu'$ is the $K$-component of $\mu$ in the $A\bar{N}K$ decomposition.

\section{Simple example on \texorpdfstring{$AdS_{2}$}{AdS2}}\label{sec_AdSdeux}

As is appendix \ref{SubsecTwoDimAdSAdGH}, we see $AdS_2$ as $\Ad(G)H$.  From definition \ref{Singular}, the singularity is the closed orbits of $AN$ and $A\bar{N}$ for the adjoint action on $AdS_2$, and the notion of fundamental field is
\begin{equation}
H^*_x=\Dsdd{\Ad(e^{-tH})x}{t}{0} =-[H,x].
\end{equation}

A basis of the Lie algebra $\sA\oplus\sN$ is given by $\{E,H\}$. So $x$ will belong to a closed orbit if and only if $E_x^*\wedge H^*_x=0$. If we put $x=x_HH+x_EE+x_FF$, the computation is
\[
E_x^*\wedge H^*_x=[E,x]\wedge[H,x]
                 =4x_Hx_F (E\wedge F)+2x_Ex_F (H\wedge E)-2x_F^2 (H\wedge F).
\]
It is zero if and only if $x_F=0$. The closed orbit of $A\bar{N}$ is given by the same computation with $H^*_x\wedge F^*_x$. The part of these orbits contained in $AdS_2$ is the one with norm $8$:
\[
B(x,x)=8(x_H^2+x_Ex_F)\stackrel{!}{=}8.
\]
In both cases, it imposes $x_H=\pm 1$, and the closed orbits in $AdS_2$ are given by
\begin{subequations}\label{EqQuatreLInesDeux}
\begin{align}		
\hS_{AN}&=\pm H+\lambda E\\
\hS_{A\bar N}&=\pm H+\lambda F,
\end{align}
\end{subequations}
with $\lambda\in\eR$. The singularity is then the union $\hS_{AN}\cup\hS_{A\bar N}$ of four lines in the hyperboloid.

\begin{proposition}
 The singularity can equivalently be described as
\begin{equation}\label{condHH}
\hS=\{x \in \Ad(G)H \tq \|H^*_x\|=0\},
\end{equation}
which has to be compared with corollary \ref{CorJannsingul}.
\label{HSing}
\label{PropAdSDeuxJannule}
\end{proposition}

\begin{proof}
The condition \eqref{condHH} on $x$ reads
\begin{equation}\label{eq:BHxHx}
B([H,x],[H,x])=0.
\end{equation}
The most general\footnote{It is actually \emph{more} than the most general element to be considered because our space is $\Ad(G)H$, which is only a part of $\sldr$.} element $x$ in $\sldr$ is $x=x_AH+x_NE+x_FF$. We have $[x,H]=-2x_NE+2x_FF$, so that the condition \eqref{eq:BHxHx} becomes $x_Nx_F=0$. The two possibilities are $x=x_AH+x_NE$ and $x=x_AH+x_FF$. The singularity in $\gsl(2,\eR)$ is composed of the planes $(H,F)$ and $(H,E)$. The intersection between the plane $(H,F)$ and the hyperboloid is given by the equation
\[
B(aH+bF,aH+bF)=8
\]
whose solutions are $a=\pm 1$. The same is also true for the plane $(H,E)$. So we find that the set \eqref{condHH} is exactly the four lines  \eqref{EqQuatreLInesDeux}.

\end{proof}

One can check that light cone of a given point of the hyperboloid is given by the two straight lines trough the point; so it automatically intersects the singularity. As conclusion, every point of $AdS_2$ belong to the black hole. For this reason we say that there is no black hole in the two dimensional case because the inclusions \eqref{EqhSssubBH} are in fact equalities.
 .

\section{Conclusions and perspectives}		\label{SecConcPerspAd}

Higher-dimensional generalizations of the BTZ construction have been studied in the physics' literature, by classifying the one-parameter subgroups of $\Iso(AdS_l)=\SO(2,l-1)$, see \cite{Figueroa,AdSBH,Madden,Banados:1997df,Aminneborg,HolstPeldan}.  Nevertheless, the approach we adopt here is conceptually different. We first reinterpret the non-rotating BTZ black hole solution using symmetric spaces techniques and present an alternative way to express its singularity.  We saw the latter as the union of the closed orbits of Iwasawa subgroups of the isometry group.  As shown, this construction extends straightforwardly to higher dimensional cases, allowing to build a non trivial black hole on anti de Sitter spaces of arbitrary dimension $l\geq 3$.  From this point of view, all anti de Sitter spaces of dimension $l\geq 3$ appear on an equal footing.

A natural question arising from this analysis is the following: \emph{given a semisimple symmetric space, when does the set of closed orbits of the Iwasawa subgroups of the isometry group, seen as singularity, define a non-trivial causal structure ?} We answered this question in the case of anti de Sitter spaces, using techniques allowing in principle for generalization to any semisimple symmetric space.

We also proved that performing a discrete quotient along the orbits of $J_1$ makes the resulting space causally inextensible (closed space-like curves appear in the singular part of the space), but we did not address  questions like: are there other vector fields defining singularities (in the three dimensional case, we know that the answer is positive) ? Can we identify a mass and an angular momentum from these hypothetic vectors ? Are \emph{all} BTZ black holes obtainable in this way in higher dimensions ?

\chapter{Deformation of anti de Sitter spaces}   \label{ChDefoBH}
\begin{abstract}

We are now going to apply deformation theory to the physical part of the $AdS_4$ black hole. The first idea was to deform the $AN$ of $\SO(2,l-1)$ and to deform $AdS_l$ by action of this group (see section \ref{SecDefAction} for an introduction to deformation by group action). We show in section \ref{SecUnifSOdn} that this procedure is possible.

Instead of that, we will only deform an open orbit of $AN$ in the four dimensional case\footnote{But the structure of algebra \eqref{EqTableSOIwa}, promises easy higher dimension generalisation.}. There are two reasons for that. First a physical domain of the black hole is contained in an open orbit of $AN$; and second it reveals possible to deform such a domain by action of a four dimensional group. Deforming a four dimensional space by a four dimensional group instead of a six dimensional one is a matter of ``no waste'' of dimensions.

The main lines of the construction are the following:
\begin{itemize}
\item We pick an open orbit $\mU$ of $AN$ in $AdS_4$, and we select a point $[u]\in\mU$.
\item We compute the stabilizer $S$ of $[u]$ in $AN$, in such a way that, as homogeneous space, $\mU=AN/S$. We consider the ``remaining group'' $R'$ of $R$ when one removes $S$ from $R$. 
\item We prove that $R'$ acts freely and transitively on $\mU$, so that $\mU$ is globally of group type (definition \ref{DefGlobGpType}).
 \item It turns out that $R'$ does not accept any symplectic structure; hence we will search for other groups acting transitively on $\mU$, and show that one and only one of them accepts a symplectic structure.
\item The latter group turns out to be a split extension of an Heisenberg group (see appendix \ref{SecExtHeiz}) for which we know a deformation.
\end{itemize}
\end{abstract}

\section{Group structure on the open orbits}
\label{SecGpStructOuvertOrb}

 \subsection{Global structure}

The following definition formalises the idea for a manifold to be ``like a group''. We will prove that the physical domain of our anti de Sitter black hole is of this type.
\begin{definition}
A $m$-dimensional homogeneous space $M$ is \defe{locally of (symplectic) group type}{Locally group type} if there exists a Lie (symplectic) subgroup $R$ of the group of automorphisms of $M$ which acts freely on one of its orbits in $M$. The homogeneous space is \defe{globally of group type}{Globally group type} if $R$ has only one orbit. In this case, for every choice of a base point $\mfo$ in $M$, the map $R\to M:g\mapsto g.\mfo$ is a diffeomorphism.
			\label{DefGlobGpType}
\end{definition}

Lie groups are themselves examples of symmetric spaces (globally) of group type.  In the symplectic situation, however, a symplectic symmetric Lie group must be abelian (\cite{ThzPierre}, page 12). We will see in what follows other non-abelian examples.  

In the context of our anti de Sitter black hole (in particular when one has causal issues in mind), it is not important to deform the whole space but it is sufficient\footnote{We will however point out in the perspectives (section \ref{SecConcPerspAd}) that a quantization of the whole space has a real interest.} to only deform one open orbit of $AN$. Indeed, if an observer begins his life somewhere in the physical space (hence in an open orbit of $AN$), he will never exit the orbit because one open orbit of $AN$ is bounded by closed orbit of $AN$ which are singular.

Let us recall that the solvable part of the Iwasawa decomposition of $\so(2,3)$ may be realized with a nilpotent part $\sN$ and an abelian one $\sA$  with elements
\begin{subequations}
\begin{align} 
\sA&=\{ J_1, J_2\}	&\sN&=\{W,V,M,L\}
\end{align}
\end{subequations}
and the commutator table 
\begin{subequations}	\label{EqTableRappelSO}
\begin{align}
[V,W]&=M &[V,L]&=2W\\
[ J_1,W]&=W       &[ J_2,V]&=V\\
[ J_1,L]&=L           &[ J_2,L]&=-L\\
[ J_1,M]&=M           &[ J_2,M]&=M,
\end{align}
\end{subequations}
where we know that that $W$, $J_{1}\in\sH$, and $J_{2}\in\sQ$. We pick the point
\[
  [u]=
\Bigg[
\begin{pmatrix}
0&1\\
-1&0\\
&&\mtu_{3\times 3}
\end{pmatrix}
\Bigg].
\]
This is an element of $K$ which, as already mentioned in page \pageref{PgNoticeKpassung}, therefore does not belong to a closed orbit of $AN$, neither to a one of $A\bar N$. Hence $[u]$ lies in the physical part of $AdS_4$. We denote by $\mU$ the $AN$-orbit of $[u]$.

Elements of the stabilizer of $[u]$ in $\SO(2,3)$ are elements $r$ such that $r\cdot[u]=[u]$, i.e. elements for which there exists  $h$ in $H$ such that $ru=uh$.
It is easy to check that $r_{1}=e^{aJ_{2}}$ and $r_{2}=e^{bV}$ are solutions by noticing that the action of $u^{-1}r_{i}u$ leaves unchanged the first basis vector\footnote{We give in a complement a more intrinsic way to prove that result, see page \pageref{subSecAltreintr}.}.

The stabilizer cannot contain more than two generators because an open orbit must be four dimensional. The stabilizer of $[u]$ in $G=\SO(2,3)$ is thus the group generated by $ e^{aJ_{2}}$ and $ e^{bV}$ plus eventually a discrete set making $S$ non connected.  The group $S$ is in fact connected because
\begin{equation}
   S=\{ r\in R\tq r\cdot[u]=[u] \}
	=\{ r\in R\tq \AD(u^{-1})r\in H \}.
\end{equation}
Since $R$ is an exponential group, we have $S=\exp\sS$ where
\[ 
  \sS=\{ X\in\sR\tq\Ad(u^{-1})X\in\sH \}=\Ad(u)\sH.
\]
The set $\sS$ being connected (because it is the image by a continuous map of the connected set $\sH$), $S$ is connected too. 

 The open orbit that we are studying is thus realised as the homogeneous space $\mU=AN/S$. An important result is the fact that what we obtain by simply removing $J_2$ and $V$ from the table \eqref{EqTableRappelSO} is still an algebra. The orbit $\mU$ is therefore isomorphic to the group $R'$ generated by the Lie algebra $\sR'=\{ J_1,W,M,L \}$. The table of $\sR'$ is 
\begin{subequations}
\begin{align}
  [J_{1},W]&=-W\\
[J_{1},L]&=L\\
[J_{1},M]&=M.
 \end{align}
\end{subequations}
From construction, $R'\cap S=\{ e \}$. Unfortunately, using the conditions \eqref{EqDefAlgSymple}, we find that in order to be compatible with the Lie algebra structure, the form $\omega$ of the algebra must satisfy $\omega(W,M)=\omega(W,L)=\omega(M,L)=0$, so that it is degenerate.  The action of $R'$ on $\mU$ enjoys however some remarkable properties.  
\begin{proposition}
The action
\begin{equation}
\begin{aligned}
 \tau\colon R'\times \mU&\to \mU \\ 
r[r_0u]&= [rr_0u] 
\end{aligned}
\end{equation}
is free and simply transitive. 
\label{PropURsimptra}
\end{proposition}

\begin{proof}

First, we prove that the action of $R'$ is transitive. As an algebra, $\sR$ is a split extension $\sR=\sS\oplus_{\ad}\sR'$.  Hence, as group, $R=S R'$, or equivalently $R=R'S$. That proves that the action is transitive.  

If the action is not simply transitive, there exists $x\in\mU$ and $r$, $r'\in R'$ such that $\tau_rx=\tau_{r'}x$. Since the action of $R'$ is transitive, we have a $r_1\in R'$ such that $x=r_1[u]$. In this case, the element $r_1^{-1}r^{-1}r'$ of $R'$ fixes $[u]$, but $R'\cap S=\{ e \}$. Then one deduces that $r'=rr_1$, so that $[rr_1u]=[r'r_1u]=[rr_1^2u]$. It follows that $r_1$ fixes $[u]$, and thus that $r_1=e$, so that $r=r'$.

For freeness remark that, in a neighbourhood of $e$, the neutral $e$ itself is the only element trivially acting on $[u]$.

\end{proof}

As corollary, the orbit $\mU$ is locally of group type $R'$.

\begin{proposition}		\label{PropmUsimpl}
The orbit $\mU$ is simply connected.
\end{proposition}

\begin{proof}
The fibration $S\to R\to\mU$ induces the long exact sequence of cohomology groups
\[ 
  H^{0}(\mU)\to H^{0}(R)\to H^{0}(S)\to H^{1}(\mU)\to H^{1}(R)\to\ldots
\]
The group $R$ being connected and simply connected, the sequence shows that $H^{0}(S)\simeq H^{1}(\mU)$, but we already mentioned that $S$ is connected, so $H^1(\mU)=0$.

\end{proof}

\begin{corollary}
The open orbit $\mU$ is globally of group type.
\label{CormUgloGppasSym}
\end{corollary}

\begin{proof}
It is immediately apparent from proof of proposition \ref{PropURsimptra}. Since
\[ 
  R'[u]=R'S[u]=R[u]=\mU,
\]
the group $R'$ acts freely on $\mU$ and has only one orbit.
\end{proof}

Remark  that it remains to be proved that $\mU$ is globally of \emph{symplectic} group type. For that, there should be a symplectic form on $R'$. Exploiting the fact that $\Span\{W,M,L \}$ is a three-dimensional abelian subalgebra of $\sR'$, it is easy to see that $\sR'$ does not accepts a symplectic form. Hence  corollary  \ref{CormUgloGppasSym} does not proves that $\mU$ is globally of symplectic group type. The lack of symplectic form on the algebra reflects on $\mU$ as manifold by the following lemma, and motivates the search for other four-dimensional groups than $R'$ acting transitively on $\mU$.

\begin{lemma}
The open orbit $\mU=R\cdot [u]$ does not admit any $R$-invariant symplectic form.
\end{lemma}

\begin{proof}
Let $\omega^{\mU}$ be such an invariant symplectic form and $\omega^{R'}$ be the pull-back of $\omega^{\mU}$ by the action: $\omega^{R'}=\tau^*\omega^{\mU}$.
We have $d\tau\circ dL_{r'}=dL_{r'}\circ d\tau$ because $\tau(r'X(t))=[r'X(t)u]=r'\tau(X(t))$, thus
\[ 
  L_{r'}^*\omega^{R'}=(\tau\circ L_{r'})^*\omega^{\mU}=(L_{r'}\circ\tau)^*\omega(\mU)=\omega^{R'},
\]
so that $\omega^{R'}$ is a $R'$-invariant symplectic form on $R'$. But we saw that such a form does not exist.
\end{proof}

\begin{proposition}
The $R$-homogeneous space $\mU$ admits a unique structure of globally group type symplectic symmetric space. The latter is isomorphic to $(\SUR_0,\omega,s)$ described in appendix \ref{SubsecDefSURme}.
\label{GT}
\end{proposition}
The next few pages are dedicated to prove this proposition and to give explicit algebra whose group gives the answer. We are searching for $4$-dimensional groups $\tilde R$ which
\begin{itemize} 
\item has a free and simply transitive action on $\mU$, i.e. $\tilde{R}[u]=R[u]$,
\item admits a symplectic structure,
\end{itemize}
and we want it to be unique.  As already mentioned, the algebra $\sR'$ fails to fulfil the symplectic condition. The algebra $\tilde{\sR}=\Span\{ A,B,C,D \}$ of a group which fulfils the first condition must at least act transitively on a small neighbourhood of $[u]$ and thus be of the form
\begin{subequations}  \label{EqGeneAlgabcd}
\begin{align} 
 A&=J_{1}+aJ_{2}+a'V\\
 B&=W+bJ_{2}+b'V\\
 C&=M+cJ_{2}+c'V\\
 D&=L+dJ_{2}+d'V.
\end{align}
\end{subequations}
Indeed, in a first attempt, we choose an algebra for which each of $A$, $B$, $C$ and $D$ contains a combination of $J_{1}$, $W$, $M$ and $L$. We consider the matrix of coefficients of $J_{1}$, $W$, $M$ and $L$ in $A$, $B$, $C$ and $D$. If the determinant of this matrix is zero, then one of the lines can be written as combination of the three others. In this case the action can even not be locally  transitive because the algebra only spans three directions actually acting ($J_{2}$ and $V$ have no importance here). So the determinant is non vanishing. In this case, the inverse of this matrix is a change of basis which puts $A$, $B$, $C$ and $D$ under the form  \eqref{EqGeneAlgabcd}.

The problem is now to fix the parameters $a,a',b,b',c,c',d,d'$ in such a way that the space $\Span\{ A,B,C,D \}$ becomes a Lie algebra (i.e. it closes under the Lie bracket) which admits a symplectic structure and whose group acts transitively on $\mU$. We will begin by proving that the surjectivity condition imposes $b=c=d=0$. Then the remaining conditions for $\tilde \sR$ to be an algebra are easy to solve by hand.

First, remark that $A$ acts on the algebra $\Span\{ B,C,D \}$ because $J_1$ does not appears in $[\sR,\sR]$. Hence we can write $\tilde{\sR}=\eR A\oplus_{\ad}\Span\{ B,C,D \}$ and, a subalgebra of a solvable exponential Lie algebra being a solvable exponential algebra, a general element of the group $\tilde{R}$ reads $\tilde r (\alpha,\beta,\gamma,\delta)= e^{\alpha A} e^{\beta B+\gamma C+\delta D}$. Our strategy will be to split this expression in order to get a product $S R'$ (which is equivalent to a product $R'S $). As Lie algebra, $\Span\{ B,C,D \}\subset\eR J_2\oplus_{\ad}\{ W,M,L,V \}$. Hence there exist functions $w$, $m$, $l$, $v$ and $x$ of $(\alpha,\beta,\gamma,\delta)$ such that
\begin{equation} \label{EqGeneRi}
 e^{\beta B+\gamma C+\delta D}= e^{xJ_2} e^{wW+mM+lL+vV}.
\end{equation}
We are now going to determine $l(\alpha,\beta,\gamma,\delta)$ and study the conditions needed in order for $l$ to be surjective on $\eR$. Since $J_2$ does not appear in any commutator, the Campbell-Baker-Hausdorff formula yields $x=\beta b+\gamma c+\delta d$. From the fact that $[J_2,L]=-L$, we see that the coefficient of $L$ in the left hand side of \eqref{EqGeneRi} is $-l(1- e^{-x})/x$. The $V$-component in the exponential can also get out without changing the coefficient of $L$. We are left with $\tilde r(\alpha,\beta,\gamma,\delta)= e^{\alpha A} e^{xJ_2} e^{yV} e^{w'W+m'M+lL}$ where $w'$ and $m'$ are complicated functions of $(\beta,\gamma,\delta)$ and $l$ is given by
\begin{equation}
l(\beta,\gamma,\delta)=\frac{ -\delta (\beta b+\gamma c+\delta d) }
{ 1- e^{-\beta b-\gamma c-\delta d} },
\end{equation}
which is only surjective when $b=c=d=0$. Taking the inverse, a general element of $\tilde R[u]$ reads $\big[ e^{ -wW-mM-lM} e^{j_1J_1}u \big]$, where the range of $l$ is not the whole $\eR$. Since the action of $R'$ is \emph{simply} transitive, $\tilde R$ is not surjective on $R[u]$ when $l(\alpha,\gamma,\delta)$ is not surjective on $\eR$.

When $b=c=d=0$, the conditions for \eqref{EqGeneAlgabcd} to be an algebra are easy to solve, leaving only two \emph{a priori} possible two-parameter families of algebras. 

{ \renewcommand{\theenumi}{\arabic{enumi}.}
\begin{enumerate}
\item 
The first one is the following:

\label{PgAlgUn}
\begin{subequations}
 \begin{align*}
A&=J_{1}+\frac{ 1 }{2}J_{2}+sV	&[A,B]&=B+sC\\
B&=W				&[A,C]&=\frac{ 3 }{2}C\\
C&=M				&[A,D]&=2sB+\frac{ 1 }{2}D\\
D&=L+rV				&[B,D]&=-rC.
\end{align*}
\end{subequations}
with $r\neq 0$. The general symplectic form on that algebra is given by
\begin{equation}   
\omega_{1}=\begin{pmatrix}
0	&-\alpha		&-\beta	&-\gamma\\
\alpha	&0			&0	&\frac{ 2\beta r }{ 3 }\\
\beta	&0			&0	&0\\
\gamma	&-\frac{ 2\beta r }{ 3 }	&0	&0
\end{pmatrix},
\end{equation} 

\[
\det\omega=\left( \frac{ 2\beta r }{ 3 } \right)^{2},
\]
Conditions: $\beta\neq 0$, $r\neq 0$. That algebra will be denoted by $\sR_{1}$. The analytic subgroup of $R$ whose Lie algebra is $\sR_1$ is denoted by $R_1$.  One can eliminate the two parameters in algebra $\sR_{1}$ by the isomorphism  
\begin{equation}		\label{EqIsomRUnrs}
\phi=
\begin{pmatrix}
1&0&0&0\\
0&1&0&4s\\
0&2sr&1/r&4s^{2}/r\\
0&0&0&1
\end{pmatrix}
\end{equation}
which fixes $s=0$ and $r=1$. The algebra $\sR_1$ is thus isomorphic to
\begin{equation}
\begin{aligned}[]
 [A,B]&=B			 &[A,C]&=\frac{ 3 }{ 2 }C\\
[A,D]&=\frac{ 1 }{ 2 }D		 &[B,D]&=-C.
\end{aligned}
\end{equation}

Comparing with equation \eqref{EqrhoBmudz}, one recognizes the one-dimensional extension of Heisenberg algebra with parameters $d=3/4$, $\mu=0$ and $\BX=\begin{pmatrix}
1&0\\0&1/2
\end{pmatrix}$. Hence $R_1$ is isomorphic to $\SUR_0$ and, by the way,  we have a product on that group (see appendix \ref{SecExtHeiz}).
\item 
The second algebra whose group acts simply transitively on $\mU$ is:
\begin{align*} 
A&=J_{1}+rJ_{2}+sV		&[A,B]&=B+sC\\
B&=W				&[A,C]&=(r+1)C\\
C&=M				&[A,D]&=2sB+(1-r)D.\\
D&=L
\end{align*}
There is no way to get a nondegenerate symplectic form on that algebra. 

\end{enumerate}
}		%

From proposition \ref{PropSymplestarEG}, the symplectic structure to be chosen on $\sR_1$ is $\delta C^*$ and lemma \ref{LemJumpCoadOrb} shows that we are able to quantize\footnote{by opposition to \emph{deform}: there are no symplectic condition in deformation.} $\sR_1$ with any symplectic form in the coadjoint orbit $\delta\big( C^*\circ\Ad(g) \big)$ with $g\in R_1$. The coadjoint adjoint action of $R_1$ on $\sR_1$ can be computed using the fact that $\sR_1$ splits into four parts; the non trivial results are
\begin{align*}
\Ad( e^{dD}A)&=A-\frac{ d }{2}D			&\Ad( e^{aA})B&= e^{a}B\\
\Ad( e^{cC})A&=A-\frac{ 3c }{2}C		&\Ad( e^{aA})C&= e^{3a/2}C\\
\Ad( e^{bB})D&=D-bC				&\Ad( e^{aA})D&= e^{a/2}D\\
\Ad( e^{bB})D&=D-bC.
\end{align*}
Direct computations show that
\begin{equation}
\begin{split}
	\Ad\big(  e^{aA} e^{bB} e^{cC} e^{dD} \big)(x_AA+x_BB&+x_CC+x_DD)\\
			&=x_AA+ e^{a}(x_B-x_Ab)B\\
			&\quad + e^{3a/2}\Big( x_C-\frac{ 3x_Ac }{2}-bx_D +\frac{ bdx_A }{2} \Big)C\\
			&\quad+ e^{a/2}\Big( x_D-\frac{ dx_A }{2} \Big)D,
\end{split}
\end{equation}
so that, with more compact notations, 
\begin{align}
\big( C^*\circ\Ad(g) \big)(X)=\big( x_C-\frac{ 3x_Ac-bx_D+\frac{ bdx_A }{2} }{2} \big) e^{3a/2},
\end{align}
and the symplectic forms that we are able to deform are given by $\delta\big( C^*\circ\Ad(g) \big)$. It provides a two-parameter familly of symplectic forms
\begin{equation}
\omega_1^g=
\begin{pmatrix}
0&0&\beta&\gamma\\
0&0&0&-2\beta/3\\
-\beta&0&0&0\\
-\gamma&2\beta/3&0&0
\end{pmatrix},
\end{equation}
\[ 
  \det\omega_1^g=\frac{ 4\beta^4 }{ 9 }.
\]

It turns out that the action of the group $R_1$ has good properties that are given in the following proposition.
\begin{proposition}   
The action of $R_1$ on $\mU$ is free and simply transitive.
\label{PropCRunXXX}
\end{proposition}

\begin{proof}

First remark that the algebra $\sR_1$ can be written as a split extension:
\[ 
  \sR_1=\eR A\oplus_{\ad} \eR D\oplus_{\ad}\Span\{ B,C \},
\]
hence a general element of $R_1$ reads
 \begin{equation}
r_1(a,b,c,d)= e^{aA} e^{dD} e^{bW} e^{cM}.
\end{equation}
The work is now to expand it by replacing $A$, $B$, $C$, $D$ by their values in function of $J_1$, $W$, $M$, $L$, $J_2$ and $V$, and then to try to put all elements of $\sS$ on the left. This is done by virtue of Campbell-Baker-Hausdorff formula. The fact that $\Span\{ W,M,N,V \}$ is nilpotent dramatically reduces the difficulty. We have
\[ 
  \ln( e^{drV} e^{dL} e^{wW} e^{mM})= drV+dL+(d^{2}r+w)W+(\frac{1}{ 6 }d^{2}r^{2}+m+\frac{ 1 }{2}drw)M.
\]
One can find $m$ and $w$ (functions of $d$) such that the right hand side reduces to $drV+dL$. Hence we have, for some auxiliary functions $w$ and $m$,
\[ 
   e^{drV+dL}= e^{drV} e^{dL} e^{w(d)W} e^{m(d)M}
\]
and a general element of $R_1$ reads
\begin{equation}
 e^{asV+\frac{ a }{2}J_2} e^{drV} e^{aJ_1} e^{dL} e^{\big( w(d)+b \big)W} e^{\big( m(d)+c \big)M}=s(a,d)r'(a,b,c,d)
\end{equation}
with $s\in S$ and $r'\in R'$ which defines a bijective map $r_1\mapsto r'$ from $R_1$ to $R'$. This proves the transitivity of the action of $R_1$.

 For freeness, just remark that in a neighbourhood of $e$, no element of $R_1$ (but $e$) leaves $[u]$ unchanged.
\end{proof}

The conclusion is that $R_1$ is the group $\tilde{R}$ that we were searching for and that it is unique (up to the two-parameter isomorphism \eqref{EqIsomRUnrs}) as symplectic subgroup of $AN$ acting transitively on $\mU$. It concludes the proof of proposition \ref{GT}.

\subsection{Alternative more intrinsic proofs}		\label{subSecAltreintr}

\begin{proposition}
Let $J\in Z(K)$ whose associated conjugation coincides with the Cartan involution: $\AD(J)=\theta$ and $u\in\SO(2,l-1)$ such that $u^2=J$ and $u\in e^{\sQ}\cap K$. Then the $AN$-orbit of $[u]$ is open.
\end{proposition}

\begin{proof}
Let us find the Lie algebra $\sS$ of the stabilizer $S$ of $[u]$. First, the Cartan involution $X\mapsto -X^{t}$ is implemented as $\AD(J)$ with
\[ 
  J=\begin{pmatrix}
-\mtu_{2\times 2}\\
&\mtu_{3\times 3}
\end{pmatrix}
\]
which satisfies $u^{2}=J$ and $\sigma(u)=u^{-1}$ because $u\in\sQ$. Now, $\AD(u^{-1})r\in H$ if and only if $\sigma\Big( \AD(u^{-1})r \Big)=\AD(u^{-1})r$. Using the fact that $\sigma$ is an involutive automorphism, we see that this condition is equivalent to
\begin{equation}
\theta\sigma r=r.
\end{equation}

On the one hand the Cartan involution $\theta$ restricts on $\sA$ to $\theta|_{\sA}=-\id$ because $\sA\subset\sP$; and on the other hand, $\sigma(\sA)=\sA$ because $J_1\in\sH$ and $J_2\in\sQ$. So $\sigma$ splits $\sA$ into two parts:  $\sA=\sA^{+}\oplus\sA^{-}$ with $\sA^+=\sA\cap\sH=\eR J_1$ and $\sA^-=\sA\cap\sQ=\eR J_2$. Let $\beta_{1},\beta_{2}\in\sA^*$ be the dual basis: $\beta_{i}(J_{j})=\delta_{ij}$.  We know that $W\in\sG_{\beta_{1}}$, $V\in\sG_{\beta_{2}}$, $L\in\sG_{\beta_{1}-\beta_{2}}$, and $M\in\sG_{\beta_{1}+\beta_{2}}$.  The set of simple roots is given by
\[ 
  \Delta=\{ \alpha=\beta_{1}-\beta_{2},\,\beta=\beta_{2} \},
\]
and the positive roots are
\[ 
 \Phi^{+}=\{ \alpha,\beta,\alpha+\beta,\alpha+2\beta \},
\]
in terms of whose, the space $\sN$ is given by
\begin{align*}
W&\in\sG_{\alpha+\beta}&V&\in\sG_{\beta}\\
L&\in\sG_{\alpha}&M&\in\sG_{\alpha+2\beta}.
\end{align*}
Since $(\sigma^*\beta)(h_1J_1+h_2J_2)=\beta_{2}(h_1J_1-h_2J_2)=-h_2$, we find
 we find
\begin{align*}
 \sigma^*\beta&=-\beta\\
\sigma^*\alpha&=\alpha+2\beta\\
\sigma^*(\alpha+\beta)&=\alpha+\beta\\
\sigma^*(\alpha+2\beta)&=\alpha.
\end{align*}
We are now able to identify the set $\sS=\{ X\in\sA\oplus\sN\tq \sigma\theta X=X \}$. Let us take $X\in\sR=\sA\oplus\sN$ and apply $\sigma\theta$:
\begin{equation}
\begin{split}
  X&=X^++X^-+X_{\alpha}+X_{\beta}+X_{\alpha+\beta}+X_{\alpha+2\beta},\\
\theta X&=-X^+-X^-+Y_{-\alpha}+Y_{-\beta}+Y_{-\beta}+Y_{-(\alpha+\beta)}+Y_{-(\alpha+2\beta)},\\
\sigma\theta X&=-X^++X^-+Z_{-(\alpha+2\beta)}+Z_{\beta}+Z_{-(\alpha+\beta)}+Z_{-\alpha}
\end{split}
\end{equation}
where $X_{\varphi}$, $Y_{\varphi}$ and $Z_{\varphi}$ denote elements of $\sG_{\varphi}$, and $X^{\pm}$ denote the component $\sA^{\pm}$ of $X$.

It is immediately apparent that $\sigma\theta X^-=X^-$, so that $X^-\in\sS$.  The only other component common to $X$ and $\sigma\theta X$ is in $\sG_{\beta}$, but   it is \emph{a priori} not clear that $X_{\beta}=Z_{\beta}$. We know however that $\sigma\theta V=\alpha V$ because $\sG_{\beta}$ has only one dimension. Using the fact that $\sigma$ and $\theta$ are commuting involutions, it is apparent that $\alpha=\pm 1$. Decomposing $V$ into $V=V_{\sH}+V_{\sQ}$, we have $\theta\sigma V=\theta(V_{\sH}+V_{\sQ})$ which has to be equal to $V$ or $-V$. Thus there are only two possibilities
\begin{align*}
  \theta V_{\sH}&=V_{\sH}&\text{or}	&&\theta V_{\sH}&=-V_{\sH}\\
\theta V_{\sQ}&=-V_{\sQ}&		&&\theta V_{\sQ}&=V_{\sQ}.
\end{align*}
If one compares the commutator table of $\SO(2,3)$ with the one of $\SO(2,2)$, one sees that $V$ is not present in $\SO(2,2)$. Since $\sH$ possesses every purely spatial rotation generators, the orthogonal complement $\sQ$ contains the time-time rotation as only rotations. Other components of $\sQ$ are boost. In particular, $V_{\sQ}$ is zero or a boost generator. In the latter case, $\theta V_{\sQ}=-V_{\sQ}$, and the conclusion is that $\sigma\theta V=V$. In the other case, $V_{\sQ}=0$ implies that $V\in\sH$ which is impossible because $[J_2,V]=0$ while $J_2\in\sQ$ and $[\sQ,\sH]\subset\sQ$.

The stabilizer $S$ is thus generated by $J_2$ and $\sG_{\beta}=\eR V$, i.e.
\begin{equation}
\sS=\Span\{ J_2,V \}.
\end{equation}
The stabilizer of $[u]$ being two-dimensional, the orbit of $[u]$ is four-dimensional and is then open in $AdS_4$.

\end{proof}

Notice that in contrast to the first way to find $\sS$, this time we have no eventually double covering problems. 

Let $\tilde R$ be a subgroup of $R$ whose Lie algebra is a complement of $\sS$ in $\sR$, i.e. $\tilde{\sR}\oplus\sS=\sR$. This group does not act transitively on $\mU$ if and only if the boundary of $\tilde R[u]$ is non empty. Let $x_0=\tau_{r'_0}[u]$ belong to that boundary with $r'_0\in R'$. On that point, the fundamental fields of $\tilde R$ are not surjective on the tangent space of $\mU$:
\[
\begin{split}
\ker\big[ \tilde{\sR}&\to T_{x_0}\mU \big] \neq\{ 0 \}\\
 Y&\mapsto Y^*_{x_0}.
\end{split}
\]
Let $Y\in\tilde R$ belongs to this kernel: $Y^*_{x_0}=0$. Since the linear map $\big( d\tau_{r_0'^{-1}} \big)_{x_0}$ is nondegenerate, $Y^*_{x_0}$ vanishes if and only if $\big( d\tau_{r_0'^{-1}} \big)_{x_0}(Y^*_{x_0})=0$, but
\begin{equation}
\begin{split}
 \big( d\tau_{r_0'^{-1}} \big)_{x_0}(Y^*_{x_0})=-\Dsdd{ r_0'^{-1} e^{tY}r_0'[u] }{t}{0}
		&=\big( \Ad(r_0'^{-1})Y \big)^*_{[u]}\\
		&=\pr_{\tilde{\sR}}\big( \Ad(r_0'^{-1})Y \big)
\end{split}
\end{equation}
because, on the point $[u]$, the action to take the fundamental field is nothing else than the projection parallel to $\sS$. Hence the group $\tilde R$ is not surjective if and only if
\[ 
  \big( \Ad(R')\sS \big)\cap\tilde{\sR}\neq\{ 0 \}.
\]
We are now going to determine $\Ad(R')\sS$. Let $X=X^{-}+X_{\beta}\in\sS$ and act with an element of $R'=\exp\big( \sA^{+}\oplus\sG_{\alpha}\oplus\sG_{\alpha+\beta}\oplus\sG_{\alpha+2\beta}\big)$:
\[ 
\begin{aligned}
   \Ad( e^{H^{+}+Y_{\alpha}+Y_{\alpha+\beta}+Y_{\alpha+2\beta}})(X^{-}+X_{\beta})
		&=X^{-}+X_{\beta}\\
		&+\underbrace{\big[ H^{+}+Y_{\alpha}+Y_{\alpha+\beta}+Y_{\alpha+2\beta},X^{-}+X_{\beta} \big]}_{N'}\\
		&\qquad+\frac{ 1 }{2}\big[  H^{+}+Y_{\alpha}+Y_{\alpha+\beta}+Y_{\alpha+2\beta},N'\big]\\
		&\qquad+\ldots
\end{aligned}
\]
The computation of $N'$ is as follows:
\begin{align*}
[H^+,X^-]&=0					&[H^{+},X_{\beta}]&=0\\
[Y_{\alpha},X^{-}]&=-\alpha(X^{-})Y_{\alpha}	&[Y_{\alpha},X_{\beta}]&=Z_{\alpha+\beta}\\
[Y_{\alpha+\beta},X^-]&=-(\alpha+\beta)(X^{-})Y_{\alpha+\beta}=0	&[Y_{\alpha+\beta},X_{\beta}]&=Z_{\alpha+2\beta}\\
[Y_{\alpha+2\beta},X^-]&=-(\alpha+2\beta)(X^{-})Z'_{\alpha+2\beta}	&[Y_{\alpha+2\beta},X_{\beta}]&=0,
\end{align*}
so $N'=-\alpha(X^{-})Y_{\alpha}+Z_{\alpha+\beta}-Z_{\alpha+2\beta}-(\alpha+2\beta)(X^{-})Z'_{\alpha+2\beta}$. Since $\beta(H^{+})=0$, the computation of $[H^{+},N']$, produces terms like $[H^{+},X_{\alpha+\beta}]=(\alpha+\beta)(H^{+})X_{\alpha+\beta}=\alpha(H^{+})X_{\alpha+\beta}$. Therefore, $[H^{+},N']=\alpha(H^{+})N'$ and
\begin{equation}  \label{EqElsInterditsSu}
\begin{split}
  \Ad( e^{H^{+}+Y})(X^{-}+X_{\beta})&=X+N'+\sum_{k\geq1}\frac{1}{ (k+1)! }\alpha(H^{+})^{k}N'\\
		&=X+\frac{  e^{\alpha(H^{+})}-1 }{ \alpha(H^{+}) } N'
\end{split}
\end{equation}
What we have proven is the following result.

\begin{proposition}
The group $\tilde R$ acts transitively on $\mU$ if and only if the Lie algebra $\tilde{\sR}$ does not contains elements of the form
\[
  X+\frac{  e^{\alpha(H^+}-1 }{ \alpha(H^+) }N'
\]
with $X\in\sS$ and $Y\in\sR'$; the element $N'$ being given by
\[
  N'=-\alpha(X^{-})Y_{\alpha}+Z_{\alpha+\beta}-Z_{\alpha+2\beta}-(\alpha+2\beta)(X^{-})Z'_{\alpha+2\beta}
\]
where $X=X^-+X_{\beta}$ is the decomposition of an element of $\sS$ and $Z_{\varphi}$ are elements of their respective root spaces $\sG_{\varphi}$.
\end{proposition}

One can distinguish three case: the first is $X=X^{-}\in\sA^{-}$, the second is $X=X_{\beta}\in\sG_{\beta}$ and the last one is $X=X^{-}+X_{\beta}$ ($X^{-}\neq  0\neq X_{\beta}$).

In the first case, formula \eqref{EqElsInterditsSu} forbids $\tilde{\sR}$ to contains elements of the form
\begin{equation}
  X^{-}+\sG_{\alpha}\oplus_{\alpha+2\beta}.
\end{equation}
The second case forbids elements of the form
\begin{equation}
  X_{\beta}+\sG_{\alpha+\beta}\oplus\sG_{\alpha+2\beta},
\end{equation}
and the third case forbids
\begin{equation}
  X^{-}+X_{\beta}\oplus\sG_{\alpha}\oplus\sG_{\alpha+\beta}\oplus\sG_{\alpha+2\beta}.
\end{equation}
We can extract constraints on the coefficients of algebra \eqref{EqGeneAlgabcd} from that analysis. The third interdiction makes that a linear combination of $J_2$ and $V$ in an element of $\tilde{\sR}$ can only occur in $A$, so that
\[ 
  bb'=cc'=dd'=0.
\]
The second interdiction says that $B$ and $C$ cannot contain $V$ alone, so $b'=c'=0$. Finally, the first condition imposes $c=d=0$ because $C$ and $D$ cannot contain $J_2$ alone. The remaining constraints for \eqref{EqGeneAlgabcd} to be an algebra are easy to solve by hand. The results are the same two algebras as previously found.

\subsection{Local group structure}

We saw in proposition \ref{PropCRunXXX} that the open orbit $\mU$ can be identified with the group generated by the algebra $\{ J_{1},W,M,L \}$. 

We want to find an algebra (whose group is) acting transitively on a neighbourhood of $[u]$ in the $AN$ orbit of $[u]$ and which admits a symplectic form. Let $A,B,C,D$ be a basis of this algebra. For local transitivity, each of them must contains at least one of $J_{1},W,M$ and $L$. As in the previous case, the most general algebra to be studied is
\begin{subequations}   \label{EqAlgGEnennsy}
\begin{align}
 A&=J_{1}+aJ_{2}+a'V\\
 B&=W+bJ_{2}+b'V\\
 C&=M+cJ_{2}+c'V\\
 D&=M+dJ_{2}+d'V.
\end{align}
\end{subequations}
Among such algebras, we will have to check surjectivity of the action and the possibility to endow with a symplectic form.

If we impose that $\Span\{ A,B,C,D \}$ is a subalgebra for the bracket inherited from $\mathfrak{so}(2,3)$, we find a lot of conditions on the coefficients $a$, $b$, $c$, $d$, $a'$, $b'$, $c'$ and $d'$. If, for example, we look at $[A,B]$, we find
\[
  [A,B]=W+a'M+ab'V-a'bV.
\]
In this combination, the coefficient of $W$ is $1$ and the one of $M$ is $a'$, so the only possibility is $[A,B]=B+a'C$. This leads to the following conditions (equating coefficients of $J_{2}$ and $V$):
\begin{subequations} \label{SubEqSystemeAlgTN}
\begin{align}
 b+a'c&=0\\
b'+a'c'&=ab'-a'b.
\end{align}
Proceeding in a similar way for the six different commutators, we find:

For $[A,C]$
\begin{align}
ac'-a'c&=(a+1)c'\\
(a+1)c&=0,
\end{align}
for $[A,D]$
\begin{align}
(1-a)d+2a'b&=0\\
ad'-a'd&=2a'b'+(1-a)d',
\end{align}
for $[B,C]$
\begin{align}
(b-c')c&=0\\
(b-c')c'&=bc'-b'c,
\end{align}
for $[B,D]$
\begin{align}
-d'c+2b'b-bd&=0\\
-d'c'-bd'+2(b')^{2}&=bd'-b'd,
\end{align}
for $[C,D]$
\begin{align}
cd'&=c'b'\\
cd&=c'b.
\end{align}
\end{subequations}
Solutions of these equations\footnote{from now until the determination of symplectic forms, all results are computed by Maxima~\cite{Maxima}.}, parametrized by reals $r$ and $s$ and the corresponding commutators are listed below.

 The next step is to determine which of these algebras admit a compatible symplectic structure.
 For this, we just have to consider a general skew-symmetric matrix
\[
  \omega=
\begin{pmatrix}
0&-\alpha&-\beta&-\gamma\\
\alpha&0&-\delta&-\sigma\\
\beta&\delta&0&-\epsilon\\
\gamma&\sigma&\epsilon&0
\end{pmatrix}
\]
 and, for each algebra, solve the four constrains. In the first algebra (see below), we find for example
\[
  \omega_{1}([A,B],C)+\omega_{1}([B,C],A)+\omega_{1}([C,A],B)=-\frac{ -5\omega_{1}(C,B) }{ 2 }\stackrel{!}{=}0,
\]
so that $\omega_{1}(C,B)=0$. Full results are listed below (the symplectic matrices are written in the basis $\{ A,B,C,D \}$).  We see in particular that most of the solutions reduce to the \emph{canonical algebra}, $\sR_c$ given by
\begin{align*}  
[a,b]&=b
&[a,c]&=2c
&[c,d]&=c.
\end{align*}

\let\ANCtheenumi\theenumi
 \renewcommand{\theenumi}{\arabic{enumi}.} 
\begin{enumerate} 
\item As first solution, we find of course the same algebra $\sR_1$ as the one of page \pageref{PgAlgUn}.
\item The second solution is also the same as the previously found one.
\item The third solution is
\begin{align*}
A&=J_{1}+J_{2}+sV		&[A,B]&=B+sC\\
B&=W-\frac{ r }{2}V		&[A,C]&=2C\\
C&=M				&[A,D]&=2sB\\
D&=L+rJ_{2}			&[B,D]&=-rB\\
 &				&[C,D]&=-rC,
\end{align*}
\begin{equation}
\omega=\begin{pmatrix}
0	&-\alpha		&-\beta	&-\gamma\\
\alpha	&0			&0	&\frac{ \beta rs-2\alpha r  }{ 2 }\\
\beta	&0			&0	&\frac{ \beta r }{2}\\
\gamma	&-\frac{ \beta rs-2\alpha r  }{ 2 }	&-\frac{ \beta r }{2}	&0
\end{pmatrix},
\end{equation}
\[
  \det\omega=\frac{ \beta^{4}r^{2}s^{2} }{ 4 }-\frac{ \alpha\beta^{3}r^{2}s }{ 2 }+\frac{ \alpha^{2}\beta^{2}r^{2} }{ 4 },
\]
Conditions: $r\neq 0$, $\beta\neq 0$ and $\alpha\neq\beta r$. The map $\phi_3\colon \sR_3\to \sR_c$
\[ 
  \phi_3=
\begin{pmatrix}
1	&	0	&	0	&	r\\
0	&	1/2s	&	0	&	1\\
0	&	s	&	1	&	s^2\\
0	&	0	&	0	&	r
\end{pmatrix}.
\]
($\det\phi_3=r/2s$) provides an isomorphism between $\sR_3$ and the canonical algebra.
\item The fourth solution is
\begin{align*}
A&=J_{1}+J_{2}+sV		&[A,B]&=B+sC\\
B&=W				&[A,C]&=2C\\
C&=M				&[A,D]&=2sB\\
D&=L+rJ_{2}+rsV			&[B,D]&=-rsC\\
 &				&[C,D]&=-rC,
 \end{align*}
\begin{equation}
\omega=\begin{pmatrix}
0	&-\alpha		&-\beta	&-\gamma\\
\alpha	&0			&0	&\frac{ \beta r s  }{ 2 }\\
\beta	&0			&0	&\frac{ \beta r }{2}\\
\gamma	&-\frac{ \beta r s  }{ 2 }	&-\frac{ \beta r }{2}	&0
\end{pmatrix},
\end{equation}
\[
\det\omega=\frac{ \beta^{4}r^{2}s^{2} }{ 4 }-\frac{ \alpha\beta^{3}r^{2}s }{ 2 }+\frac{ \alpha^{2}\beta^{2}r^{2} }{ 4 }.
\]
Conditions:  $r\neq 0$, $\beta\neq 0$ and $\alpha\neq\beta s$.  The map $\phi_4\colon \sR_4\to \sR_c$
\[
  \phi_{4}=
\begin{pmatrix}
2&0&0&-r\\
0&1&1/s&rs\\
0&1&0&2rs\\
3&0&0&-r
\end{pmatrix}
\]
($\det\phi_4=-r/s$)  provides an isomorphism between $\sR_{4}$ and the canonical algebra.
\item The fifth solution is
\begin{align*}
A&=J_{1}-J_{2}+sV		&[A,B]&=B+sC\\
B&=W-rsJ_{2}+rs^{2}V		&[A,D]&=2sB+2D\\
C&=M+rJ_{2}-rsV			&[B,D]&=2rs^{2}B+rs^{3}C+rsD\\
D&=L+rs^{2}J_{2}-rs^{3}V	&[C,D]&=-2rsB-rs^{2}C-rD,\\
\end{align*}
\begin{equation}
\omega=\begin{pmatrix}
0&-\alpha&-\beta&\frac{ \beta r s^{2}+2\alpha r s+2\epsilon }{ r }\\
\alpha&0&0&\epsilon s\\
\beta&0&0&-\epsilon\\
-\frac{ \beta r s^{2}+2\alpha r s+2\epsilon }{ r }&-\epsilon s&\epsilon&0
\end{pmatrix}
\end{equation}
\begin{equation}
\det\omega=\beta^{2}\epsilon^{2}s^{2}+2\alpha\beta\epsilon^{2}s+\alpha^{2}\epsilon^{2},
\end{equation}
Conditions:  $r\neq 0$, $\epsilon\neq 0$, $\alpha\neq -\beta s$.  The map $\phi_5\colon \sR_5\to \sR_c$
\begin{equation}
\phi_{5}=
\begin{pmatrix}
1&0&0&0\\
0&1&0&2s\\
0&0&0&-1\\
0&-rs&r&rs^{2}
\end{pmatrix}
\end{equation}
($\det\phi_{5}=r$) provides an isomorphism between $\sR_{5}$ and the canonical algebra.
\item The sixth solution is
\begin{align*}
A&=J_{1}+J_{2}		&[A,B]&=B\\
B&=W			&[A,C]&=2C\\
C&=M			&[C,D]&=-rC,\\
D&=L+rJ_{2}		
\end{align*}
\begin{equation}
\omega=\begin{pmatrix}
0	&-\alpha		&-\beta	&-\gamma\\
\alpha	&0			&0	&0\\
\beta	&0			&0	&\frac{ \beta r }{2}\\
\gamma	&0			&-\frac{ \beta r }{2}	&0
\end{pmatrix},
\end{equation}
\begin{equation}
\det\omega=\left( \frac{ \alpha\beta r }{ 2 } \right)^2,
\end{equation}
Conditions: $\alpha\neq 0$, $\beta\neq 0$ and $r\neq 0$. This algebra is isomorphic to the next one.
\item The seventh solution is
\begin{align*}
A&=J_{1}-J_{2}		&[A,B]&=B\\
B&=W			&[A,D]&=2D\\
C&=M+rJ_{2}		&[C,D]&=-rD,\\
D&=L		
\end{align*}
\begin{equation}
\omega=\begin{pmatrix}
0	&-\alpha		&-\beta	&-\gamma\\
\alpha	&0			&0	&0\\
\beta	&0			&0	&\frac{ \gamma r }{2}\\
\gamma	&0			&-\frac{ \gamma r }{2}	&0
\end{pmatrix},
\end{equation}
\begin{equation}
\det\omega=\pm\left( \frac{ \alpha\gamma r }{ 2 } \right),
\end{equation}
and the conditions are $\alpha\neq 0$, $\gamma\neq 0$, $r\neq 0$.  The map $\phi_7\colon \sR_7\to \sR_c$
\[
  \phi_{7}=
\begin{pmatrix}
1&0&0&0\\
0&1&0&0\\
0&0&0&1\\
0&0&r&0
\end{pmatrix}
\]
($\det\phi_{7}=-r$) provides an isomorphism between $\sR_{7}$ and the canonical algebra.
\end{enumerate}
\let\ANCtheenumi\theenumi
All these algebras are solvable of order two (the commutators of commutators vanish) --- but not nilpotent.

\section{Isospectral deformations of \texorpdfstring{$M$}{M}}

In this section, we present a modified version of the oscillatory integral product \eqref{PRODUCT} leading to a left-invariant associative algebra structure on the space of square integrable functions on $\SUR_0$. Why is it better that the initial product defined over smooth compactly supported functions ? The motivation of considering the square integrable functions is the fact that the spectral triple defined in non commutative geometry contains the space of square integrable spinors (see the book \cite{ConnesNCG}). The fact to stabilize the space of square integrable functions is then an indispensable step in order to put our results in the framework of spectral geometry.

\begin{theorem} 

 Let $u$ and $v$ be smooth compactly supported functions on $R_0$. Define the following three-point functions:
 \begin{equation}
\begin{split}
S:=& S_V\big(\cosh(a_1-a_2)x_0, \cosh(a_2-a_0)x_1, \cosh(a_0-a_1)x_2\big)\\
&-\bigoplus_{0,1,2}\sinh\big(2(a_0-a_1)\big)z_2;
\end{split}
\end{equation}
and 
\[ 
\begin{split}
A:= \Big[&\cosh\big(2(a_1-a_2)\big)\cosh(2(a_2- a_0))\cosh(2(a_0-a_1))\\
& \big[\cosh(a_1-a_2)\cosh(a_2- a_0)\cosh(a_0-a_1)\big]^{\dim R_0-2}\Big]^{\frac{1}{2}}.
\end{split}
\]
Then the formula                                                       
\begin{equation}\label{HILB}
u\star^{(2)}_\theta v:=\frac{1}{\theta^{\dim R_0}}
\int_{R_0\times R_0}Ae^{\frac{2i}{\theta}S}u\otimes v
\end{equation}
extends to $L^2(R_0)$ as a left-invariant associative Hilbert algebra structure. In particular, one has the strong closedness property:
\begin{equation*}
\int u\star^{(2)}_\theta v=\int uv.
\end{equation*}
\label{thmL2}
\end{theorem}
\begin{proof}
The oscillatory integral product \eqref{PRODUCT} may be obtained by intertwining the Weyl product on the Schwartz space $\swS$ (in the Darboux global coordinates \eqref{DARBOUX}) by the following integral operator \cite{Biel-Massar}:
\begin{equation*}
\tau:=F ^{-1}\circ(\phi_\theta^{-1})^\star\circ F,
\end{equation*}
$F $ being the partial Fourier transform with respect to        the central variable $z$:
\begin{equation*}
F (u)(a,x,\xi):=\int e^{-i\xi z}u(a,x,z){\rm d}z;
\end{equation*}
and $\phi_\theta$ the one parameter family of diffeomorphism(s):
\begin{equation*}
\phi_\theta(a,x,\xi)=(a,\frac{1}{\cosh(\frac{\theta}{2}\xi)}x,
\frac{1}{\theta}\sinh(\theta\xi)).
\end{equation*}
Set $\mJ:=|(\phi^{-1})^\star\mbox{Jac}_\phi|^{-\frac{1}{2}}$ and observe that for all $u\in C^\infty\cap L^2$, the function $\mJ\,(\phi^{-1})^\star u$ belongs to $L^2$.  Indeed, one has
\begin{equation*}
\int|\mJ\,(\phi^{-1})^\star u|^2=\int|\phi^\star \mJ|^2\,|\mbox{Jac}_\phi|\,|u|^2=\int|u|^2.
\end{equation*}
Therefore, a standard density argument yields the following isometry:
\begin{equation*}
T_\theta:L^2(R_0)\longrightarrow L^2(R_0):
u\mapsto F ^{-1}\circ m_\mJ\circ(\phi^{-1})^\star\circ F (u),
\end{equation*}
where  $m_{\mJ}$ denotes the multiplication by $\mJ$.  Observing that $T_\theta=F ^{-1}\circ m_{\mJ}\circ F \circ\tau$, one has $\star^{(2)}_\theta=F ^{-1}\circ m_{\mJ}\circ F (\star_\theta)$.  A straightforward computation (similar to the one in \cite{StrictSolvableSym})     then yields the announced formula.
\end{proof}

Let us point out two facts with respect to the above formulas:
\begin{enumerate}
\item The oscillating three-point kernel $A\exp \big (\frac{2i}{\theta}\,S\big)$ is symmetric under cyclic permutations.
\item The above oscillating integral formula gives rise to a strongly closed, symmetry invariant, formal star product on the symplectic symmetric space $(R_0,\omega,s)$.
\end{enumerate}

\begin{proposition}
The space $L^2(R_0)^{\infty}$ of smooth vectors in $L^2(R_0)$ of the left regular representation closes as a subalgebra of $(L^2(R_0),\star^{(2)}_\theta)$.
\end{proposition}

\begin{proof}
First, observe that the space of smooth vectors may be described  as the intersection of the spaces $\{V_n\}$ where  $V_{n+1}:=(V_n)_1$,   with $V_0:=L^2(R_0)$ and $(V_n)_1$ is defined as the space of elements $a$ of $V_n$ such that, for all $X\in\sR_0$, $X.a$ exists as an element of $V_n$ (we endow it with the projective limit topology).                             

Let thus $a,b\in V_1$. Then, $(X.a)\star b+a\star(X.b)$ belongs to $V_0$.  Observing that $D \subset V_1$ and approximating $a$ and $b$ by sequences $\{a_n\}$ and $\{b_n\}$ in $D $, one gets (by continuity of $\star$): $(X.a)\star b+a\star(X.b)=\lim(X.a_n\star b_n+a_n\star(X.b_n))= \lim X.(a_n\star b_n)=X.(a\star b)$. Hence $a\star b$ belongs to $V_1$.  One then proceeds by induction.
\end{proof}

\section{Spin structure and Dirac operator}	\label{SecDirADs}

Construction of the frame bundle is a straightforward adaptation of theorem 2.2 (chapter II) in \cite{AnnikFranc}, while connection issues are adapted from proposition 1.3 (chapter III).  According proposition \ref{PropGHconn}, notations $G$ and $H$ stand for the identity components of $\SO(1,l-1)$ and $\SO(2,l-1)$.

\subsection{Frame bundle and spin structure}

An element of the frame bundle is a map from $\sQ$ to $T(G/H)$ of the form\footnote{See \ref{SubSechoappahomsp} for notations.} $d\mu_g\circ A$ where $g\in G$ and $A\in \SO_0(\sQ)$. By proposition \ref{PropSOADHequal}, there exists a $h\in H$ for which $A=\Ad(h)$ for every $A\in\SO_0(\sQ)$ so we have
\[ 
 d\mu_g\circ A=d\pi\circ dL_g\circ\Ad(h)=d\pi\circ dL_g \circ dL_h\circ dR_h=d\pi\circ dL_g \circ dL_h=d\mu_{gh}
\]
hence in fact every element in the frame bundle reads $d\mu_g$ for some $g\in G$. We conclude that the fibre $B_{[g]}$ over $[g]$ is made of maps of the form $d\tau_k$ with $k\in[g]$. The action of $H$ on the frame bundle is given by 
\[ 
  (d\mu_g)\cdot h=d\mu_g\circ\Ad(h).
\]

\begin{proposition}
The map
\begin{equation}
\begin{aligned}
 \beta\colon G&\to B \\ 
g&\mapsto d\mu_g 
\end{aligned}
\end{equation}
is a principal bundle isomorphism between the frame bundle and the principal bundle
\begin{equation}		\label{PrincHGGH}
\xymatrix{%
   	G\ar[d]^{\pi}& H\ar@{~>}[l]		\\
   					G/H
 }
\end{equation}
where $\pi$ is the natural projection, the action of $H$ is the right one and the wavy line means ``acts on''.
\end{proposition}

\begin{proof}
Surjectivity of $\beta$ is clear. For injectivity, suppose $d\mu_g=d\mu_{g'}$. In order for the two target spaces to be equal, one needs $g'=gh$ for a $h\in H$. Now we have, for all $q_j\in\sQ$, 
\begin{equation}
  d\mu_gq_j=d\mu_{gh}q_j=d\pi dR_{h^{-1}}dL_gdL_hq_j
		=d\pi dL_g\big( \Ad(h)q_j \big),
\end{equation}
but $d\pi$ is an isomorphism from $\sQ_g$, so we deduce that $q_j=\Ad(h)q_j$. Since we are using the connected component of $\SO(\sQ)$, that implies that $h=e$, and thus that $g=g'$. The following proves that $\beta$ is a morphism:
\[ 
  \beta(gh)=d\pi dL_g dL_h=d\pi dL_g dL_h dR_{h^{-1}}=d\pi dL_g \Ad(h)=\beta(g)\cdot h.
\]
\end{proof}

The following lemma provides  a convenient way to express the tangent bundle over $G/H$ as an associated bundle to the principal bundle \eqref{PrincHGGH}. We denote by $G\times_{\rho} \sQ$ the quotient of $G\times\sQ$ by the equivalence relation $(g,X)\sim(gh,\Ad(h^{-1})X)$ for all $h\in H$.

\begin{lemma}  
The map 
\begin{equation} 
\begin{aligned}
 \beta\colon G\times_{\rho}\sQ&\to TM \\ 
[g,X]&\mapsto d\tau_gd\pi X 
\end{aligned}
\end{equation}
with $\rho(h)X=\Ad(h^{-1})X$ is diffeomorphic. 
\label{LemBazHGGH}
\end{lemma}

\begin{proof}

 In order to check that $\beta$ is well defined, first compute
\[ 
  \beta[gh,\Ad(h^{-1})X]=d\tau_{gh}d\pi\Ad(h^{-1})X=d\pi dL_{gh}\Ad(h^{-1})X,
\]
and then using the fact that $d\pi dR_h=d\pi$, the latter line reduces to $d\pi dL_gX=\beta(g,X)$. For injectivity, let $\beta[g,X]=\beta[g',X']$. In order for these two to be vectors on the same point, there must exists a $h\in H$ such that $g'=gh$. The equality becomes $d\pi dL_g dL_h X'=d\pi dL_gX$. Commuting $d\pi$ with $dL_g$ and using the fact that $d\tau_g$ is an isomorphism, we are left with the condition $d\pi dL_h X'=d\pi X$.

An element of $\sG/\sH$ is an equivalence class which contains exactly one element of $\sQ$. In the right hand side of the condition, this element is $X$ while the element of $\sQ$ in the class $d\pi dL_h X$ is $\Ad(h)X'$. Equating these two elements, we find the condition $X'=\Ad(h^{-1})X$, which proves that $[g,X]=[g',X']$ and concludes the proof of the injectivity of $\beta$.
\end{proof}

The following proposition will prove useful in order to identity the spin structure over $AdS_4$.

\begin{proposition}
If $G$ is a connected Lie group and if $Z$ is the center of $G$, then
\begin{enumerate}
\item $\Ad_G$ is an analytic homomorphism from $G$ to $\Int(G)$, with kernel $Z$,
\item the map $[g]\to\Ad_G(g)$ is an analytic isomorphism from $G/Z$ to $\Int(\lG)$ (the class $[g]$ is taken with respect to $Z$).
\end{enumerate}
\end{proposition}
On the one hand that proposition among with the fact that $Z\big( \SP(2,\eR) \big)$ proves that the quotient $\SP(2,\eR)/\eZ_2$ is isomorphic to $\Int\big(\gsp(2,\eR)\big)$. On the other hand one knows that $\SO_0(2,3)$ has no center, so that $\SO_0(2,3)\simeq\Int(\so(2,3))$. But the subsection \ref{SubSecIsosp} provides an isomorphism between $\so(2,3)$ and $\gsp(2,\eR)$. Thus we have
\begin{equation}
\SP(2,\eR)/\eZ_2\simeq \SO_0(2,3).
\end{equation}
We denote by $\varphi\colon \SP(2,\eR)\to \SO_0(2,3)$ the corresponding homomorphism with kernel $\eZ_2$. In particular the restriction $\varphi|_{\SL(2,\eC)}$ is a double covering of $\SO_0(1,3)$. But $\chi$ is the same kind of double covering, so universality of $\SL(2,\eR)$ on $\SO_0(1,3)$ provides an automorphism $f\colon \SL(2,\eC)\to \SL(2,\eC)$ such that $\varphi=\chi\circ f$. The spin structure to be considered on $AdS_4$ is
\[ 
\xymatrix{%
   \Spin(1,3) \ar@{~>}[r]		&\SP(2,\eR)  \ar[rd] \ar[rr]^{\displaystyle\varphi} && 	\SO_0(2,3)\ar[ld]&\SO_0(1,3)\ar@{~>}[l]\\
   &&	AdS_4 
}
\]
where the action of $\Spin(1,3)$ on $\SP(2,\eR)$ is given by $a\cdot s=af^{-1}(s)$ where we identified $\Spin(1,3)$ with $\SL(2,\eC)$ as subgroup of $\SP(2,\eR)$. One immediately has $\varphi(a\cdot s)=\varphi(a)\chi(s)$.

\subsection{Connection}

There are a lot of ways to express a vector field $X\colon G/H\to T(G/H)$. From the identification $T(G/H)=G\times_{\rho}\sQ$, one has $X\colon G/H\to G\times_{\rho}\sQ$. As section of an associated bundle, $X$ can be expressed by an equivariant function $\hat{X}\colon G\to \sQ$ such that $X_{[g]}=[g,\hat{X}(g)]$. The $H$-equivariance of $\hat X$ means that $\hat{X}(gh)=\Ad(h^{-1})\hat{X}(g)$.  Let $X\in\sG$ and consider the function
\begin{equation}		\label{EqDefhatAcol}
\begin{aligned}
 \hat A_X\colon G&\to \sQ \\ 
g&\mapsto \big( \Ad(g^{-1})X \big)_{\sQ} 
\end{aligned}
\end{equation} 
which is equivariant because the decomposition $\sG=\sH\oplus\sQ$ is reductive. The corresponding vector field is
\[ 
  A_X[g]=\big[ g,\big( \Ad(g^{-1})X \big)_{\sQ} \big];
\]
or
\[ 
  A_X[g]=d\tau_gd\pi\big( \Ad(g^{-1})X \big)_{\sQ}=d\pi dL_g\big( \Ad(g^{-1})X \big)
\]
because $d\pi X_{\sQ}=d\pi X$. It is easy to check that the form
\[ 
  \omega_g(X)=-\big( dL_{g^{-1}}X \big)_{\sH}
\]
is a connection form on the principal bundle \eqref{PrincHGGH}.  We are going to determine the associated covariant derivative of this connection on the tangent space, and prove that it is torsion free. The horizontal lift of $A_X[g]$ is 
\begin{equation}	\label{EqovlAprQhor}
  \overline{ A }_X(g)=dL_g\big( \Ad(g^{-1})X \big)_{\sQ}=\Dsdd{ g e^{t\pr_{\sQ}\Ad(g^{-1})X} }{t}{0}.
\end{equation}
The equivariant function associated with the covariant derivative of $A_Y$ in the direction of $A_X$ is given by $(\overline{ A }_X)_g\hat A_Y$. Using expressions  \eqref{EqDefhatAcol} and  \eqref{EqovlAprQhor} of $\hat A_Y(g)$ and $\overline{ A }_X(g)$, we have
\[ 
\begin{split}
  (\bar A_X)_g\hat A_Y	&=\Dsdd{ \hat A_Y\big( g e^{t\pr_{\sQ}\Ad(g^{-1})X} \big)_{\sQ} }{t}{0}\\
			&=\Dsdd{ \left( \Ad\big(  e^{-t\pr_{\sQ}\Ad(g^{-1})X}g^{-1}  \big)Y \right)_{\sQ} }{t}{0}\\
			&=\Big( \ad\big( -\pr_{\sQ}\Ad(g^{-1})X \big)\Ad(g^{-1})Y \Big)_{\sQ}\\
			&=-\left[ \big( \Ad(g^{-1})X \big)_{\sQ},\Ad(g^{-1})Y  \right]_{\sQ}.
\end{split}  
\]
This commutator is an expression of the form $[Z_{\sQ},Z'_{\sQ}+Z'_{\sH}]_{\sQ}$. Using reducibility we find
\begin{equation}
  (\overline{ A }_X)_g\hat A_Y=-\Big[ \big( \Ad(g^{-1})X \big)_{\sQ},\big( \Ad(g^{-1})Y \big)_{\sH} \Big].
\end{equation}
The commutator produces
\[ 
  (\overline{ A }_X)_g\hat A_Y-(\overline{ A }_Y)g\hat A_X=-\hat A_{[X,Y]}(g),
\]
which by construction the equivariant function associated with the vector field $\nabla_{A_X}A_Y-\nabla_{A_Y}A_X$; so on the one hand we have
\[ 
\begin{split}
(\nabla_{A_X}A_Y-\nabla_{A_Y}A_X)[g]=-d\tau_gd\pi\hat A_{[X,Y]}(g)
		&=-d\tau_gd\pi\big( \Ad(g^{-1})[X,Y] \big)_{\sQ}\\
		&=-d\pi dR_g[X,Y].
\end{split}  
\]
On the other hand,
\[ 
  [A_X,A_Y][g]=d\pi[dR_g X,dR_gY]=-d\pi dR_g[X,Y],
\]
which proves that the connection is torsion free.

We are now going to study the horizontal vector fields on $\SP(2,\eR)$ with this connection and the homomorphism $h^{-1}$ of equation \eqref{Eqdefhspsl}. We have to study for which elements $\Sigma_a\in \SP(2,\eR)$ the expression
\begin{equation}		\label{EqAtrouverdhemu}
  \omega_a(\Sigma_a)=\omega_{h^{-1}(a)}\big( (dh^{-1})_a\Sigma_a \big)=-\Big( dL_{h^{-1}(a)^{-1}}dh^{-1}\Sigma_a \Big)_{\sH}
\end{equation}
vanishes. Every such element can of course be written under the form $\Sigma_a=dL_a\psi X$ for some $X\in \so(2,3)$. So we are lead to consider the expression
\begin{equation} 		\label{Eqdhemuconide}
  (dh^{-1})_a(dL_a)_e\psi X.
\end{equation}
It is easy to deal with that expression in the case of $a=e$:
\[ 
  (dh^{-1})_e\psi(X)=\psi^{-1}\psi X=X.
\]
In particular, if $\Sigma\in\mT$, then $dh^{-1}\Sigma\in\sQ$ and when $\Sigma\in\mI$, we have $dh^{-1}\Sigma\in\sH$. This result propagates to other elements $a\in \SP(2,\eR)$ using the general result
\[ 
  df\circ dL_g=\big( dl_{f(g)} \big)\circ df
\]
which holds for any group homomorphism $f$. Using that property with $h^{-1}$ on the point $a\in \SP(2,\eR)$, we find $(dh^{-1})_a\circ(dL_a)_e=\big( dL_{h^{-1}(a)} \big)\circ (dh^{-1})_e$, and the expression \eqref{EqAtrouverdhemu} becomes
\[ 
  \omega_a(dL_a\psi X)=\Big( dL_{\big(h^{-1}(a)\big)^{-1}}dh^{-1} dL_a\psi X \Big)_{\sH}=X_{\sH}.
\]
It is zero if and only if $X\in\sQ$, so that the horizontal vectors on $a$ are exactly the ones of $dL_a\psi\sQ=\mT_a$.

\subsection{Dirac operator}

When $\hat{s}\colon \SP(2,\eR) \to \Lambda W$ is the equivariant function associated with a spinor, the Dirac operator reads
\begin{equation}		\label{EqDiracAdsquatre}
\widehat{Ds}(a)=g_{ij}\gamma^j\widehat{\nabla_{t_i}s}(a)
		=g_{ij}\gamma^j\overline{ t }_i(a)\hat{s}
		=g_{ij}\gamma^j\tilde t_i(a)\hat{s}
\end{equation}
where the metric $g$ is the usual four-dimensional Minkowskian metric and the matrices $\gamma$ are the associated $4\times 4$ Dirac matrices. The elements $\tilde t_i(a)=dL_at_i=dL_a\psi(q_i)$ span the natural basis of $\mT_a$, see appendix \ref{SubSecRedspT}. The matrices $\gamma^i$ are the usual $4\times 4$ Dirac matrices for the $4$-dimensional Minkowskian metric.

One can find a change of basis which express the Dirac operator in terms of vectors of $R_1$. For that, let $\{ X_i \}$ be a basis of $R_1$. We have 
\[ 
  X^*_i[u]=\Dsdd{ [ e^{-tX_i}u] }{t}{0}=-d\pi dR_uX_i
\]
that is necessarily decomposable by corollary \ref{ICordpiietwii}as combinations of vectors of the form $d\pi dL_uq_i$ because $[u]$ belongs to an open orbit of the action of $R_1$. That defines a matrix $B$ by
\[ 
  d\pi dL_uq_i=B_{ij}d\pi dR_uX_j,
\]
and then a vector $Y\in \sH$ by
\begin{equation}
q_i=B_{ij}\Ad(u^{-1})X_j+Y.
\end{equation}
Now we have $\tilde t_i(a)=dL_a\psi\big( \Ad(u^{-1})B_{ij}X_j+Y \big)$. We can go further using the fact that $\psi\Big( \Ad\big( h^{-1}(a) \big)X \Big)=\Ad(a)\psi(X)$ for every $a\in\SP(2,\eR)$ and $X\in\sG$. Defining the vectors $s_i=\Ad\big( h(^{-1}) \big)\psi X_i$ we find
\begin{equation}
\tilde t_i(a)=B_{ij}\tilde s_i(a)+\widetilde{\psi(Y)}(a).
\end{equation}

\section{Perspectives}

A main achievement of spectral non-commutative geometry is the ability of retrieving the original Riemannian manifold from the data of the spectral triple. Such a result does not exist in the case of $AdS$ because the latter is a non-compact \emph{pseudo}-Riemannian manifold. The main lines of such a reconstruction method can however be foreseen in the case of anti de Sitter space.
\begin{itemize}
\item Knowing the family of products $\star^{(2)}_{\theta}$, we know in particular the usual commutative product of functions. That should allow us to find back the manifold $AdS_4$.
\item It is possible to extract the data of the curvature of the manifold from the data of its Dirac operator as the non-differential part of its square. That part will of course appear to be constant and negative (because we know that we were starting from anti de Sitter).
\end{itemize}

We only quantized an open orbit of $AdS_4$ because it is a whole physical domain. Quantization of the full space could be very interesting because of a special effect of the noncommutative product: two functions with disjoint supports can have a non vanishing product. What about the physical significance of that property when one multiplies a function supported in the singularity by a function supported in the physical part ?

There is another reason to study the quantization of the full space. We will show in section \ref{SecUnifSOdn} that a deformation of the full space by action of the Iwasawa component of $\SO(2,l-1)$ is possible. That quantization has the advantage of deforming the space by the action of the group which is precisely defining the singularity. In other words the \emph{same} group can describe a singularity and a quantization. A work to be done is to try and recover the special causal structure from the data of the quantized manifold. That structure must be in some way contained in the spectral triple.

\chapter{Two notes for further developments}		\label{ChapNoteDev}
\begin{abstract}
This chapter contains two directions that were explored during my thesis and that were not finished for different reasons. 

In the first section, we state a result of Unterberger in \cite{UnterD} which provides a deformation of the complex half-plane, and we show how to translate it as a new noncommutative product on the group $ax+b$, i.e. the Iwasawa subgroup of $\SL(2,\eR)$. The technique of deformation by group action described in appendix \ref{SecDefAction} then induces a deformed product on the dual of its Lie algebra. We do not study the properties (symmetries, maximal functional space of convergence, symplectic condition to be a true quantization, \ldots) of this product, but we show that Unterberger's result assures the existence of at least one good functional space. Unfortunately the formula reveals not to be universal; we show the lack of universality on two examples of actions of the group $ax+b$ on $AdS_2$. The failure is due to divergences of the derivatives of the functions $z_i$ (see equation \eqref{EqDefziDefA}). 

This study is motivated by the fact that recent work (not published yet) of P. Bieliavsky provides an universal deformation of the $AN$ of $\SL(2,\eR)$. We are thus allowed to say that the latter new product is ``better'' than the one of Unterberger. We do not address the question to know the precise point that makes the lack of universality in Unterberger.

The second section is an application of the extension lemma (lemma \ref{EXT}). We show that all the ingredients needed to deform the $AN$ of $\SO(2,n)$ are present. The idea was to deform the $AdS$ black hole using the action of the so-deformed $AN$. That should provide an alternative way to deform $AdS$ to the one presented in chapter \ref{ChDefoBH}, and a quantization of $AdS_l$ using the same group as the group which defines a black hole. That method would use the deformation by group action machinery described in appendix \ref{SecDefAction}. The arising question is naturally to know if that quantization is in some sense equivalent to the one given in chapter \ref{ChDefoBH} or not. That question is not answered yet.

\end{abstract}

\section{Formula of Unterberger on \texorpdfstring{$\SL(2,\eR)$}{SL2R} }	\label{SecEplolUnter}

The following results come from \cite{UnterD} (from page 1219) and provide a deformation\footnote{or, at least, a new noncommutative product.} of the half-plane
\[
  D=\{ (\xi,\eta)\tq \eta>0 \}\subset\eR^2.
\]

Before to give the precise statement that will be used, we need some definitions. A first product is defined by (we will precise the functional space later):
\begin{equation}
(f\circ g)(\xi,\eta)=\sum_{\alpha,\beta}\frac{ (-1)^{\alpha} }{ \alpha !\beta! }(4i\pi)^{-\alpha-\beta}(\partial^{\alpha}_q\partial^{\beta}_p\tilde f)(0,0)(\partial^{\beta}_q\partial^{\alpha}_p\tilde g)(0,0)
\end{equation}
where
\[ 
  \tilde f(p,q)=f\Big( p+\xi\big(q+\sqrt{1+q^2}\big),\eta\big( q+\sqrt{1+q^2} \big) \Big),
\]
and the same for $\tilde g$. In particular,
\[
  (f\circ g)(0,1)=4\int f\big( \Psi(q_{1},p_{1}) \big)g\big( \Psi(q_{2},p_{2}) \big) e^{ -4i\pi(-q_{1}p_{2}+q_{2}p_{1})}\,dq_{1}\,dp_{1}\,dq_{2}\,dp_{2}.
\]
with
\[
  \Psi(p,q)=\big(p,q+\sqrt{1+q^{2}}\big)
\]

\begin{definition}
Let $r_1$, $r_2$ and $n$ be real numbers with $r_1\geq 0$. We denote by $\Sigma^{n}_{r_1,r_2}$ the space of functions $f\in  C^{\infty}(D)$ such that for all  $(j,k)\in\eN\times\eN$, there exists a $C>0$ such that
\begin{equation} 
\Big| \left( \frac{ \partial }{ \partial\xi } \right)^{j}\left( \eta\frac{ \partial }{ \partial\eta } \right)^{k}f(\xi,\eta)  \Big|
\leq C\eta^{r_1}(1+\eta)^{r_2}(1+| \xi |)^{n-j}.
\end{equation} 
\end{definition}

Now, theorem 8.2 in \cite{UnterD} states
\begin{theorem}
Let $f\in\Sigma ^n_{r_1,r_2}$ and $g\in\Sigma^{n'}_{r_1'nr_2'}$. For each $N\in\eN_0$, the function
\begin{equation}
h_N=f\circ g-\sum_{\alpha+\beta\leq N-1} \frac{ (-1)^{\alpha} }{ \alpha !\beta ! }(4i\pi)^{-\alpha-\beta}\sum_{j,k,j',k'}C_{\beta,\alpha}^{j,k}C_{\alpha,\beta}^{j',k'}(e_1^je_2^kf)(e_1^{j'}e_2^{k'}g)
\end{equation}
belongs to the space $\Sigma_{r_1+r_1',r_2+r'_2}^{n+n'-N}$ if constants $C_{\alpha,\beta}^{j,k}$ are defined by the requirement that
\[ 
  (\epsilon_2^{\beta}\epsilon_1^{\alpha}f)(\xi_0,\eta_0)=\sum_{j,k}C_{\alpha,\beta}^{j,k}(e_1^je_2^kf)(\xi_0,\eta_0)
\]
for every smooth function $f$ and $(\xi_0,\eta_0)\in D$ when $j+k\leq \alpha+\beta$ and $j\geq\alpha$ and $C_{\alpha,\beta}^{j,k}=0$ otherwise.  The operators $\epsilon_i$ are defined by $\epsilon_1=e_1=\partial_{\xi}$ and $\epsilon_2=2\big[ 1+\big( \frac{ \eta_0 }{ \eta } \big)^2 \big]^{-1}(\xi_0\partial_{\xi}+\eta_0\partial_{\eta})$.
				\label{ThoUnterSigmaStable}
\end{theorem}

For our purpose, the point is that there exists a product on $\Fun(D)$ and that theorem \ref{ThoUnterSigmaStable} provides a functional space stabilized by the product. We are now going to translate this result in terms of the Iwasawa subgroup $R=AN$ of $\SL(2,\eR)$ that is parametrized (see \eqref{EqParmalSL}) by
\[
  (a,l)=\begin{pmatrix}
 e^{a}&l e^{a}\\0& e^{-a}
\end{pmatrix}.
\]
 The map 
\begin{equation}
\begin{aligned}
 j\colon R&\to R' \\ 
(a,l)&\mapsto ( e^{2a},l e^{2a}) 
\end{aligned}
\end{equation}
provides an isomorphism between $R$ and the group
\[
  R'=\left\{ (\alpha,\beta)= \begin{pmatrix}
\alpha&\beta\\0&1
\end{pmatrix},\,\alpha>0 \right\}.
\]
The inverse of $j$ is  $j^{-1}(\alpha,\beta)=(\ln\alpha^{1/2},\beta\alpha^{-1})$.  The group $R'$ acts on $D$ by
\begin{equation}
(\alpha,\beta)\cdot(\xi,\eta)=(\xi+\beta\alpha^{-1}\eta,\alpha^{-1}\eta)
\end{equation}
which is a freely transitive action. For each choice of ``reference point'' $(\xi_0,\eta_0)\in D$ we  build an identification $i\colon D\to R'$ by the requirement $ i(\xi,\eta)\cdot(\xi_0,\eta_0)=(\xi,\eta)$, that is
\begin{equation}
i(\xi,\eta)=\left( \frac{ \eta_0 }{ \eta },\frac{ \xi-\xi_0 }{ \eta } \right).
\end{equation}

Now we can identify $D$ to $R$ by $k\colon D\to R$, $k=j^{-1}\circ i$. For the choice $(\xi_0,\eta_0)=(0,1)$, we find $k(\xi,\eta)=(k_a(\xi,\eta),k_l(\xi,\eta))$ where
\begin{align}\label{Eqkaklexp}
	k_a(\xi,\eta)&=-\frac{ 1 }{2}\ln\eta,
  &k_l(\xi,\eta)&=\xi
\end{align}
and the function $f$ on $R$ corresponds to the function $\tilde f=f\circ k$ on $D$.

The result of Unterberger is that the function $f$ ``can be quantized'' if
\begin{equation}  \label{eq_condeUnterD}
\left|  \left( \frac{ \partial }{ \partial\xi } \right)^j\left( \eta\frac{ \partial }{ \partial\eta } \right)^k  \tilde f(\xi,\eta)  \right|
\leq C\eta^{r_1}(1+\eta)^{r_2}(1+| \xi |)^{n-j}
\end{equation}
where $n$, $r_1$ and $r_2$ are real numbers and $r_1\geq0$. We want to see what condition has to be imposed on $f$  in order for $\tilde f$ to fulfil this condition. In other words, we want to express the operator
\[ 
  A_{ij}=\left( \frac{ \partial }{ \partial\xi } \right)^j\left( \eta\frac{ \partial }{ \partial\eta } \right)^k
\]
in terms of the coordinates on $R$. For that we compute $\partial_{\xi}\tilde f$ and $(\eta\partial_{\eta})\tilde f$ in terms of $\partial_lf$ and $\partial_af$.

Let us precise that, when we write expressions like $\eta\partial_{\eta}$, we mean for example
\[
  (\eta\partial_{\eta}\tilde f)(\xi,2)=2(\partial_{\eta}f)(\xi,2).
\]
For $\partial_{\xi}\tilde f$ we have:
\[ 
\begin{split}
  (\partial_{\xi}\tilde f)(\xi,\eta)&=(\partial_lf)\circ k(\xi,\eta)(\partial_{\xi}k_k)(\xi,\eta)\\
				&\quad+(\partial_af)\circ k(\xi,\eta)(\partial_{\xi}k_a)(\xi,\eta),
\end{split}  
\]
using the formula \eqref{Eqkaklexp}, we find $(\partial_{\xi}\tilde f)(\xi,\eta)=(\partial_lf)\circ k(\xi,\eta)$ and we conclude that
\begin{equation}
\partial_{\xi}\tilde f=(\partial_lf)\circ k.
\end{equation}
For $(\eta\partial_{\eta})\tilde f$, we find
\[ 
\begin{split}
(\eta\partial_{\eta})(f\circ k)(\xi,\eta)&=\eta\big( \partial_{\eta}(f\circ k) \big)(\xi,\eta)\\
		&=\eta(\partial_af)\circ k(\xi,\eta)(\partial_{\eta}k_a)(\xi,\eta)\\
		&\quad+\eta(\partial_lf)\circ k(\xi,\eta)(\partial_{\eta}k_l)(\xi,\eta)\\
		&=-\frac{ 1 }{2}(\partial_af)\circ k(\xi,\eta),
\end{split}  
\]
and we conclude that
\begin{equation}
  (\eta\partial_{\eta})(f\circ k)=-\frac{ 1 }{2}(\partial_af)\circ k.
\end{equation}
So the operator $A_{ij}$, expressed on $R$, reads
\begin{equation}
  A_{ij}(f\circ k)=\left( -\frac{ 1 }{2} \right)^j(\partial_a^j\partial_l^if)\circ k,
\end{equation}
and condition \eqref{eq_condeUnterD}, with $(\xi,\eta)=k^{-1}(a,l)=(l, e^{-2a})$ reads now
\begin{equation}  \label{eq_condUR}
 \left|  \frac{1}{ 2^k }(\partial_a^k\partial_l^jf)(a,l)   \right|\leq C e^{-2r_1a}(1+ e^{-2a})^{r_2}(1+| l |)^{n-j}
\end{equation}
with $r_1\geq 0$ and $r_2$, being any real number.
From now on this regularity condition will be referred as the \emph{Unterberger's condition}. That condition characterises a stable functional space for the Unterberger product on $R$.

We want now  to test the deformation of manifold by action of $R$. A somewhat deceiving result that will be shown is that Unterberger's deformation of $R$ is not an universal deformation in the sense that we will find some action of $R$ on manifold for which the action deformation does not provides a deformation of the manifold.

\subsection{Action on the dual of its Lie algebra}

The action if given by
\[ 
  (a,l)\cdot\xi=(y_H+2y_El)H^*+y_E e^{-2a}E^*
\] 
where $\xi=y_HH^*+y_EE^*$ is any point in $\sR^*$.  The question is to know if the product $(u\star_{\sR^*}v)$ makes sense when $u$ and $v$ are compactly supported smooth functions on $\sR^*$. In order to address this question, we have to check if for every $\xi$ in $\sR^*$, the function
\[ 
   (\alpha^{\xi}u)(a,l)=u\big( (a,l)^{-1}\cdot\xi \big)
		=u\big( (y_H-2y_Ee^{-2a}l),y_E e^{2a} \big)
\]
fulfils condition \eqref{eq_condUR}. So we consider 
\[ 
  f(a,l)=u\big( \underbrace{(y_H-2y_Ee^{-2a}l)}_{z_H(a,l)},\underbrace{y_Ee^{-2a}}_{z_E(a,l)} \big),
\]
and we compute
\[ 
\begin{split}
  (\partial_lf)(a,l)&=(\partial_Hu)(z_H,z_E)\partial_l(y_H-2y_Ee^{-2a}l)\\
		&\quad+(\partial_Eu)(z_H,z_E)\partial_l(y_Ee^{-2a})\\
		&=(\partial_Hu)(z_H,z_E)(-2y_Ee^{-2a}),
\end{split}  
\]
so
\begin{equation}
(\partial_l^jf)(a,l)=(\partial_H^ju)(z_H,z_E)(-2y_Ee^{-2a})^j.
\end{equation}
The combination $y_E e^{-2a}$ which goes out is precisely $z_E$ which remains in the derivative of $u$. But the derivative of $u$ has compact support. Hence, in fact, the coefficient $y_E e^{-2a}$ remains constrained in the domain where the derivative of $u$ does not vanishes. The point is that the coefficient which go out with derivatives is exactly made of $z_H$ and $z_E$.

So $\sR^*$ is as deformable as $D$. More precisely, a deformation of $\sR^*$ by action of $R$ is induced by the deformation of $D$ by Unterberger.

\subsection{First action on the two dimensional anti de Sitter space}

We see $AdS_2$ as in \ref{SubsecGpAdsDeux} and we consider the following action of $AN$ on $AdS_2$:
\[ 
  r\cdot\Ad(g)H=\Ad(gr^{-1})H.
\]
It is easy to see what does this action become in terms of the cylinder:
\[ 
  \big(  e^{y_AH} e^{y_NE} \big)\cdot \Ad\big(  e^{x_KT} e^{x_NE} \big)
	=\Ad\big(   e^{x_KT} e^{x_NE-y_NE} e^{y_AH}  \big)H
\]
where the adjoint action of $ e^{y_AH}$ on $H$ is of course trivial. Thus we have
\begin{equation}
(y_A,y_N)\cdot (x_K,x_N)=(x_K,x_N-y_N).
\end{equation}
Notice that only one dimension of $AN$ really acts. This action is thus not the most natural one, but is gives an interesting first toy model. Using the notations of coordinates \eqref{EqCylAdSDeux}, we consider $x=\phi(\theta,h)\in AdS_2$ and $u\in  C^{\infty}_c(Cyl)$, a compact supported function on $AdS_2$ and we compute
\[ 
(\alpha^xu)(a,l)=u\big( (a,l)^{-1}\cdot x \big)\\
		=u\big( (-a,-l  e^{2a})\cdot x \big)\\
		=u\big( \theta,h+l e^{2a} \big),
\]
so that if we pose $f(a,l)=u\big( \theta,h+l e^{2a} \big)$, we have
\[ 
  (\partial_lf)(a,l)= e^{2a}(\partial_2u)(\theta,l e^{2a}).
\]
When one makes $a\to\infty$ and $l\to 0$ in such a way that $l e^{2a}$ remains constant, the function $(\partial_lf)$ diverges in an exponential way with respect to $a$. It contradicts Unterberger's condition \eqref{eq_condUR}.
  
\subsection{Second action on the two dimensional anti de Sitter space}

Let us now study the more natural action
\begin{equation}
  r\cdot \Ad(g)H=\Ad(rg)H.
\end{equation}
It is in general very difficult to find, for given $y_A$, $y_N$, $x_K$ and $x_N$, the numbers (unique by construction) $z_K$ and $z_N$ such that 
\[ 
  \Ad( e^{y_AH} e^{y_NE} e^{x_KT} e^{x_NE})H=\Ad( e^{z_KT} e^{z_NE})H.
\]
In order to simplify the computations, we use the lemma \ref{LemUnPtParOrbite} which states that we only have to perform the computation for one $(x_K,x_N)$ in each orbit. We begin by $x_K=x_N=0$, i.e. the orbit of $H$ itself. First, computations show that 
\[ 
\begin{split}
\Ad( e^{z_KT} e^{z_NE})H=&\big( \cos(2z_K)-\sin(2z_K)z_N \big)H\\
			&-2\big( \cos(2z_K)z_N+\sin(2z_K) \big)E\\
			&+\Big( \big( \cos(2z_K)-1 \big)z_N+\sin(2z_K) \Big)T.
\end{split}  
\]
Next, 
\[ 
  \Ad( e^{aH} e^{lE})H=\begin{pmatrix}
1&-2 e^{2a}l\\
0& -1
\end{pmatrix}
=H-2 e^{2a}lE.
\]
Comparing with the general form, we find that
\begin{equation}
\Ad( e^{aH} e^{lE})H=\Ad( e^{l e^{2a}E})H,
\end{equation}
or $(a,l)\cdot (0,0)=(0,l e^{2a})$. What is important in our deformation problem is the function
\[ 
  (\alpha^Hu)(a,l)=u\big( (-a,-l e^{2a}))\cdot H \big)=u(0,-l).
\]
This function of course satisfies the Unterberger condition when $u$ has a compact support.

The second orbit that we study is the one of $V=\Ad( e^{\pi T/4})H=-2E+T$. One has
\[ 
\begin{split}
\Ad( e^{aH} e^{lE})V=-lH + e^{-2a}( e^{4a}l^2- e^{4a}-1)E+ e^{-2a}T.
\end{split}  
\]
If we pose $c=\cos(2z_K)$, $s=\sin(2z_K)$ and $b= e^{2a}$, we have to solve the system
\begin{subequations}
	\begin{align}
		c-sz_N=-l  \\
		-2cz_N-2s=\frac{1}{ b }(b^2l^2-b^2-1)\\
		(c-1)z_N+s=\frac{1}{ b }\\
		c^2+s^2=1
	\end{align}
\end{subequations}
with respect to $c$, $s$ and $z_N$. One can check that the following is a solution:
\begin{subequations}
\begin{align}
 c&=\frac{ b^2l^2-2b^2l+b^2-1 }{ b^2l^2-2b^2l+b^2+1 }\\
s&=\frac{ 2b(1-l) }{ b^2l^2-2b^2l+b^2+1 }\\
z_N&=\frac{ b^2(1-l^2)-1 }{ 2b }.
\end{align}
\end{subequations}
If we pose $\bar{z}_K(a,l)=z_K(-a,-l e^{2a})$ and $\bar{z}_N(a,l)=z_N(-a,-l e^{2a})$, we have
\begin{align}
\bar{z}_K(a,l)&=\frac{ 1 }{2}\arcsin\left(\frac{ 2 e^{-2a}(1+l e^{2a}) }{ l^2+2l e^{-2a}+ e^{-4a}+1 }\right)\\
2\bar{z}_N(a,l)&= e^{-2a}- e^{2a}(l^2+1).
\end{align}
The principle of deformation by action of group leads us to  deal with the function
\[ 
  f(a,l)=u(\bar{z}_K(a,l),\bar{z}_N(a,l)),
\]
which should satisfies Unterberger's condition when $u$ is compactly supported. Notice that $z_K$ is a compact variable, so that $u$ can be non vanishing for all values of $z_K$ without violate the compact support requirement.  The derivative of $f$ with respect to $a$ uses the chain rule, and it is apparent the higher order derivatives have to use the Leibnitz formula:
\[ 
\frac{ \partial f }{ \partial a }(a,l)=(\partial_1u)(\bar z_K,\bar{z}_N)\frac{ \partial \bar{z}_K }{ \partial a }(a,l)
		+(\partial_2u)(\bar{z}_K,\bar{z}_N)\frac{ \partial \bar{z}_N }{ \partial a }(a,l).
\]
In order to give an idea of what is going on, here is the first derivative of $\bar{z}_K$ with respect to $a$:
\[ 
  (\partial_a\bar{z}_K)(a,l)=\frac{ 2 e^{2a} }{  e^{4a}l^2+2 e^{2a}l+ e^{4a}+1 }.
\]
Let us look at the limit $a\to-\infty$ on the line $l= e^{-2a}$. If one performs multiple derivatives of $f(a,l)$ with respect to $a$, Leibnitz rules yields a lot of terms of the form
\begin{equation}		\label{EqTermGeederrau}
  (\partial_1^p\partial_2^qu)\big(\bar{z}_K(a,l),\bar{z}_N(a,l)\big)(\partial_a^i\bar{z}_K)(a,l)^j(\partial_a^k\bar{z}_N)(a,l)^m.
\end{equation}
On the line $l= e^{-2a}$, the numerator of $(\partial_a^i\bar{z}_K)(a,l)$ is $( e^{4a}+4)^{2i}$ while the numerator is a sum and product of monomials of the form $( e^{4a}+N)$ with $N > 0$. At the limit, this factor in \eqref{EqTermGeederrau} goes to a finite number. The factor $(\partial_a^k\bar{z}_N)(a,l)$ is very different because  
\[ 
  (\partial_a^k\bar{z}_N)(a,l)=(-1)^k2^{k-1} e^{-2a}-2^{k-1} e^{2a}(l^2+1).
\]
which becomes
\[ 
  2^{k-1} e^{-2a}\big( (-1)^k-1 \big)-2^{k-1} e^{2a}
\]
on $l= e^{-2a}$. It goes to zero when $a\to-\infty$ and $k$ is even, but is goes to $-\infty$ at the same limit when $k$ is odd. The highest divergence in all the terms of type \eqref{EqTermGeederrau} in $(\partial^n_af)$ is expected for maximal $m$, so when $i=j=0$. This is a divergence as
\[ 
  x\mapsto e^{2(n-1)x}.
\]
Notice that this divergence increases when the order of derivative increases. Hence it contradicts Unterberger's condition which works with parameters $r_1$ and $r_2$ who are \emph{constant} with respect to the order of the derivative.
\section{Deformation of \texorpdfstring{$\SOdn$}{SO2n}}   \label{SecUnifSOdn}

\subsection{Decomposition as split extension}

We try to decompose $\sA\oplus\sN$ as symplectic sum in order to use the extension lemma. In that purpose, let us consider the change of basis \eqref{EqChmHJ} in $\sA$: $H_1=J_1-J_2$ and $H_2=J_1+J_2$.

If $H_1\in\mathfrak{s}_2$, then $L,V_i,W_i\in\mathfrak{s}_2$ because $\mathfrak{s}_1$ must act on $\mathfrak{s}_2$. Hence $M\in\mathfrak{s}_2$ and $H_2$ remains alone in $\mathfrak{s}_1$. That proves that $H_1\in\mathfrak{s}_1$. If we suppose that $H_2\in\mathfrak{s}_2$, we find
\begin{equation}  \label{eq_HLss}
 \begin{split}
  \mfs_1&=\{ H_1,L  \}\\
\mfs_2&=\{ H_2,V_i,W_j,M  \}.
\end{split}   
\end{equation}
The case $H_1,H_2\in\mathfrak{s}_1$ leads to
\begin{equation}   \label{Eq_HHVaMLW}
 \begin{split}
	\mfs_1&=\{ H_1,H_2,\overbrace{V_a,\ldots V_b}^{\text{even}}  \}\\
	\mfs_2&=\{ M,L,W_i,V_{\text{others}}  \}.
\end{split}   
\end{equation}
The symplectic condition excludes the second decomposition. Indeed for each $s$ such that $[H_1,s]=\alpha s$ (i.e. $s=V_i,W_j,L$), we have
\[ 
  \Omega_2\big(  e^{ad H_1}M, e^{\ad H_1}s \big)= e^{\alpha}\Omega_2(M,s)\stackrel{!}{=}\Omega_2(M,s).
\]
Hence $\Omega_2(M,s)=0$. This proves that the decomposition \eqref{Eq_HHVaMLW} imposes the symplectic form $\Omega_2$ to be degenerate. We are left with decomposition \eqref{eq_HLss}.

Root space decomposition of $SU(1,n)$ can be found on pages 314--315 of \cite{Knapp}: it has $\dim\sA=1$, $\dim\sG_{2f}=1$ and $\dim\sG_f=2(n-1)$. In $\mathfrak{s}_2$, we have $V_i\in\sG_1$, $W_j\in\sG_1$, $M\in\sG_2$, and when we look at $AdS_l=\SO(2,l-1)/\SO(1,l-1)$, we have $l-3$ matrices $V_i$ and $W_j$. Therefore $\mathfrak{s}_2$ is nothing else than the $\sA\oplus\sN$ of $\mathfrak{su}(1,l-2)$ (recall $l\geq3$). The analysis shows that $\mathfrak{s}_1$ is the $\sA\oplus\sN$ of $\mathfrak{su}(1,1)$.

\subsection{Conclusion and perspectives}

For our $AdS_l$ black hole, the algebra of the group which defines the singularity is the split extension 
\[ 
  (\sA\oplus\sN)_{\so(2,l-1)}=(\sA\oplus\sN)_{\mathfrak{su}(1,1)}\oplus_{\ad}(\sA\oplus\sN)_{\mathfrak{su}(1,l-2)}.
\]
A deformation of the corresponding groups is given in the article \cite{Biel-Massar}, and the extension lemma \ref{EXT} yields now an oscillatory integral universal deformation formula for proper actions of the Iwasawa subgroup of $\SO(2,l-1)$. That remark provides an alternative way to deform the black hole to the one presented in section \ref{SecGpStructOuvertOrb}.

The availability of a quantization of $AdS_l$ by action of $AN$ is an opportunity to embed our black hole toy model in the framework of noncommutative geometry. Indeed, the quantization of $AdS_l$ is the data of the anti de Sitter manifold and the action of the group $AN$; that is precisely the data which defines the black hole of chapter \ref{ChapAdS}. So we would be able to ``see'' the causal issue from the data of the deformed spectral triple. Remark that a causal structure (in the physical meaning of the term) is a special property of \emph{pseudo}-riemannian manifolds for which spectral geometry does not exist yet.

An important remaining problem with that method is the fact that the extension lemma does not assure the existence of a stable functional space for the new product. So there is still a lot of analytic work to be done.

\appendix

\setcounter{chapter}{0}
\setcounter{section}{0}
\renewcommand{\theequation}{\Alph{chapter}.\arabic{equation}}
\renewcommand{\thenumtho}{\Alph{chapter}.\arabic{numtho}}
\renewcommand{\thechapter}{\Alph{chapter}}

\chapter{Deformations}		\label{ChapDefo}

\begin{abstract}
Deformation is a main theme of research in the present work. We begin here to describe WKB quantization and a general method to guess deformations of function algebras. The role of Darboux charts and momentum maps appears clearly. A careful example is given by the deformation of $\SL(2,\eR)$.

We prove a useful result (from \cite{articleBVCS}), the extension lemma, which allows to deform a split extension when one knows a deformation of the two components of the extension. The kernel is simply the product of the two kernels.

Then we see the principle of deformation by group action: when a Lie group is deformable, one can find a deformation of any manifold on which the group acts. Universal formulas exist in some cases. This is why deformations of groups are studied. An application of that extension lemma to the Iwasawa subgroup of $\SO(2,n)$ is given in chapter \ref{ChapNoteDev}.

\end{abstract}

\section{WKB quantization}\label{subsec:WKB}

More details can be found in the article \cite{lcBBM}. A manifold $M$ is given with its usual commutative and associative algebra $(C^{\infty}(M),\cdot)$ of smooth functions. A \defe{deformation}{Deformation}, or a \emph{quantization}\index{Quantization}\footnote{In fact, we make a difference between these two words. A \emph{deformation} is only the fact to find a new product from an old one; the new product depends on a parameter and has to reduce to the old one when the parameter goes to zero. A \emph{quantization} is a deformation in which the first order term (whatever it means) of the new product contains the symplectic structure as in condition \eqref{EqExigSymplePremOrd} below.}, of $M$ is the data of a new product $\star^M_{\hbar}$ on a functional space over $M$. 

Let $G$ be a Lie group acting on a manifold $M$. We consider $\Fun(M,\eC)$\nomenclature{$\Fun(M)$}{Functions on the manifold $M$}, the space of all the maps from $M$ to $\eC$, without any regularity conditions. The \defe{regular left representation}{Regular!representation}\index{Representation!regular left} of $G$ on $M$ is the representation of $G$ on $\Fun(M)$ given by
\begin{equation}
   [ L^*_g(a) ](h)=a(gh)
\end{equation}
for all $a\in\Fun(M)$, $g$, $h\in G$.

\subsection{Definitions and general setting}

Let $(M,\omega,\nabla)$ be an affine symplectic manifold, i.e. a $2n$-dimensional symplectic manifold $(M,\omega)$ endowed with a torsion-free connection $\nabla$ such that $\nabla\omega=0$. The \defe{automorphism}{Automorphism!of affine symplectic manifold} group $\Aut(M,\omega,\nabla)$\nomenclature{$\Aut(M,\omega,\nabla)$}{Automorphism group of an affine symplectic manifold} is defined as
\[ 
  \Aut(M,\omega,\nabla)=\gpAff(\nabla)\cap\gpSymp(\omega)
\]
where $\gpAff(\nabla)$ is the group of affine transformations of the affine manifold $(M,\nabla)$ and $\gpSymp(\omega)$ is the group of symplectomorphisms of $(M,\omega)$.

Let $R$ be a subgroup of $\Aut(M,\omega,\nabla)$. The following definition of a $R$-invariant WKB quantization can be found in \cite{StrictSolvableSym}.
\begin{definition}
A $R$-invariant \defe{WKB quantization}{WKB quantization} of $(M,\omega,\nabla)$ is the data of a product 
\begin{equation}	\label{EqFormeWKBdel}
(u\star_{\theta}v)(x)=\frac{1}{ \theta^{2n} }\int_{M\times M} a_{\theta}(x,y,z) e^{\frac{ i }{ \theta }S(x,y,z)}u(y)v(z)\,dy\,dz
\end{equation}
(where $dy\,dz$ is the Liouville measure $\omega^{n}/n!$) with the following constrains:
\begin{enumerate}
\item For each $\theta$, we have a space $A_{\theta}$ containing the space $ C^{\infty}_{c}(M)$ of compactly supported smooth functions. The product $\star_{\theta}$ extends to $A_{\theta}$ in such a way that $(A_{\theta},\star_{\theta})$ becomes a one-parameter family of associative $*$-algebras.
\item The product $\ast_{0}$ on $A_{0}$ is the usual pointwise product  and $(A_{0},\ast_{0})$ is a Poisson subalgebra of $ C^{\infty}(M)$ for the induced Poisson structure from the symplectic form $\omega$.
\item $\forall \theta\geq 0$, the space $A_{\theta}$ is a $*$-vector subspace of $ C^{\infty}(M)$ such that \[ 
   C^{\infty}_{c}(M)\subset A_{0}\subset A_{\theta}
\]
where the involution $*$ on $  C^{\infty}(M)$ is the usual complex conjugation.
\item $S$ is a real valued smooth function $S\colon M\times M\times M\to \eR$ such that for all $x_{0}\in M$, the function $S(x_0,.,.)\in C^{\infty}(M\times M)$ has a nondegenerate critical point at $(x_0,x_0)$.
\item The functions $a_{\theta}$ are positive real-valued: 
\[ 
  a_{\theta}\colon M\times M\times M\to \eR^{+}.
\]
\item The functions $S$ and $a_{\theta}$ are invariant under the diagonal action of $R$ on $M\times M\times M$.
\item $\forall\, x\in M$ and $\forall\,u,v\in C^{\infty}_{c}(M)$ with support in a suitably small neighbourhood of $x$, a stationary phase method yields the extension
\begin{equation}  
  (u\star_{\theta} v)(x)\sim u(x)v(x)+\frac{ \theta }{ i }c_{1}(u,v)(x)+o(\theta^{2})
\end{equation}
where $c_{1}$ satisfies
\begin{equation}		\label{EqExigSymplePremOrd}
  c_{1}(u,v)-c_{1}(v,u)=2\{ u,v \}.
\end{equation}

\end{enumerate}
\label{DefWKBCompl}
\end{definition}

We emphasize the fact that the functional space $A^M$ is stable under $\star_{\theta}$: this is a \emph{strict} quantization in contrast to a \emph{formal} star product which only stabilises the space of formal power series of $\theta$.

An example of WKB quantization is the Weyl product which is nothing but an integral reformulation of the Moyal star product:
\[ 
  (f\star^W_{\hbar}g)(x)=\frac{1}{ \hbar^{2n} }\int_{\mathbb{R}^{2n}\times\mathbb{R}^{2n}}  e^{\frac{ 2i }{ \hbar }S^0(x,y,z)}f(y)g(z)\,dy\,dz
\]
where $S^0(x,y,z)=\Omega(x,y)+\Omega(y,z)+\Omega(z,x)$, and $\Omega$ denotes the usual symplectic form on $\mathbb{R}^{2n}$.

The function $K=a_{\theta} e^{\frac{ i }{ \theta }S}$ is the \defe{kernel}{Kernel of a WKB quantization} of the product $\star_{\theta}$. The \defe{associativity}{Associativity!of a WKB quantization} of $\star_{\theta}$ on the functional space $A_{\theta}$ is the fact that the equality
\[ 
  \big( (u\star_{\theta}v)\star_{\theta}r \big)(x)=\big( u\star_{\theta}(v\star_{\theta}r) \big)(x)
\]
holds for every $u$, $v$, $r\in A_{\theta}$ and $x\in M$.  That condition translates under an integral form to the following relation
\begin{equation}\label{EqCondAssoc}
\begin{split}
&\int_{M\times M}K(x,y,z)\left[ \int_{M\times M}K(y,t,s)u(t)v(s)\mu_M(t,s) \right] r(z)\mu_M(y,z)\\
&=\int_{M\times M}K(x,y,z)u(y)\left[ \int_{M\times M}K(z,t,s)v(t)r(s)\mu_M(t,s) \right]\mu_M(y,z)
\end{split}
\end{equation}
where $\mu_M(y,z)=\mu_M(y)\mu_M(z)$ is the Liouville measure on $M$. Performing formal manipulations (such as a Fubini theorem), one can express this condition as
\begin{equation}		\label{EqAssosssens}
\int_{M}K(x,y,t)K(t,p,q)\mu(t)=\int_{M}K(x,\tau,q)K(\tau,y,p)\,\mu(\tau).
\end{equation}
That form is easier to handle and to check, but it is meaningless in general.

The fact to have a \defe{left-invariant kernel}{Left-invariant!kernel} on a group $G$ means that the kernel $K\colon G\times G\times G\to \eC$ has the property $L_g^*K=K$, or
 \begin{equation}
K(gh_{1},gh_2,gh_{3})=K(h_1,h_2,h_{3})
\end{equation}
for every $g\in G$.  The following lemma allows us to use group isomorphisms to push forward a kernel from a group to another.
\begin{lemma}   
Let $G_{1}$ and $G_{2}$ be two symplectic Lie groups and $K_1$, a left-invariant kernel on $G_{1}$ which provides an associative product on the functional space $A_1$. Let $\phi\colon G_{2}\to G_{1}$ be a symplectic Lie group isomorphism. Then the kernel $K_2=\phi^*K_1$ is invariant and gives rise to an associative product on $A_2=\phi^*A_1$.
\label{LemKerINvarIsom}
\end{lemma}

\begin{proof}
By definition,
\[ 
(\phi^*K_1)(h_1,h_2,h_3)=K_1\big( \phi(h_1),\phi(h_2),\phi(h_3)\big).
\]
Therefore, using the left-invariance of $K_1$, we have
\[ 
\begin{split}
	L_{g}^*\phi^*K_2=(\phi\circ L_{g})^*K_2=(L_{g}\circ\phi)^*K_2=\phi^*L_{\phi(g)}^*K_1=\phi^*K_1.
\end{split}  
\]
That proves left-invariance of $\phi^*K_1$ on $G_{2}$.  Now we prove the associativity of $K_2$, this is to check condition  \eqref{EqCondAssoc}. We have
\[
\begin{split}
 &\int_{G_2\times G_2}K_2(x,y,z)\Bigg[ \int_{G_2\times G_2}K_2(y,t,s)(\phi^*u)(t)(\phi^*v)(s)\mu_2(t,s) \Bigg]\\
&\qquad(\phi^*r)(z)\mu_2(y,z)\\
=&\int_{G_2\times G_2}K_1(\phi x,\phi y,\phi z)\Bigg[  \int_{G_2\times G_2}K_1(\phi y,\phi t,\phi s)u(\phi t)v(\phi s)\mu_2(t,s)  \Bigg]\\
&\qquad r(\phi z)\mu_2(y,z).
\end{split}
\]
We perform in this integral the change of variables $\tau_y=\phi y$, $\tau_t=\phi t$, $\tau_z=\phi z$ and $\tau_s=\phi s$. This does not affect the measure because  $\phi$ is a symplectomorphism and $\mu_i$ are the Liouville measures on $G_i$, so that for example,  $\mu_2(t)=\mu_2(\phi^{-1}\tau_t)=\mu_1(\tau_t)$. The previous integral becomes
\[ 
\begin{split}
  &\int_{G_1\times G_1}K_1(\phi x,\tau_y,\tau_z)\Bigg[  \int_{G_1\times G_1}K_1(\tau_y,\tau_t,\tau_s)u(\tau_t)v(\tau_s)\mu_1(\tau_t,\tau_s)  \Bigg]\\
&\qquad r(\tau_z)\mu_1(\tau_y,\tau_z).
\end{split}
\]
Using now the associativity of $K_1$ on $G_1$ and performing the inverse change of variables, we find
\[ 
\begin{split}
\int_{G_2\times G_2}K_2(x,y,z)(\phi^*u)(y)\Bigg[ \int_{G_2\times G_2}K_2(z,t,s)(\phi^*v)(t)(\phi^*r)(s)&\mu_2(t,s)   \Bigg]\\
					&\mu_2(y,z),
\end{split}  
\]
which proves the associativity of $K_2$ on $\phi^*A_1$.

Notice that condition \eqref{EqAssosssens} can be checked in much the same way.

\end{proof}

 It is worth noticing that lemma \ref{LemKerINvarIsom} needs a group isomorphism while one often only has a Lie algebra isomorphism. Due to Campbell-Backer-Hausdorff formula, it may be very difficult to find a group isomorphism from an algebra one. 

\begin{remark}
Most of the time, the symplectic condition \eqref{EqExigSymplePremOrd} does not have to be checked because we just define the symplectic form $\omega_2$ on $G_2$ as $\omega_2=\phi^*\omega_1$ where $\omega_1$ is the symplectic form on $G_1$.
\end{remark}

\begin{definition}
When $\alpha\colon G\times A\to A$ is an action of a Lie group $G$ on a vector space $A$, one says that the element $a\in A$ is a \defe{differentiable vector}{Differentiable vector} of $\alpha$ if the map $g\mapsto\alpha_g(a)$ is a differentiable map from $G$ into $A$.
\end{definition}

We are interested in the regular left representation $L\colon R\times A_{\theta}\to  A_{\theta}$ defined by $\big( L_r(u) \big)(x)=u(r\cdot x)$. A function $u\in A_{\theta}$ is a differentiable vector of $L$ when the map 
\begin{equation}
\begin{aligned}
 \alpha_u\colon R&\to  A_{\theta} \\ 
r&\mapsto L_r(u) 
\end{aligned}
\end{equation}
is differentiable. The differential of $\alpha_u$ is what we will denote by $dL$ in the next few pages: $dL(X)u=(d\alpha_u)_eX$. By definition, 
\[ 
  (d\alpha_u)_eX=\Dsdd{ \alpha_u( e^{tX}) }{t}{0}=\Dsdd{ L_{ e^{tX}}(u) }{t}{0},
\]
and the element $(d\alpha_u)_eX\in A_{\theta}$ applied to $x\in M$ is
\begin{equation}
 \big( dL(X)u \big)(x) =\Big( (d\alpha_u)X \Big)(x)=\Dsdd{ L_{ e^{tX}}(u)x }{t}{0}=\Dsdd{ u( e^{tX}\cdot x) }{t}{0}.
\end{equation}
We denote by $ A_{\theta}^{\infty}$ the space of differentiable vectors of the representation $L$.

If one particularises to  $ A_{\theta}\subset  C^{\infty}(R)$ (the manifold $M$ being $R$ itself), the vector fields of $R$ naturally act on $ A_{\theta}$. In particular, if $u\colon R\to \eC$ and $X\in \mR$ we have
\[ 
  \big( X^*(u) \big)(r)=X^*_r(u)=\Dsdd{ u\big(  e^{-tX}r \big) }{t}{0}=\big( dL(-X)u \big)(r), 
\]
so that
\begin{equation}
   dL(X)=-X^*
\end{equation}
holds on the space of differentiable vectors $ A_{\theta}^{\infty}$.

\begin{definition}
A formal star product $\dpt{\ast_G}{\Cinf(M)\dcr{\nu}\times \Cinf(M)\dcr{\nu}}{\Cinf(M)\dcr{\nu}}$ is said to be \defe{$\mG$-covariant}{Covariant star product} if for all $X$, $Y\in\mG$,
\begin{equation}
[\lambda_X,\lambda_Y]_{\ast_G}=2\nu\{\lambda_X,\lambda_Y\}
\end{equation}
where $[\lambda_{X},\lambda_{Y}]_{\ast_{G}}:=\lambda_X\ast_G\lambda_Y-\lambda_Y\ast_G\lambda_X$. In other words the start product is $\mG$-covariant when the expected terms of higher order in the right hand side are zero. 
\end{definition}

A crucial use of $\mG$-covariance will be done in proposition \ref{Proprhonureprez} in order to build a map $\rho_{\nu}$ that fulfils the following proposition (instead of $dL$ itself). 
\begin{proposition}
In the setting of definition \ref{DefWKBCompl}, the map $dL$ is a representation  by derivation of $\mR$  on~$ A_{\theta}^{\infty}$.
\label{prop:dL_reprez}
\end{proposition}

\begin{proof}
We will not pay attention on the domain $A_{\theta}$. Its definition will come later.  First, we prove that $\dpt{dL}{\mR}{\End{ A_{\theta}^{\infty}}}$ is a representation. Indeed,
\begin{equation}
\begin{split}
  dL([X,Y])u=\Dsddp{L^*_{\exp(-t[X,Y])}u}{t}{0}
            &=\Dsddp{ [ L^*_{\exp(-tX)},L^*_{\exp(-tY)}  ] u}{t}{0}\\
	    &=\big[dL(X),dL(Y)\big]u.
\end{split}
\end{equation}
Next, $L_R$-invariance of $\ast_{\theta}$ yields 
\[ 
  \big( L^*_{\exp -tX}u \big)\ast_{\theta}\big( L^*_{\exp -tX}v \big)=L^*_{\exp -tX}(u\ast_{\theta} v).
\]
If we derive this equality with respect to $t$ at $t=0$, we find
\[
   dL(X)u\ast_{\theta} v+u\ast_{\theta} dL(X)v=dL(X)(u\ast_{\theta} v).
\]
\end{proof}

\subsection{Deformation of Iwasawa subgroups}	

The motivation in deforming (or quantizing) groups resides in the method of deformation by group action (appendix \ref{SecDefAction}) which states that if one can deform a group, one can write a formula for a deformed product on any manifold on which the group acts.

Let first describe the next few steps in the construction of WKB quantizations of groups. Let $G$ be a semisimple Lie group with its Iwasawa decomposition $G=ANK$. The group $R=AN$ is solvable and can be seen as the homogeneous space $R=G/K$. We consider the canonical multiplicative action $\dpt{\tau}{G\times R}{R}$ which we restrict to $\dpt{\tau}{R\times R}{R}$. We are interested in a $R$-invariant quantization of $R$. Here is a summary of the notations that will be used.

\begin{itemize}
\item $\stM$ is the Moyal star product on $\eR^n$ endowed with its canonical symplectic form,
\item $\star^R_{\theta}$ is the product we are searching for. It has to be defined at least on $ C_c^{\infty}(R)$ and should be extended to $ C^{\infty}(R)$,
\item $A^R\subset\Fun(R,\eC)$ must contain $ C^{\infty}_c(R)$. The purpose is $(A^R,\star^R_{\theta})$ to be an associative algebra and $A^R$ to be invariant under the left regular representation of $R$,
\item $\eA_{\nu}= C^{\infty}(R)[ [\nu]]$ is an intermediary space which serves to guess $\star^R_{\theta}$ and perform formal manipulations with $\rho_{\nu}$ and $dL$,
\item $\ast_M^R$ is the pull-back of Moyal to $\eA_{\nu}$. It serves to formal manipulations in order to guess the twist that defines $\ast^R_{\nu}$,
\item $\ast^R_{\nu}$ is the product on $\eA_{\nu}$. The problem of determining that product is formal. When this problem is solved, we have to prove that in a well chosen $A^R$, taking $\ast_{\nu}^R\to\star^R_{\theta}$ yields a solution to the problem. As previously noticed, in order to make sense, one has to apply $dL$ on the subspace $\eA_{\nu}^{\infty}$ of differentiable vector of the regular left representation. We will however not take care of this issue in the formal manipulations.
\end{itemize}

The main steps are the following:
{\renewcommand{\theenumi}{\arabic{enumi}.}
\begin{enumerate}
\item In the case of a WKB product we show in proposition \ref{prop:dL_reprez} that $dL$ is a representation of $\mR$ on $\eA_{\nu}^{\infty}$. Hence we will try to build a formal product for which $dL$ is a representation by derivation. From this point of view, the manipulation with $\rho_{\nu}$ is only a trick designed to guess a product formula.

\item We suppose that the group $R$ ---the one that we are trying to quantize--- has a symplectic structure $\omega$ and we consider $\phi\colon \eR^{2n}\to R$, a Darboux chart; i.e. $\omega=\phi^*\Omega$ where $\Omega$ is the canonic symplectic form on $\eR^{2n}$.

\item We suppose that the left action of $R$ on itself is strongly hamiltonian and we denote by $\lambda_X$ the momentum maps. We suppose that the Moyal product is $\mG$-covariant\footnote{In fact, we only need the $\mR$-covariance.}.

\item We pose $\rho_{\nu}(X)=\frac{1}{ 2\nu }\ad_{\ast_M^R}(\lambda_X)$. The $\mR$-covariance of $\ast_M^R$ is used in order to prove that $\rho_{\nu}$ is a  representation by derivations of $\mR$ on $(\eA_{\nu},\ast_M^R)$. 
\item If one can find an intertwining operator between $dL$ and $\rho_{\nu}$ (i.e. if they are equivalent representations), we define $\ast_{\nu}^R$ as the pull-back of $\ast_M^R$ by this intertwining operator. In this case, we prove that $dL$ is a representation by derivations of the product $\ast_{\nu}^R$.
\end{enumerate}

}		%

We try now to find a formal product $\ast_{\nu}^R$ on $\eA_{\nu}^{\infty}$ such that $dL$ is a representation by derivations. For this purpose we suppose $R$ to accept a symplectic structure $\omega$ and $\phi\colon \eR^{2n}\to R$ to be a Darboux chart, i.e. $\omega=\phi^*\Omega$ where $\Omega$ denotes the canonical symplectic form on $\eR^{2n}$. Then we bring the Moyal product of $\eR^{2n}$ to $R$ by the usual formula
\begin{equation}   
	(u\ast_M^R v)=(u\circ\phi\ast_M v\circ\phi)\circ\phi^{-1}.
\end{equation}
We suppose that product to be $\mG$-covariant\footnote{Only the $\mR$-covariance will be actually used.}:
\begin{equation}
  [\lambda_X,\lambda_Y]_{\ast_M^R}=2\nu \{ \lambda_X,\lambda_Y \}_R.
\end{equation}
 Now we consider the left action of $R$ on itself and we suppose that this is an Hamiltonian action for the symplectic structure $\omega=\phi^*\Omega$ with dual momentum maps $\dpt{\lambda_X}{R}{\eC}$. We define, for each $X\in\mR$, a linear map, $\dpt{\rho_{\nu}(X)}{\eA_{\nu}}{\eA_{\nu}}$ by
\begin{equation}
\begin{aligned}
 \rho_{\nu}\colon \mR&\to \End{\eA_{\nu}} \\ 
X&\mapsto\us{2\nu}\ad_{\ast_M^R}(\lambda_X)
\end{aligned}
\end{equation}
Notice that the formal series of $[\lambda_X,u]_{\ast_M^R}$ begins with order one, so the division by $\nu$ make sense in the space of formal series.  The main interest of $\rho_{\nu}$ is to be as we want $dL$ to be. So it will be used to guess how to twist the product in order to make $dL$ work as $\rho_{\nu}$.

\begin{proposition}
The map $\rho_{\nu}$ is a representation of $\mR$ on $\eA_{\nu}$, and $\rho_{\nu}(X)$ is a derivation of $(\eA_{\nu},\ast_M^R)$ for each $X\in\mR$.
\label{Proprhonureprez}
\end{proposition}

\begin{proof}
The proof  that $\rho_{\nu}$ is a representation is only to check that the relation $[\rho_{\nu}(X),\rho_{\nu}(Y)]f=\rho_{\nu}([X,Y])f$ holds for any $X$, $Y\in\mR$ and $f\in\eA_{\nu}$. Using the $\mG$-covariance and the Jacobi identity,
\begin{equation}
\begin{split}
  \rho_{\nu}([X,Y])f&=\us{4\nu^2}\ad_{\ast_M^R}(2\nu\lambda_{[X,Y]})f
     		=\us{4\nu^2}\ad_{\ast_M^R}([\lambda_X,\lambda_Y]_{\ast_M^R})f\\
		&=\us{4\nu^2}[[\lambda_X,\lambda_Y]_{\ast_M^R},f]_{\ast_M^R}\\
	     	&=\frac{1}{ 4\nu^2 }(\ad_{\ast_M^R}\lambda_X\circ\ad_{\ast_M^R}\lambda_Y-\ad_{\ast_M^R}\lambda_Y\circ\ad_{\ast_M^R}\lambda_X)f\\
		&=[\rho_{\nu}(X),\rho_{\nu}(Y)]f.
\end{split}
\end{equation}
It remains to check that $\rho_{\nu}(X)(u\ast_M^R v)=\rho_{\nu}(X)u\ast_M^R v+u\ast_M^R\rho_{\nu}(X)v$ for every $X\in\mR$. This is once again just a computation.
\begin{equation}
\begin{split}
   \rho_{\nu}(X)u\ast_M^R v+u\ast_M^R\rho_{\nu}(X)v&=\us{2\nu}(\lambda_X\ast_M^R u-u\ast_M^R\lambda_X)\ast_M^R v\\
                           &\quad+\us{2\nu}u\ast_M^R(\lambda_X\ast_M^R v-v\ast_M^R\lambda_X)\\
                           &=\us{2\nu}\ad_{\ast_M^R}\lambda_X(u\ast_M^R v).
\end{split}
\end{equation} 
\end{proof}
Notice that the $\mG$-covariance of $\ast_M^R$ was used to prove that $\rho_{\nu}$ is a representation.  Now, if we could show that $\rho_{\nu}=dL$, then the answer to our deformation problem would be $A_{\theta}=\eA^{\infty}_{\nu}$ and $\stt=\ast_M^R$. But instead of that we have $\rho_{\nu}=dL+o(\nu)$ because
\begin{equation}
\begin{split}
  \rho_{\nu}(X)u&=\us{2\nu}[\lambda_X,u]_{\ast_M^R}
         =\us{2\nu} 2\nu\{\lambda_X,u\}+o(\nu)
	 =X^*(u)+o(\nu)\\
	&=-dL(X)u+o(\nu)
\end{split}
\end{equation}
where the notion of fundamental field $X^*$ is taken for the regular left representation (which is Hamiltonian). That shows that $\rho_{\nu}$ is something like a deformation of $dL$. As a consequence, one has $dL(X)=X_{\lambda_X}$, or
 \begin{equation}\label{eq:dL_et_Poisson}
 dL(x)u=X_{\lambda_X}(u)=\{\lambda_X,u\}
 \end{equation}
(see subsection \ref{app:ham_act}).
 
Since $\rho_{\nu}$ is not $dL$, the hope is to see if $\rho_{\nu}$ and $dL$ should be \emph{equivalent} representations. As next proposition shows, the fact to find an equivalence between $\rho_{\nu}$ and $dL$ actually solves the problem to find a product for which $dL$ is a representation by derivation. 
\begin{proposition}
Let $\dpt{\mT}{\eA_{\nu}}{\eA_{\nu}}$ be an intertwining operator between $dL$ and $\rho_{\nu}$:
\begin{equation}\label{eq:TrnT} 
   \mT\rho_{\nu}(X) \mT^{-1}=dL(X).
\end{equation} 
If we define the star product $\ast_{\nu}^R$ by
\begin{equation}	\label{Eq_candprodANSL}
   u\ast_{\nu}^R v=\mT\bnu(\mT\bnu^{-1} u\ast^{R}_M \mT\bnu^{-1} v),
\end{equation}
$dL$ becomes a derivation of $\ast_{\nu}^R$.
\label{prop:def_stn}
\end{proposition}

\begin{proof}
If we develop the expression of $dL(X)(u\ast^R_M v)$, we find $\mT\rho_{\nu}(X)(\mT^{-1} u\ast^R_M \mT^{-1} v)$, using the fact that $\rho_{\nu}$ is a derivation of $\ast^R_M$, one easily finds $dL(X)u\ast^R_M v+u\ast_{\nu}^R dL(X)v$.
\end{proof}
\section{Deformation of \texorpdfstring{$\SL(2,\eR)$}{SL2R}}		\label{sec:unifsl}

\begin{abstract}
This section shows in some detail an instructive example of deformation of an Iwasawa subgroup: the Iwasawa subgroup of $\SL(2,\eR)$.
In this section we will use the parametrization \eqref{EqParmalSL} of $\SL(2,\eR)$, as well as the notations $G=\SL(2,\eR)$ and $\mG=\sldr$. Here are the main steps that will be performed:
{
\renewcommand{\theenumi}{\arabic{enumi}.}

\begin{enumerate}
\item The Iwasawa component $R=AN=G/K$ provides a double covering onto $\mO=\Ad(G)Z$ where $Z$ is any element of $\mK$ (which is one dimensional). The adjoint orbit $\mO$ being endowed with a canonical symplectic form described in subsection \ref{sub:coadjoint}, we consider on $R$ the corresponding symplectic structure.

\item The map $(a,l)\mapsto \Ad( e^{aH} e^{lE})Z$ turns out to be a global Darboux chart and induces the diffeomorphism
\[ 
  R\simeq \mO\simeq \eR^2.
\]
Under these identifications, the adjoint action of $R$ on $\mO$ becomes the simple multiplication of $R$ in itself, which is strongly hamiltonian.
\item  The Moyal product is $\gsl(2,\eR)$-covariant for the action of $\SL(2,\eR)$ on $\eR^2$.

\item We explicitly build the intertwining operator between $\rho_{\nu}$ and $dL$ and we write down a product (see proposition \ref{prop:def_stn}).

\item A theorem is stated in which we list the properties of the so constructed product.

\end{enumerate}

}
\end{abstract}

\subsection{Actions and Symplectic structure}

Here, in contrast with the case studied in \ref{sub:coadjoint}, we are working with adjoint orbits (and not the \underline{co}adjoint orbits), so the subalgebra to be studied is no more $\widetilde{\mO}$ but
\[
    \mO=\Ad(G)Z,
\]
and the symplectic form is not exactly \eqref{eq:omega_Gs}, but
\begin{equation}\label{eq:omega_G}
  \omega_X(A^*,B^*)=B(X,[A,B]).
\end{equation}
The action of $G$ on $\mO$ is $g\cdot X=\Ad(g)X$. The corresponding notion of fundamental field is given by
\[
   X^*_{\phi(a,l)}=\Dsdd{ \Ad(e^{-tX})\phi(a,l) }{t}{0}.
\]
The Iwasawa theorem \ref{ThoIwasawaVrai} claims that $G/K=AN$ and that we have global diffeomorphism $\mA\oplus\mN\to AN$, $(a,n)\to e^ae^n$; $\mA\to A$, $a\to e^a$; $\mN\to N$, $n\to e^n$. We define $\mR=\mA\oplus\mN$ and the global diffeomorphism 
\begin{equation}	\label{EqDefphiaHlEZ}
\begin{aligned}
 \phi\colon \mA\oplus\mN&\to \mO \\ 
 aH+lE&\mapsto \Ad( e^{aH} e^{lE})Z.
\end{aligned}
\end{equation}
 That map can also be seen as 
\begin{equation}
\begin{aligned}
 \phi\colon\eR^2&\to \mO \\ 
(a,l)&\mapsto\Ad(e^{aH}e^{lE})Z.
\end{aligned}
\end{equation}
 In this way, we identify $\mA\oplus\mN$ and $\eR^2$ as two dimensional space.
\begin{proposition}
As homogeneous space, there is a double covering
\begin{equation}
\begin{aligned}
 \psi\colon G/K&\to \mO \\ 
[g]&\mapsto \Ad(g)Z. 
\end{aligned}
\end{equation}

\end{proposition}
\begin{proof}
The map $\psi$ is well defined and injective (up to the double covering) because the stabilizer of $\mK$ is $K$. The surjective condition is clear. The \emph{double} covering is expressed by the fact that $\psi([g])=\psi([g'])$ if and only if $g=\pm g'$.
\end{proof}

The symplectic $2$-form $\omega$ on $\mO$ induces a symplectic form  
\[
  \Omega=\phi^*\omega
\]
 on $\mA\oplus\mN\simeq\eR^2$.

\begin{proposition}
   The $2$-form $\phi^*\omega$ is constant and its value is
\[
              \Omega:=\phi^*\omega=-2B(F,E)da\wedge dl=\beta da\wedge dl;
\]
in other words, $\phi$ is a $\emph{global}$ Darboux chart for $\mO$.
\label{prop:Omega}
\end{proposition}

\begin{proof}
We have to compute
\[
   \Omega_{(a,l)}(\partial_a,\partial_l)=\omega_{\phi(a,l)}\big( (d\phi)_{(a,l)}\partial_a,(d\phi)_{(a,l)}\partial_l\big).
\]
First, we show that $d\phi(\partial_a)=-H^*_{\phi}$: 
\[
\begin{split}
   d\phi_{(a,l)}\partial_a&=\Dsdd{\phi(a+t,l)}{t}{0}
                    	=\Dsdd{\Ad(e^{(a+t)H}e^{lE})Z}{t}{0}\\
                    	&=\Dsdd{ \Ad(e^{tH}e^{aH}e^{lE})Z  }{t}{0}
		    	=\Dsdd{ \Ad(e^{tH})\phi(a,l) }{t}{0}\\
		    	&=-H^*_{\phi(a,l)}.
\end{split}
\]
In the same way, we find $d\phi(\partial_l)=\big( \Ad(e^{aH})E\big)^*_{\phi}$:
\[
 \begin{aligned}
   d\phi_{(a,l)}\partial_l&=\Dsdd{ \Ad(e^{aH}e^{l+tE})Z }{t}{0}
                    =\Dsdd{ \Ad(e^{aH}e^{tE} e^{-aH}e^{aH}e^{lE} )Z }{t}{0}\\
		    &=\Dsdd{ \Ad(e^{aH}e^{tE}e^{-aH})\phi(a,l) }{t}{0}
		    =\Dsdd{ \Ad( e^{t\Ad(e^{aH})E} )\phi(a,l) }{t}{0}\\
		    &=-\big( \Ad(e^{aH})E\big)^*_{\phi(a,l)}.
\end{aligned}
\]
Using formula \eqref{eq:omega_G} for the symplectic form, 
\begin{equation}
\begin{split}
  \Omega_{(a,l)}(\partial_a,\partial_l)&=B\big( \phi(a,l),[-H,-\Ad(e^{aH})E]\big) \\
                           &=B\big( \Ad(e^{aH})\Ad(e^{lE})Z,\Ad(e^{aH})[H,E]\big)\\
			   &=2B\big( Z,\Ad(e^{-lE})E \big)\\
			   &=2B(Z,E).
\end{split}
\end{equation}
Defining $\beta=-2B(E,F)$ we write it as
\begin{equation} 
  \Omega=\phi^*\omega=-2B(F,E)da\wedge dl=\beta da\wedge dl.
\end{equation}
\end{proof}
So, as symplectic manifold, $(\mO,\omega)$ is nothing but $(\eR^2,da\wedge dl)$, the diffeomorphism being $\phi$. The symplectic structure $\Omega$ induces a Poisson structure $P$ given by equation \eqref{eq:def_Poisson}. In the present case, it reads
\begin{align}
(\Omega_{ij})&=\beta\begin{pmatrix}
0 & 1 \\
-1 & 0
\end{pmatrix}
&(P)&=\beta^{-1}\begin{pmatrix}
0 & -1 \\
1 & 0
\end{pmatrix}
\end{align}
 and
\begin{equation}\label{eq:Poisson}
  \{f,g\}=\beta^{-1}(\partial_lf\partial_ag-\partial_af\partial_lg).
\end{equation}

The action of $G$ on $\mO$ can be turned into an action on $\eR^2$ using the chart $\phi$. It is done by defining $\dpt{\tau}{G\times\eR^2}{\eR^2}$,
\begin{equation}
   \tau=\phi^{-1}\circ \Ad\circ\phi,
\end{equation}
or $\tau_g(a,l)=\phi^{-1}\big( \Ad(g)\phi(a,l)\big)$.  The notion of fundamental field\index{Fundamental!vector field!on $\eR^2$} at $x=(a,l)\in\eR^2$ is thus given by
\begin{equation}
  X^*_x=\Dsdd{e^{-tX}\cdot x}{t}{0}
       =\Dsdd{ \phi^{-1}\big( \Ad(e^{-tX})\phi(a,l)  \big) }{t}{0},
\end{equation}
for which we will often use the path representation
\[
   X^*_x(t)=\phi^{-1}\big( \Ad(e^{-tX})\phi(a,l)  \big).
\]
From $\Ad$-invariance of $\omega$,
\[
   \tau^*\Omega=\tau^*\phi^*\omega
               =(\phi\circ\phi^{-1}\circ \Ad\circ\phi)^*\omega
	       =\phi^*(\Ad)^*\omega
	       =\phi^*\omega
	       =\Omega.
\]
Thus the symplectic form is $G$-invariant:
\begin{equation}     \label{eq:tau_s_Omega}
  \tau^*\Omega=\Omega,
\end{equation}
That implies in particular that $\tau$ satisfies theorem \ref{tho:equiv_Poisson}.

\begin{proposition}
The action $\tau$ of $G$ on the symplectic space $(\eR^2,\Omega)$ is Hamiltonian and the dual momentum maps $\dpt{\lambda'_X}{\eR^2}{\eR}$ are given by (cf .\ref{def:app_mom_mom_duale})
\begin{equation}
  \lambda'_X(a,l)=-B\big(X,\phi(a,l)\big)
\end{equation}
for each $X\in\mG$.
\label{prop:lambda_X}
\end{proposition}

\begin{proof}
We have first to check the identity $i(X^*)\Omega=i(X^*)(\phi^*\omega)=d\lambda'_X$. Let us apply both sides on the vector $A^*_x$, with $A\in\mG$ and $x=(a,l)\in\eR^2$. On the one hand
\[
  i(X^*_x)\Omega_x(A^*_x)=\omega_{\phi(x)}\big(   d\phi_xX^*_x,d\phi_xA^*_x   \big),
\]
but  
\begin{equation}
  d\phi_xX^*_x=\Dsdd{\phi(X^*_x(t))}{t}{0}
              =\Dsdd{ \Ad(e^{-tX})\phi(aH,lE) }{t}{0}
	      =-X^*_{\phi(a,l)}.
\end{equation}
The same being true for $A$,
\[
  i(X^*_x)\Omega_x(A^*_x)=\omega_{\phi(x)}(X^*_{\phi(x)},A^*_{\phi(x)})=B(\phi(x),[X,A]).
\]
On the other hand,
\begin{equation}
\begin{aligned}
    (d\lambda'_X)_x(A^*_x)&=\Dsdd{ (\lambda'_X\circ\phi^{-1})\Big(   \Ad(e^{tA})\phi(a,l)   \Big) }{t}{0}\\
                         &=\Dsdd{  B\Big(X,\Ad(e^{tA})\phi(a,l) \Big)  }{t}{0} \\
			 &=B\Big(  \Dsdd{\Ad(e^{tA})\phi(x)}{t}{0},X   \Big)    &\text{$B$ is linear}\\
			 &=B\Big(  (\ad A)\phi(x),X   \Big)\\
			 &=-B\big(\phi(x),(\ad A)X\big) &\text{$B$ is $\Ad$-invariant}\\
			 &=B(\phi(x),[X,A]).
\end{aligned}
\end{equation}
That proves that $i(X^*)\Omega=d\lambda'_X$.  The second part of the proof is to see that condition \eqref{eq:hamil} holds.  Using the fact that $X_{\lambda'_Y}=Y^*$, we find
\[ 
\begin{split}
  \{ \lambda'_X,\lambda'_Y \}(a,l)&=-\Omega(X_{\lambda'_X},X_{\lambda'_Y})
		=-\Omega_{(a,l)}(X^*,Y^*)\\
		&=-\omega_{\phi(a,l)}(X^*,Y^*)
		=-B([X,Y],\phi(a,l))\\
		&=\lambda'_{[X,Y]}(a,l)
\end{split}
\]
where the star refers to the action on $\mO$. Explicit computations of Poisson bracket between $\lambda'_X$'s at page \pageref{pg:explic_com_lamb} will confirm that result.

\end{proof}

We are now able to furnish explicit formulas for $\lambda'_H$, $\lambda'_E$ and $\lambda'_F$ by virtue of the latter proposition.  The first computation is:
\begin{equation}
\begin{aligned}
  \lambda'_H(a,l)&=-B(H,\Ad(e^{lE})Z)
                =-B( \Ad(e^{-lE})H,Z )\\
		&=-B(H+[-lE,H]+\ldots,Z)
		=-B(H,Z)+B([-lE,H],Z)\\
		&=-2lB(E,F),
\end{aligned}
\end{equation}
so
\begin{equation}   \label{EqlamHal}
  \lambda'_H(a,l)=-\beta l.
\end{equation}
Second,
\begin{equation}
\begin{split}
  \lambda'_E(a,l)=-B(\Ad(e^{-aH})E,\Ad(e^{lE}))
		&=-e^{-2a}B(\Ad(e^{-lE})E,Z)\\
		&=-\frac{\beta}{2}e^{-2a}.
\end{split}
\end{equation}
Then,
\begin{equation}  \label{EqlamEal}
\lambda'_E(a,l)=-\frac{\beta}{2}e^{-2a}.
\end{equation}
The last one is
\begin{equation}
\begin{aligned}
\lambda'_F(a,l)&=-B\big(  \Ad(e^{lE})Z,e^{-aH}F  \big)
              =-e^{2a}B\big(Z,  \Ad(e^{-lE})F   \big)\\
	      &=-e^{2a}B\big(Z, F-l[E,F] +\frac{l^2}{2} [E,[E,F]]+\ldots  \big)\\
	      &=-e^{2a}\left[     B(Z,F)-lB(Z,H)-\frac{l^2}{2}B(Z,2E)        \right]\\
	      &=-e^{2a}\big(  B(Z,F)+l^2B(F,E)   \big)\\
	      &=-e^{2a}\big(  -\frac{\beta}{2}-l^2\frac{\beta}{2}   \big)
	      =e^{2a}\frac{\beta}{2}(l^2+1).
\end{aligned}
\end{equation}
Finally,
\begin{equation}  \label{EqlamFal}
\lambda'_F(a,l)=\frac{\beta}{2}e^{2a}(l^2+1).
\end{equation}
Using formula \eqref{EqPoisson} for the Poisson bracket, one can check that the required relations \eqref{eq:hamil} are satisfied:
\begin{subequations}  \label{pg:explic_com_lamb}
\begin{align}
  \{\lambda'_H,\lambda'_E\}&=2\lambda'_E\\
  \{\lambda'_H,\lambda'_F\}&=-2\lambda'_F\\
  \{\lambda'_E,\lambda'_F\}&=\lambda'_H.
\end{align}
\end{subequations}
This confirms the fact that our action of $\SL(2,\eR)$ on $AN$ is Hamiltonian.

Using the global diffeomorphism \eqref{EqDefphiaHlEZ}, and the map 
\begin{equation}
\begin{aligned}
 j\colon AN&\to \mO \\ 
r&\mapsto \Ad(r)Z 
\end{aligned}
\end{equation}
we identify
\[
   R\simeq\mO\simeq\eR^2.
\]
The action of $R$ on itself induced from the adjoint action of $R$ on $\mO$ is 
\[ 
  r\cdot s=j^{-1}\big( r\cdot j(s) \big)=j^{-1}\big( \Ad(rs)Z \big)=rs.
\]
It is the left multiplicative action required in definition \ref{DefWKBCompl}. The Lie group $R$ is endowed with the symplectic form 
\[ 
\omega^R=j^*{\phi^{-1}}^*\Omega.
\]
The notion of fundamental vector for the action of $R$ on itself is given by
\begin{equation}
  X^*_r=\Dsdd{e^{-tX}\cdot r}{t}{0} 
		=\Dsdd{ j^{-1}\big( e^{-tX}\cdot j(r) \big) }{t}{0}
		=dj^{-1} X^*_{j(r)},
\end{equation}
but we know that 
\[ 
  e^{-tX}\cdot j(r)=\Ad(e^{-tX r})Z=[\phi\circ \tau(e^{-tX}r)\circ \phi^{-1}]Z,
\]
 then
\[
X^*_r=dj^{-1}\circ d\phi X^*_{r\cdot \phi^{-1}(Z)}.
\]
If $r=e^{aH}e^{lE}$, then $r\cdot \phi^{-1}(Z)=(a,l)$ and
\begin{equation}
   X^*_r=(dj^{-1}\circ d\phi) X^*_{(a,l)}
\end{equation}
where the fundamental field of the right hand side is taken in the sense of the action of $R$ on $\eR^2$. 

The following proposition shows that the explicit form of $\lambda$ and $\lambda'$ are the same up to natural identifications.

\begin{proposition}  
The left multiplicative action of $R$ on itself is Hamiltonian and the dual momentum maps are given by  $\lambda_X\colon R\to \eC$,
\begin{equation}
\lambda_X= \lambda'_X\circ\phi^{-1}\circ j.
\end{equation}
for each $X\in\mR$.
\label{PropMomslR}
\end{proposition}

\begin{proof}
Once again, the proof is just a verification of the two properties of a momentum map. The first one is
\begin{equation}
\begin{split}
  i(X^*_r)\omega^R Y&=\omega^R_r(dj^{-1} d\phi X^*_{(a,l)},Y)
		=\Omega_{(\phi^{-1}\circ j)(r)}\big( X^*_{(a,l)},d\phi^{-1} dj_r Y \big)\\
		&=(\lambda'_X\circ d\phi^{-1}\circ dj)Y
		=d\lambda_X Y.
\end{split}
\end{equation}
For the second condition, we consider $r=e^{aH}e^{lE}$ and
\begin{equation}
\begin{split}
  \{ \lambda_X,\lambda_Y \}(r)&=X^*_r(\lambda_Y)
		=(dj^{-1} d\phi X^*_{(a,l)})(\lambda'_Y\circ \phi^{-1}\circ j)\\
		&=X^*_{(a,l)}(\lambda'_Y)
		=\lambda'_{[X,Y]}(a,l)\in\eC
\end{split}
\end{equation}
while 
\[ 
  \lambda_{[X,Y]}(r)=\lambda'_{[X,Y]}\circ\phi^{-1}\circ j(r)=\lambda'_{[X,Y]}(a,l).
\]
\end{proof}

\subsection{Guessing the star product}

The Moyal star product is invariant under the action of $\eR^2$ on itself $L_xy=x+y$ in the sense that if we pose $(L_y^*f)(x)=f(x+y)$ it is clear that
\begin{equation}
\begin{split}
    (L_s^*f\ast_M L_s^*g)(x)&=
    \exp\left[{\displaystyle\frac{\nu}{2}P^{ij}(\partial_{y^i}\wedge\partial_{z^j})}\right]f(y+s)g(z+s)|_{y=z=x}\\
	                    &=L^*_s(f\ast_M g)(x).
\end{split}
\end{equation}
We are however not interested by that action on $\eR^2$. The action which we look at is the one of $\SL(2,R)$.

\begin{proposition}
The product $\ast_M$ is $\gsl(2,\eR)$-invariant at order $0$ and $1$.
\end{proposition}

\begin{proof}
The invariance at order zero is given with some concise notations by
\[
 (gu)(gv)(x)=u(gx)v(gx)=(uv)(gx),
\]
The action $\tau_g$ of an element $g\in G$ satisfies $\tau_g^*\Omega=\Omega$ (equation \eqref{eq:tau_s_Omega}), so  theorem \ref{tho:equiv_Poisson} gives $\{u\circ\tau_g,v\circ\tau_g\}=\{u,v\}\circ\tau_g$.  Since Poisson bracket\index{Poisson bracket!and Moyal product} is the first term of the Moyal product, at first order
\[
  \tau_g^*(u\ast_M v)=\tau_g^*u\ast_M\tau_g^*v.
\]

\end{proof}

\begin{proposition}
   The product $\ast_M$ is $\sldr$-covariant for the homomorphism given by proposition \ref{prop:lambda_X} or equivalently by equations \eqref{EqlamHal}, \eqref{EqlamEal}, and \eqref{EqlamFal}.
\end{proposition}

\begin{proof}
The Moyal star product can be written as
\[
   u\ast_M v=\sum \frac{\nu^k}{k!}P_k(u,v)
\]
with $P_k(u,v)=\Omega^{IJ}\partial_Iu\partial_Jv$ where $I$ and $J$ are summed over $k$-uple of $0$ and $1$, including a sum over $k$ itself ($x^0=a$, $x^1=l$). For a given $I$, there is only one $J$ such that $\Omega^{IJ}\neq 0$. There are $\binom{k}{m}$ multi-indices $I$ providing the term $\partial_I=\partial_0^m\partial_1^n$ with $n+m=k$. For each of them, $\Omega^{IJ}=\me{n}$. Therefore
\begin{equation}\label{eq:P_k}
  P_k(u,v)=\sum_{m=0}^k\me{k-m}\binom{k}{m}\partial_0^m\partial_1^nu\,\partial_0^n\partial_1^mv.
\end{equation}
For example,
\[
  P_1(u,v)=-\partial_1u\partial_0v+\partial_0u\partial_1v=\{u,v\}.
\]
If $k$ is even, the expression \eqref{eq:P_k} is symmetric with respect of $u$ and $v$, so that these terms will not contribute in the computation of the commutators $[u,v]_{\ast_M}$. We are left with
\begin{equation}\label{eq:comm_lambda_X}
   [\lambda'_X,u]_{\ast_M}
         =2\sum_{k=0}^{\infty}\frac{\nu^{2k+1}}{(2k+1)!}P_{2k+1}(\lambda'_X,u).
\end{equation}

First we compute $[\lambda'_H,u]_{\ast_M}$:
\begin{equation}
   P_{2k+1}(\lambda'_H,u)=\delta_{k0}(-\partial_1\lambda'_H\partial_0 u+\partial_0\lambda'_H\partial_1 u)=\delta_{k0}\{\lambda',u\},
\end{equation}
thus
\begin{equation}
\begin{split}
  [\lambda'_H,u]_{\ast_M}=2\nu P_1(\lambda'_H,u)
                        =2\nu\{\lambda'_H,u\}
			=2\nu\beta\partial_au.
\end{split}
\end{equation}
By the way, we point out the relation
\[
\ad_{\ast_M}\lambda'_H=2\nu\beta\partial_a.
\]

Now, we turn our attention to the commutator $[\lambda'_E,u]_{\ast_M}$:
\begin{equation}
\begin{split}
  P_{2k+1}(\lambda'_E,u)
      &=-\sum_{n=0}^{k}\me{m}\binom{2k+1}{m}\binom{2k+1}{m}(\partial_0^m\partial_1^n\lambda'_E)\,(\partial_0^n\partial_1^mu)\\
      &=\partial_a^{2k+1}\big( -\frac{\beta}{2}e^{-2a} \big)\partial_l^{2k+1}u 
	=\beta 2^{2k}e^{-2a}\partial_l^{2k+1}u,
\end{split}
\end {equation}
thus
\begin{equation}  \label{eq:comm_lambda_E}
\begin{split}
  [\lambda'_E,u]_{\stM}&=2\sum_{k=0}^{\infty}\frac{\nu^{2k+1}}{(2k+1)!}\beta
                            2^{2k}e^{-2a}\partial_l^{2k+1}u
	             =\beta e^{-2a}\sinh(2\nu\partial_l)u,
\end{split}
\end{equation}
so that
\[
     \ad_{\stM}\lambda'_E=\beta e^{-2a}\sinh(2\nu\partial_l).
\]
Last we check $[\lambda'_E,\lambda'_F]_{\stM}=2\nu\{\lambda'_E,\lambda'_F\}$. When $u=0$, the only non vanishing term in the sum \eqref{eq:comm_lambda_E} is $k=0$. Since $\partial^{3}_{l}\lambda'_{F}=0$,
\[
   [\lambda'_E,\lambda'_F]_{\stM}=2\nu\beta e^{-2a}\partial_l\lambda'_F,
\]
but
\begin{equation}
 2\nu\{\lambda'_E,\lambda'_F\}
             =2\nu(\partial_a\lambda'_E\partial_l\lambda'_F-\partial_l\lambda'_E\partial_a\lambda'_F)
	     =2\nu\beta e^{-2a}\partial_l\lambda'_F.
\end{equation}
\end{proof}

Before to go on, let us compute the operator $\ad_{\stM}\lambda'_F$ in order to complete our collection. We take once again the formula \eqref{eq:P_k}, with $\lambda'_F$ and $u$:
\begin{equation}
   P_{2k+1}(\lambda'_F,u)=-\sum_{m=0}^{2k+1}\me{m}\binom{2k+1}{m}\partial_a^m\partial_l^n\lambda'_F
                                                                \partial_a^n\partial_l^m u.
\end{equation}
It is clear that $\lambda'_F$ can be derived  only two times with respect of $l$ and as much as we want with respect of $a$. Then possible $n$ are $n=0,1,2$, whose corresponding $m$ are $2k-1$, $2k$, and $2k+1$. Some computations lead to
\begin{equation}
\begin{split}
    P_{2k+1}(\lambda'_F,u)&=-k(2k+1)\beta 2^{2k-1}e^{2a}\partial_a^2\partial_l^{2k-1}u\\
                    &\quad +(2k+1)\beta 2^{2k}l\partial_a\partial_l^{2k}u\\
                    &\quad -\beta 2^{2k}(1+l^2)e^{2a}\partial_l^{2k+1}u.
\end{split}
\end{equation}
Replacing into the series \eqref{eq:comm_lambda_X}, we find
\begin{align*}
\begin{split}
[\lambda'_F,u]_{\stM}&=e^{2a}\Big\{
   \sum\frac{\nu^{2k+1}}{(2k+1)!}\beta (-k)(2k+1)\frac{2^{2k+1}}{2}
                     \partial_a^2\partial_l^{2k-1}u\\
 &\qquad+\sum\frac{\nu^{2k+1}}{(2k+1)!}(2k+1)2^{2k+1}l\partial_a\partial_l^{2k}u\\
 &\qquad\beta\sum\frac{\nu^{2k+1}}{(2k+1)!}2^{2k+1}(1+l^2)\partial_l^{2k+1}u\Big\}
\end{split}\\
\begin{split}
 &=-\beta e^{2a}\partial_a^2\sum_{k=1}^{\infty}\frac{(2\nu)^{2k}}{(2k)!}k\nu\partial_l^{2k-1}\\
 &\quad+2\beta\nu e^{2a}\partial_a\circ\cosh(2\nu\partial_l)\\
 &\quad-\beta e^{2a}(1+l^2)\sinh(2\nu\partial_l).
\end{split}
\end{align*}
Finally,
\begin{equation}
\begin{split}
   \ad_{\stM}\lambda'_F&=-\nu^2\beta e^{2a}\partial_a^2\circ\sinh(2\nu\partial_l)\\
                     &\quad+2\nu\beta e^{2a}l\partial_a\circ\cosh(2\nu\partial_l)\\
		     &\quad-e^{2a}(1+l^2)\sinh(2\nu\partial_l).
\end{split}
\end{equation}

\begin{corollary}
The star product $\ast_M^R$ on $R$ defined for $u,v\in C^{\infty}(R)$ by
\begin{equation}
  (u\ast_M^R v)(r)=( u\circ T^{-1}\stM v\circ T^{-1} )T(r)
\end{equation}
where $T=\phi^{-1}\circ j$ is covariant for the functions $\lambda$ of proposition \ref{PropMomslR}.

\end{corollary}
Remark that from general theory of star products, the so-defined $\ast_R$ is a formal star product on $R$.

\begin{proof}
From definition of $\ast_M^R$, on the one hand
\[ 
  (\lambda_X\ast_M^R \lambda_Y)(r)-X\leftrightarrow Y=(\lambda'_X\stM\lambda'_Y)T(r)-X\leftrightarrow Y
=2\nu\{  \lambda'_X,\lambda'_Y \}_{\eR^2}T(r),
\]
while on the other hand, $\omega^R=T^*\Omega$, so that 
\[ 
  \{ \lambda'_X,\lambda'_Y \}_{\eR^2}\circ T=\{ \lambda'_X\circ T,\lambda'_Y\circ T \}_R=\{ \lambda_X,\lambda_Y \}_R.
\]
\end{proof}

All that makes the theory developed earlier, and in particular proposition \ref{prop:def_stn}, valid here. So we pose
\begin{equation}
\begin{aligned}
 {\rho_{\nu}}\colon\sR &\to {\End\big(  C^{\infty}(R)[ [\nu]] \big)} \\ 
  X &\mapsto {\frac{1}{2\nu}\ad_{\ast_M^R}(\lambda_X);} 
\end{aligned}
\end{equation}
using the explicit expressions of $ad_{\stM}(\lambda'_X)$, we find
\begin{align}
  \rho_{\nu}(H)&=\beta\partial_a,	&\rho_{\nu}(E)&=\frac{\beta}{2\nu}e^{-2a}\sinh(2\nu\partial_l).
\end{align}
Using \eqref{eq:dL_et_Poisson} with $\lambda_H=-\beta l$, it is clear that $dL(H)=-\beta\{l,u\}=\beta\partial_a u$. Therefore
\begin{equation}
   \rho_{\nu}(H)=dL(H),
\end{equation}
but the requested identity $\rho_{\nu}(E)=dL(E)$ will not hold. The problem is that $dL(X)=X_{\lambda_X}$ is a vector field, while $\rho_{\nu}(E)$ comes with (infinitely) multiple derivatives, hence this is not a vector field. Conclusion: the operator $\mT$ of equation \eqref{eq:TrnT} must not act on the variable~$a$.

First we consider a partial Fourier transform $\mF$\nomenclature{$\mF$}{Fourier transform}\index{Partial Fourier transform}:
\begin{equation}
 (\mF u)(a,\alpha)=\hu(a,\alpha):=\us{\sqrt{2\pi}}\int e^{-i\alpha l}u(a,l)dl,
\end{equation}
the inverse being given by
\begin{equation}
 (\mF^{-1}\hu)(a,l)=\us{\sqrt{2\pi}}\int e^{il\alpha}\hu(a,\alpha)d\alpha.
\end{equation}
It is clear that $\mF\rho_{\nu}(H)\mF^{-1}=\rho_{\nu}(H)$, but $\mF\rho_{\nu}(E)\mF^{-1}=\frac{\beta}{2\nu}e^{-2a}\sinh(2i\nu\alpha)$. Indeed, if we define $\hv(a,\alpha)=\sinh(2i\nu\alpha)\hu(a,\alpha)$, 
\begin{equation}
\begin{split}
 (\rho_{\nu}(E)\mF^{-1}\hu)(a,l)&=\frac{\beta}{2\nu}e^{-2a}\us{\sqrt{2\pi}}\sinh(2\nu\partial_l)\int e^{il\alpha}\hu(a,\alpha)d\alpha\\
                        &=\frac{\beta}{2\nu}e^{-2a}\us{\sqrt{2\pi}}\int e^{il\alpha} \sinh(2i\alpha\nu)\hu(a,\alpha)d\alpha\\
			&=\frac{\beta}{2\nu}e^{-2a}(\mF^{-1}\hv)(a,l).
\end{split}
\end{equation}
This is nothing but the fact that the Fourier transform turns a derivation into a multiplication.

As can be seen on an asymptotic development, the deformation $\nu$ parameter is necessarily purely imaginary, then we can here pose $\nu=i\theta$ with $\theta\in\eR$, so that
\begin{equation}\label{eq:FrnEF}
   \mF\rho_{\nu}(E)\mF^{-1}=\frac{\beta i}{2\theta}e^{-2a}\sinh(2\alpha\theta).
\end{equation}
Using \eqref{eq:dL_et_Poisson}, we find 
\begin{equation}\label{eq:dLE}
   dL(E)=\beta e^{-2a}\partial_l.
\end{equation}
 Comparing it with the expression of $\mF\rho_{\nu}(E)\mF^{-1}$, we see that (up to constant factor) we have to act in such a way that $\sinh(2\alpha\theta)$ is converted into a derivation. This is done by a Fourier transform. We pose $\xi=\sinh(2\theta\alpha)$ and
\[
  \tf(a,\xi)=\us{\sqrt{2\pi}}\int e^{i\xi p}f(a,p)dp.
\]
As usual, 
\[
   \widetilde{\partial_{\alpha}f}=-i\xi\tf.
\]
This suggests us to consider the change of variable 
\[
  \phi_{\theta}(a,\alpha)=(a,\us{2\theta}\sinh(2\theta\alpha)),
\]
and finally,
\begin{equation}
   \mT_{\theta}:=\mF^{-1}\circ\phi_{\theta}^*\circ\mF,
\end{equation}
where $\phi_{\theta}^*$ is defined by $(\phi_{\theta}^*u)(a,\alpha)=u(a,\us{2\theta}\sinh(2\theta\alpha))$. The result of our construction is the following which proves that we are in the situation of proposition~\ref{prop:def_stn}.
\begin{theorem}
\[
   \mT_{\theta}\circ\rho_{\nu}(E)\circ \mT_{\theta}^{-1}=\beta e^{-2a}\partial_l=dL(E).
\]
\end{theorem}
\begin{proof}
Notice that $(\phi_{\theta}^*\mF u)(a,\alpha)=\hu(a,\sinh(2\theta\alpha))$, and then define $\hv(a,\alpha)=\hu(a,\sinh(2\theta\alpha))$; equation \eqref{eq:FrnEF} is
\[
   (\mF\rho_{\nu}(E)\mF^{-1}\hv)(a,\alpha)=\frac{\beta i}{2\theta}e^{-2a}\sinh(2\theta\alpha)\hv(a,\alpha).
\]
Applying $(\phi^*)^{-1}$ on the right hand side, we find $\frac{\beta i}{2\theta}e^{-2a}2\theta\alpha\hu(a,\alpha)$.  This allows us to compute
\[ 
\begin{split}
  (\mT_{\theta}\rho_{\nu}(E)\mT_{\theta}^{-1})u(a,l)&=\beta ie^{-2a}\mF^{-1}(\alpha\hu)(a,l)\\
                                        &=\beta ie^{-2a}\us{\sqrt{2\pi}}\int \hu(a,\alpha)(-i)\partial_l e^{il\alpha}d\alpha\\
					&=\beta e^{-2a}(\partial_lu)(a,l).
\end{split}
\]
\end{proof}

\subsection{Formula for the product}

The fact the $\mT_{\theta}$ intertwines $\rho_{\nu}$ and $dL$ makes that the candidate to be a product on the $AN$ of $\SL(2,\eR)$ can be computed using formula \eqref{Eq_candprodANSL}. Computations are rather long and done in the articles \cite{StrictSolvableSym} and \cite{Biel-Massar} (see particularly point 4), so we will not give them here. We will also not precise the functional space of convergence for the resulting product.

In the parametrization 
\[ 
  (a,l)=\begin{pmatrix}
 e^{a}		&  e^{a}l\\
0		&	 e^{-a}
\end{pmatrix},
\]
of $R=AN$ the form $da\wedge dl$ is a left-invariant measure, so the integral of the function $f\colon R\to \eR$ on $R$ is given by
\[ 
  \int_R f=\int_{\eR^2} f(a,l)da\,dl.
\]
Remark that $da\,dl$ is the Liouville measure by proposition \ref{prop:Omega}. It is important for definition \ref{DefWKBCompl}.

We consider a subset $\eA\subset\Fun(R)$, and we define the product $\star^{R}_{\theta}$ on $\eA$ by
\begin{equation}\label{eq:star_R}
\begin{split}
(a\star^{R}t b)(a_0,l_0)
=\int_{R\times R}K^R_{\theta}\big((a_0,l_0)&,(a_1,l_1),(a_2,l_2)\big)\\
                                           &a(a_1,l_1)b(a_2,l_2)da_1dl_1da_2dl_2.
\end{split}
\end{equation}
where
\[
K_{\theta}^A(g_0,g_1,g_2)=\us{\theta^2}\mA^R(g_0,g_1,g_2)e^{i\theta \mS^R(g_0,g_1,g_2)}
\]
 with
\begin{subequations}
\begin{align}
  \mA^R(g_0,g_1,g_2)&=\bigoplus_{0,1,2}\cosh(a_1-a_2)\\
  \mS^R(g_0,g_1,g_2)&=\bigoplus_{0,1,2}\sinh(2(a_0-a_1))l_2.
\end{align}
\end{subequations}
Here, the symbol $\bigoplus_{0,1,2}$ stands for a cyclic sum over the indices $0,1,2$.

\subsubsection{Remark on (formal) star product}\label{subsec:rem_on_sp}

One can  find a definition of an asymptotic development for oscillating integrals in \cite{Dieu7} under the form
\[
   I_{\lambda}=\int e^{(i/\lambda)S(x)}\phi(x)\sim\sum_n\lambda^nc_n.
\]
It can be shown that such a  development used on \eqref{eq:star_R} gives rise of a formal star product:
\begin{equation}	\label{EqDevFedFor}
(a\star^{R} b)(g)\sim a(g)b(g)+\frac{\theta}{2i}\{a,b\}(g)+o(\theta^2).
\end{equation}

\section{Extension lemma}		\label{SubSecExtLem}

Let $(\mfs_i,\Omega_i)_{i=1,2}$ be  symplectic Lie algebras and $(S_i,\omega_i)$ the respective Lie groups with left-invariant symplectic forms: $(\omega_i)_g=(L_g)^*\Omega_i$. 
We suppose to know a homomorphism  $\dpt{\rho}{\mfs_1}{\Der(\mfs_2)\cap\mfsp(\Omega_2)}$ and a Darboux chart $\phi_i\colon \mathfrak{s}_i\to S_i$ for each of the two symplectic Lie groups. Our first purpose is to build a Darboux chart on the split extension
\[
   \mfs:=\mfs_1\oplus_{\rho}\mfs_2.
\]

\begin{remark}
Most of the time we are in the inverse situation: we have an algebra $\mathfrak{s}$ which turn out to be a split extension $\mathfrak{s}_{1}\oplus_{\ad}\mathfrak{s}_{2}$ for which we have to check that $\ad(\mathfrak{s}_{1})$ is a symplectic action of $\mathfrak{s}_{1}$ on $(\mathfrak{s}_{2},\Omega_{2})$. See the example of section \ref{SecUnifSOdn}.
\end{remark}
 
\begin{proposition}
In this setting, the map $\dpt{\phi}{\mfs}{S}$,
\begin{equation}
  \phi(X_1,X_2)=\phi_2(X_2)\phi_1(X_1)
\end{equation}
is a Darboux chart.
\label{prop:Darboux}
\end{proposition}

\begin{proof}
An element $X\in T_{\phi^{-1}(g)}(\mfs_1\oplus\mfs_2)=\mathfrak{s}_1\oplus\mathfrak{s}_2$ is denoted by $X=(X_1,X_2)$ with $X_i\in\mathfrak{s}_i$, and the symplectic form on $\mfs_1\oplus\mfs_2$ is given by
\begin{equation}		\label{eq:Omega}
   \Omega\big( (A_1,A_2),(B_1,B_2) \big)=\Omega_1(A_1,B_1)+\Omega_2(A_2,B_2)
\end{equation}
where we identify  $\mfs_i$ and $T_e\mfs_i$.  Let $A$ and $B$ belongs to $T_{\phi^{-1}(g)}(\mathfrak{s}_1\oplus\mathfrak{s}_2)$. We have to show that the quantity
\begin{equation}\label{eq:omega_g_et_omega_e}
\begin{split}
\omega_g\Big(   (d\phi)_{\phi^{-1}(g)}A&,(d\phi)_{\phi^{-1}(g)}B      \Big)\\
                                &=\omega_e\Big(  (dL_{g^{-1}})_g (d\phi)_{\phi^{-1}(g)}A,(dL_{g^{-1}})_g(d\phi)_{\phi^{-1}(g)}B      \Big)
\end{split}  
\end{equation}
does not depend on $g$.

The vector $A$ is represented by a path $A(t)=( A_1(t),A_2(t) )$ with $A_i(t)\in\mfs_i$. In order to characterise that path, we want first to know precisely what is $A_i(0)$. Since $A\in T_{\phi^{-1}(g)}\mfs$, the path must fulfil $\phi( A_1(0),A_2(0) )=g$, or
\begin{equation}
  \phi_2(A_2(0))\phi_1(A_1(0))=g.
\end{equation}
We denote $A_i(0)=G_i\in\mfs_i$ and $\phi_i(G_i)=g_i$. The relation between $g_1$ and $g_2$ is $g_2g_1=g$.  In particular, it is wrong to say ``$A_1(0)=\phi^{-1}(g)$, thus $\phi_1(A_1(0))=g$''. This point being clear, 
\begin{equation}
  (d\phi)_{\phi^{-1}(g)}A=\Dsdd{\phi(A(t))}{t}{0}=\Dsdd{ \phi_2(A_2(t))\phi_1(A_1(t)) }{t}{0}.
\end{equation}
If one particularises to the case $A\in\mfs_2$, that is $A_1(t)=cst=G_1$,
\begin{equation}
  (d\phi)_{\phi^{-1}(g)}A=(dR_{g_1})_{g_2}(d\phi_2)_{G_2}A_2.
\end{equation}

Since $g=g_1g_2$, we have $L_{g^{-1}}=L_{g_1^{-1}}\circ L_{g_2^{-1}}$, and the first argument of $\omega_e$ in equation \eqref{eq:omega_g_et_omega_e} is
\[
   (dL_{g_1^{-1}})_{g_1}(dL_{g_2^{-1}})_{g_2g_1}(dR_{g_1})_{g_2}(d\phi_2)_{G_2}A_2.
\]
If we write that in terms of the derivative of the path $A_2(t)$, what we get in the derivative is
\begin{equation}
g_1^{-1} g_2^{-1} \phi_2(A_2(t))g_1=\AD_{g_1^{-1}}\Big( g_2^{-1}\phi_2(A_2(t)) \Big).
\end{equation}
Since $g^{-1}_2\phi\big( A_2(0) \big)=g^{-1}_2\phi_2(G_2)=e$, the derivative of that term is
\begin{equation}
  (dL_{g^{-1}})_g(d\phi)_{\phi^{-1}(g)}A=\Ad_{g_1^{-1}}\Big(  (dL_{g_2^{-1}})_{g_2}(d\phi_2)_{G_2}A  \Big)
\end{equation}
with some abuse between $A\in\mfs$ and $A_2\in\mfs_2$.  Doing the same computation with $B\in\mfs_1$ (so that $B_2(t)=cst=G_2$), we find
\begin{equation}
   (d\phi)_{\phi^{-1}(g)}B=\Dsdd{ g_2\phi_1(B_1(t)) }{t}{0}
                         =(dL_{g_2})_{g_1}(d\phi_1)_{G_1}B_1,
\end{equation}
 and what appears in $\omega_e$ reads 
\begin{equation}
  (dL_{g^{-1}})_g(dL_{g_2})_{g_1}(d\phi_1)_{G_1}B=(dL_{g_1^{-1}})_{g_1}(d\phi_1)_{G_1}B.
\end{equation}
Finally, for $A\in\mfs_2$ and $B\in\mfs_1$,
\begin{equation}\label{eq:gros_omega_e}
\begin{split}
\omega_g\big( (d\phi)_{\phi^{-1}(g)}A&,(d\phi)_{\phi^{-1}(g)}B \big)\\
                                    &= \omega_e\Big(
       \Ad_{g_1^{-1}}\big[  (dL_{g_2^{-1}})_{g_2}(d\phi_2)_{G_2}A  \big],
       (dL_{g_1^{-1}})_{g_1}(d\phi_1)_{G_1}B
         \Big). 
\end{split} 
\end{equation}
The first argument belongs to $T\mfs_2$ (because $g_2\in\mfs_2$) while the second belongs to $T\mfs_1$. Hence definition \eqref{eq:Omega} makes the right hand side vanishing.

If we want to compute equation \eqref{eq:gros_omega_e} with $A$, $B\in\mfs_2$,
\begin{equation}
\begin{split}
\omega_e\Big(  
           \Ad_{g_1^{-1}}\big[&  (dL_{g_2^{-1}})_{g_2}(d\phi_2)_{G_2}A  \big],
	   \Ad_{g_1^{-1}}\big[  (dL_{g_2^{-1}})_{g_2}(d\phi_2)_{G_2}B  \big]   
        \Big)\\
  &=\Omega(\ldots)=\underbrace{\Omega_1(\ldots)}_{=0}+\Omega_2(\ldots)\\
  &=\Big( \Ad_{g_1^{-1}}^*\Omega_2  \Big)
     \Big(
        (dL_{g_2^{-1}})_{g_2}(d\phi_2)_{G_2}A,\ldots B
     \Big)
\end{split}
\end{equation}
At this point, notice that $\Ad^*_{g_1}\Omega_2=\Omega_2$. Indeed the exponential $\exp\colon \mathfrak{s}_1\to S_1$ being surjective , there exists a $X_1\in\mathfrak{s}_1$ such that $\Ad(g_1)= e^{\ad(X_1)}$. Now, $\ad(X_1)\in\gsp(\Omega_2)$ by assumption, so that $\Ad(g_1)\in\SP(\Omega_2)$. The previous expression becomes
\begin{equation}
\begin{split}
  \Omega_2( (dL_{g_2^{-1}})_{g_2}(d\phi_2)_{G_2}A,\ldots B )
  &=(\omega_2)_{g_2}( (d\phi_2)_{g_2}A,\ldots B )\\
  &=(\phi_2^*\omega_2)_{G_2}(A,B)\\
  &=\Omega_2(A,B).
\end{split}  
\end{equation}
The last line is the fact that $\phi_2$ is a Darboux chart: $ \phi_2^*\omega_2=\Omega_2$. The case with $A$, $B\in\mfs_1$ yields to compute
\[
   \omega_e\Big(
                   (dL_{g_1^{-1}})_{g_1}(d\phi_1)_{G_1}A,  (dL_{g_1^{-1}})_{g_1}(d\phi_1)_{G_1}B
           \Big).
\]
It is done by the same way as the previous cases.
\end{proof}

A direct computation shows the following \defe{extension lemma}{Extension!lemma}.
\begin{lemma}[Extension lemma]
Let $K_i\in \Fun(S_i^3) $ be a left-invariant three point kernel on $S_i$ ($i=1,2$).  Assume that $K_2\otimes1\in \Fun(S^3)$ is invariant under conjugation by elements of $S_1$.  Then $K:=K_1\otimes K_2\in\Fun(S^3)$ is left-invariant (under $S$).
\label{EXT}
\end{lemma}
 
\begin{proof}
An element of $S$ has the form $g_1g_2$ with $g_i\in S_i$, the multiplication being given by $(g_1g_2)(a_1a_2)=(g_1a_1)(g_2a_2)$. Using this rule, the definition of the tensor product, and the left-invariance of both $K_i$,
\[ 
\begin{split}
\big( L_{g_1g_2}(K_1\otimes K_2) \big)&(a_1a_2,b_1b_2,c_1c_2)\\
			&=(K_1\otimes K_2)\big( (g_1g_2)(a_1a_2),(g_1g_2)(b_1b_2),(g_1,g_2)(c_1c_2) \big)\\
			&=K_1(g_1a_1,g_1b_1,g_1c_1)K_2(g_2a_2,g_2b_2,g_2c_2)\\
			&=K_1(a_1,b_1,c_1)K_2(a_2,b_2,c_2)\\
			&=(K_1\otimes K_2)(a_1a_2,b_1b_2,c_1c_2).
\end{split}  
\]

\end{proof}

This lemma shows that if one has kernels on $S_1$ and $S_2$ satisfying the above hypotheses, their tensor product provides a kernel for an associative left-invariant kernel on $S=S_1\otimes_{\rho} S_2$.  Proposition \ref{prop:Darboux} allows us to hope that the product on $S$ will satisfy the same kind of symplectic compatibility as the products on $S_i$; in particular when the latter were constructed using Darboux chart in the same way as the product described in section  \ref{sec:unifsl}.

\section{Deformation by group action}  	\label{SecDefAction}		

The procedure of deformation by group action is described in \cite{TrsStProd}. Let $G$ be a Lie group. We suppose to know a subset $A^G$ of $\Fun(G,\eC)$  such that
\begin{enumerate}
\item $A^G$ is invariant under the left regular representation of $G$ on itself,
\item $A^G$ is provided with a $G$-invariant product $\stG$\nomenclature{$\stG$}{Star product on $G$} such that $(A^G,\stG)$ is an associative algebra. The $G$-invariance means that $\forall a,b\in A^G$,
\[
    (L^*_ga)\stG(L^*_gb)=L^*_g(a\stG b).
\]
\end{enumerate}

Notice that we do not impose any regularity condition on this product. The reason is that the deformation by group action is a formal procedure which allows to guess a product on a manifold. The ``true'' work to prove convergences and invariances has to be done on the level of the deformed manifold.

Now, let $X$ be a manifold endowed with a right action  $\dpt{\tau}{G\times X}{X}$ of $G$.  For $u\in\Fun(X)$, $x\in X$ and $g\in G$, we consider $\alpha^x(u)\in\Fun(G)$ and $\alpha_g(u)\in\Fun(X)$ defined by
\begin{equation}		\label{EqDefalphaxu}
   \alpha^x(u)(g)=\alpha_g(u)(x)=u(\tau_{g^{-1}}(x)),
\end{equation}
and the following functional space on $X$:
\[
   A^X=\{u\in\Fun(X)|\alpha^x(u)\in A^G\,\forall x\in X\}.
\]
For example, the $A^X$ corresponding to $A^G=\Fun(G)$ is the whole $\Fun(X)$.  For $u$, $v\in A^X$, we define $\dpt{\stX}{A^X\times A^X}{\Fun(X)}$\nomenclature{$\stX$}{Star product on $X$}
\begin{equation}\label{eq:def_stG}
   (u\stX v)(x)=\big(\alpha^x(u)\stG \alpha^x(v)\big)(e)
\end{equation}
where $e$ is the identity of $G$.

\begin{theorem}
The product \eqref{eq:def_stG} obtained by action of the group $G$ on the manifold $X$ fulfils the following properties:
\begin{enumerate}
\item\label{itemthostG} The operation $\alpha^x$ intertwines the products $\star^X$ and $\star^G$:
\[
   \alpha^x(u\stX v)=(\alpha^xu)\stG(\alpha^xv).
\]

\item $A^X$ is stable under $\stX$,
\item $(A^X,\stX)$ is an associative algebra.
\end{enumerate}
\end{theorem}

\begin{proof}
First remark that $\alpha^{\tau_{g^{-1}}(x)}u=L^*_g\alpha^xu$ because
\[ 
(\alpha^{\tau_{g^{-1}}(x)}u)(h)=u\big(\tau_{(gh)^{-1}}(x)\big)=(\alpha^xu)(gh)
=\big(L^*_g\alpha^xu\big)(h),
\]
It follows that
\begin{equation}
\begin{split}
\alpha^x(u\stX v)(g)&=(u\stX v)(\tau_{g^{-1}}(x))
                    =(\alpha^{\tau_{g^{-1}}(x)}u\stG\alpha^{\tau_{g^{-1}}(x)}v) (e)\\
		    &=\left[ L^*_g(\alpha^xu\stG\alpha^xv) \right](e)
		    =(\alpha^xu\stG\alpha^xv)(g).
\end{split}
\end{equation}
The first point is proved.

Using the first point, we see that $\alpha^{\tau_{g^{-1}}(x)}u$ belongs to $A^G$ because $\alpha^xu\in G^G$ and $A^G$ is stable under $L_g$. So we have the second point.  For the third one,
\[ 
\begin{split}
[(u\stX v)\stX w](x)&=\big(\alpha^x(u\stX v)\stG\alpha^x(w)\big)(e)\\
                    &=\big((\alpha^xu\stG\alpha^xv)\stG\alpha^x(w)\big)(e).
\end{split}
\]
The conclusion follows from associativity of $\star^G$.
\end{proof}

Let us summarize what was done up to now. When $G$ acts on $X$, and when we have a ``good'' product on $A^G\subset\Fun(G)$, we are able to build an associative product on $A^X\subset\Fun(X)$. The space $A^X$ is defined by $A$ and the action. So a deformation of a group gives rise to a deformation of any manifold on which the group acts. This is why we call it an \oge universal\fge\ deformation. That universal construction is the motivation to deform groups.

\begin{lemma}   
A function $u$ belongs to $A^{X}$ if and only if there exists one $y$ such that $\alpha^y(u)\in A^G$ in each $g$-orbit in $X$.
		\label{LemUnPtParOrbite}
\end{lemma}

\begin{proof}
The necessary condition is direct because, when $u\in A^X$, the function $\alpha^x(u)$ belongs to $A^G$ for every $x$. For the sufficient condition, suppose $\alpha^y(u)\in A^{G}$, then $\alpha^{g\cdot y}(u)=L_{g^{-1}}^*(\alpha^yu)\in A^{G}$ for all $g$ because $A^{G}$ is left-invariant. If it holds for a $y$ in each $G$-orbit, then $\alpha^xu\in A^{G}$ for all $x\in X$.
\end{proof}

The content of this lemma is that if one wants to check if a given function $u$ belongs to $A^{X}$, one only has to check is $\alpha^yu\in A^{G}$ for one $y$ in each $G$-orbit.

The functions $\alpha^x(u)$ are not ``gentle'' functions, even when $u$ is. Let us give two examples of pathology that can occur in $\alpha^x(u)$ without to be present in $u$. Firstly,  if the action is the identity, the support of $\alpha^x(u)$ is the whole $G$ which can be non compact. So, even when $u$ is compactly supported, there are no guarantee with respect to the support of $\alpha^x(u)$. 

Secondly, the function $\alpha^x(u)$ is of course bounded; but the derivatives are not specially such. Indeed, in order to fix ideas, suppose that the group $G$ is a two parameter group and that the manifold $X$ is a two dimensional manifold. In this case, one can write
\begin{equation}		\label{EqDefziDefA}
  f(a,l)=\alpha^x(u)(a,l)=u\big( z_1(a,l),z_2(a,l) \big)
\end{equation}
where $x$ is a parameter in the functions $z_i$. Depending on the action, the function $z$ can be very odd. In particular, the derivatives
\[ 
  (\partial_af)(a,l)=(\partial_1u)(z_1,z_2)(\partial_az_1)(a,l)+(\partial_2u)(z_1,z_2)(\partial_az_2)(a,l)
\]
in which $\partial_az_i$ can be divergent. Even worse, the degree of the divergence can increase with the degree of the derivation. Two examples of such a hill behaviour are given in section \ref{SecEplolUnter}. 
   \section[Split extensions of Heisenberg algebras]{One dimensional split extensions of Heisenberg algebras} \label{SecExtHeiz}

\subsection{Introduction}

The one dimensional extensions of Heisenberg algebras are classified by triples $(\BX,\mu,d)$. The quantization in the case $(\id,0,\mu)$  reveals to be a particular case of the one studied in \cite{Biel-Massar}, while quantization of other extensions can be found using symmetries of the kernel. Here we are reporting results of \cite{articleBVCS} and most of proofs (in particular the trick of subsection \ref{subsecTrick} which allows to extend the known product to every one dimensional split extensions) are due to Y. Voglaire. It is to be published in his future PhD thesis.

The kernel of the quantization of \cite{Biel-Massar} will be denoted by $K$.  Then we will give a way to twist $K$ in order to obtain a kernel $K'$ on any extension of the form $(\BX,0,2)$. Quantizations of other extensions can be obtained by composing with Lie group isomorphisms. The kernel for an arbitrary extension is denoted by $K_{0}(\BX,\mu,d)$, or simply $K_{0}$ when there are no possible ambiguity.

When we will deal with the anti de Sitter situation, our starting point will be this $K_{0}$ that we will have to adapt to another symplectic form that $\delta E^*$ invoking lemma \ref{LemJumpCoadOrb}.

\subsection{General definitions}

Let $\pH_{n}=V\oplus\eR E$ be the \defe{Heisenberg algebra}{Heisenberg algebra} of dimension $2n+1$, with a natural symplectic structure defined from the Heisenberg algebra structure:
 \[
[v,w]=\Omega(v,w)E
\]
for all $v$, $w\in V$. Now we consider a one dimensional algebra $\mA=\eR A$ generated by an element $A$, and we build the split extension of $\pH_{n}$ by $\mA$:
\begin{equation}
\mF(\rho)=\mA\oplus_{\rho}\pH_{n}
\end{equation}
where the split homomorphism is an action by derivation $\rho\colon \mA\to \Der(\pH_{n})$. The so obtained algebra is what we call a \defe{one dimensional extension of Heisenberg algebra}{Extension!of Heisenberg algebra}. Let us study the possibilities for $\rho(A)$. From linearity, its general form is 
\[ 
 \rho(A)(v,z)=\rho(A)(v,0)+\rho(A)(0,z)
	=(\BX v,\mu(v))+(zv_{0},2dz)
\]
with $\BX \in\End(V)$, $\mu\in V^*$, $v_{0}\in V$ and $d\in\eR$. Since $\eR E=[\pH_{n},\pH_{n}]$, the fact that $\rho(A)$ is a derivation of $\pH_{n}$, implies that $v_{0}=0$ because
\begin{equation}
   \rho(A)\eR E=\rho(A)[\pH_{n},\pH_{n}]
		=[\rho(A)\pH_{n},\pH_{n}]+[\pH_{n},\rho(A)\pH_{n}]\subset \eR E.
\end{equation}
Thus we have
\begin{equation}    \label{EqrhoBmudz}
\rho(A)(v,z)=(\BX v,\mu(v)+2dz).
\end{equation}
From commutation relations in $\pH_{n}$, we easily find 
\[ 
  [(v,z),(v',z')]=[v,v']=\Omega(v,v')E.
\]
Applying $\rho(A)$ to this equality, and using the fact that this is a derivation, we find
\[ 
  \Omega(\BX ,v')E+\Omega(v,\BX ')E=\rho(A)\Omega(v,v')E=2d\Omega(v,v')E
\]
which can be rewritten as
\begin{equation}
\Omega\big( (\BX -d\,\mtu)v,v' \big)+\Omega\big( v,(\BX -d\,\mtu)v' \big)=0.
\end{equation}
In conclusion, the endomorphism $\rho(A)$ is given by a triple $(\BX ,\mu,2d)$ with $(\BX -d\,\mtu)\in\gsp(V,\Omega)$, $\mu\in V^*$ and $d\in\eR$. Using this result, we write the general commutator on $\mR=\mA\oplus_{\rho}\pH_{n}$ under the form
\begin{equation}  \label{EqGeneExtHeizCom}
\big[ (a,v,z),(a',v',z') \big]=\big( 0,\BX (av'-a'v),\mu(av'-a'v)+2d(az'-a'z)+\Omega(v,v') \big)
\end{equation}
where we adopted the notation
\begin{equation} 
(a,v,z)=aA+v+zE.
\end{equation}

\subsection{Symplectic structure}

The following proposition gives a symplectic structure on $\mF$.

\begin{proposition}   
The algebra $(\BX,\mu,d)$ endowed with
\begin{equation}
 \Omega^{\mF}=-\delta E^*=E^*([.,.])
 \end{equation}
where the star denotes the Chevalley cocycle defined by \eqref{EqDefChevCoycl} is symplectic if and only if $d\neq 0$.
\label{PropSymplestarEG}
\end{proposition}

\begin{proof}
It is evident that $\Omega^{\mF}$ is closed because it is exact. For non-degeneracy, we compute
\[ 
\begin{split}
\Omega^{\mF}&=E^*[.,.]=a\mu(v')-a'\mu(v)+2d(az'-a'z)+\Omega(v,v')\\
	&=\begin{pmatrix}
0	&	\mu	&	2d\\
-\mu^{t}&	\Omega	&	0\\
-2d	&	0	&	0	
\end{pmatrix}
\end{split}  
\]
whose determinant is $\det\Omega^{\mF}=-4d^{2}\det\Omega$ which is non vanishing if and only if $d\neq 0$.

\end{proof}

This symplectic algebra  is denoted by $\mF_{\Omega}(\BX ,\mu,d)$, or simply $\mF$ when there are no possible confusions.

Since we are only interested in symplectic algebras, we suppose $d\neq 0$ and we look at extensions of type $(d\BX ,d\mu,2d)$ with $\BX -d\mtu\in\gsp(V,\Omega)$. The bracket is given by
\begin{equation}   \label{EqCommGeneF}
\big[ (a,v,z),(a',v',z') \big]=\big( 0,d\BX (a'v-a'v),d\mu(av'-a'v)+2d(az'-a'z)+\Omega(v,v') \big).
\end{equation}

\subsection{Isomorphisms}  \label{SubsecIsomsdX}

The extension obtained by the derivation $D=(\BX ,\mu,d)$ is \emph{a priori} not the same as the one obtained by $D'=(\BX',\mu',d')$. Two extensions are isomorphic when there exists a linear bijection $dL\colon \mF_{D}\to \mF_{D'}$ such that\footnote{the reason why we write $dL$ instead of $L$ comes from the fact that we will be interested in the corresponding group isomorphism later.}
\begin{subequations}
\begin{align}
dL\big( [X,Y]_{D'} \big)&=\big[ dL(X),dL(Y) \big]_{D'}\\
(dL)^*\Omega^{D'}&=\Omega^{D}.
\end{align} 
\end{subequations}
We find the following isomorphisms:
\begin{subequations}   \label{SubEqsIsommud}
\begin{itemize}
\item $\mF(d\BX ,d\mu,2d)\simeq \mF(\BX ,\mu,d)$ by
\begin{equation}
   dL(a,v,z)=(da,v,z),
\end{equation}
\item $\mF(\BX ,\mu,2)\simeq\mF(\BX ,0,2) $ by
\begin{equation}
   dL(a,v,z)=(a,v+au,z),
\end{equation}
where $u$ is the vector of $V$ satisfying $i(u)\Omega=\mu$,
\item $\mF( \BX ,0,2)\simeq\mF(\BX ',0,2)$ by
\begin{equation}
   dL(a,v,z)=(a,M(v),z)
\end{equation}
where $M\in\SP(V,\Omega)$ fulfills $M\BX M^{-1}=\BX'$ or, equivalently,
\[
M(\BX-\mtu)M^{-1}=\BX'-\mtu.
\]
\end{itemize}
\end{subequations}
The third isomorphism only gives the equivalence between $\BX-\mtu$ and $\BX'-\mtu$ when they belongs to the same orbit of the adjoint action of $\SP(V,\Omega)$. In particular, there are no isomorphisms between the identity and anything else.

\subsection{Reminder about a previous deformation}   \label{SubsecDefSURme}

Before going on with the construction of a deformation of one dimensional split extensions of Heisenberg algebras, we have to recall a result on deformation in $\SU(1,n)$. The product on the extension of Heisenberg algebra will be nothing else than a transport of this one.

The article \cite{Biel-Massar} provides a formal universal deformation formula for the actions of the Iwasawa component $\SUR_0:=\SUA_0\SUN_0$ of $\SU(1,n)$ under an oscillatory integral form.  It turns out (see \cite{lcBBM}) that this deformation formula is in fact non-formal for proper actions on topological spaces. 

Here is the precise result. The Iwasawa decomposition of $\SU(1,n)$ induces the identification $\SUR_0=\SU(1,n)/U(n)$. The group $\SUR_0$ is endowed with a (family of) left-invariant symplectic structure(s)\footnote{This is done using the hermitian symmetric structure, cf proposition 1.1 in \cite{Biel-Massar}.} $\omega$.  If we denote by $\sR_0=\sA_0\oplus\sN_0$ the Lie algebra of $\SUR_0$, the map
\begin{equation}  \label{DARBOUX}
\begin{aligned}
 \phi_0\colon \sR_0&\to \SUR_0 \\ 
(a,n)&\mapsto \exp(a)\exp(n) 
\end{aligned}
\end{equation}
reveals to be a global Darboux chart for $(\SUR_0,\omega)$.  The nilpotent component appears to accept a decomposition $\sN_0=V\times\eR Z$ in which the Lie bracket reads
\[ 
[(x,z)\,,\,(x',z')]=\Omega_V(x,x')\,Z; 
\]
the full Iwasawa component is now parametrized by $\sR_0=\{(a,v,z)\,|\,,a,z\in\eR;x\in V\}$. The interest of this situation resides in the fact that the algebra $\sR_0$ turns out to be a one dimensional split extension of an Heisenberg algebra; namely, 
\[ 
\sR_0=\mF(\mtu,0,2).
\]
The deformation result is the following.

\begin{theorem}
For every non-zero $\theta\in\eR$, there exists a Fréchet function space $\swE_\theta$ satisfying the inclusions $C^\infty_c(\SUR_0)\subset\swE_\theta\subset C^\infty(\SUR_0)$, such that, defining for all $u,v\in C^\infty_c(\SUR_0)$  
\begin{equation}  \label{PRODUCT}
\begin{split}
(u\star_\theta v)(a_0,x_0,z_0)
		:=\frac{1}{\theta^{\dim \SUR_0}} \int_{ \SUR_0\times \SUR_0}& \cosh(2(a_1-a_2))\\
		&[\cosh(a_2- a_0)\cosh(a_0-a_1)\,]^{\dim \SUR_0-2}\\
&\times \exp \Big( \frac{2i}{\theta}\varphi(r_1,r_2,r_3)\Big)\\
		&\times u(a_1,x_1,z_1)\,v(a_2,x_2,z_2)\, da\,dx\,dz;
\end{split}
\end{equation}
where 
\[ 
\begin{split}
  \varphi(r_1,r_2,r_3)=&S_V\big(\cosh(a_1-a_2)x_0, \cosh(a_2-a_0)x_1, \cosh(a_0-a_1)x_2\big)\\
			&-\bigoplus_{0,1,2}\sinh(2(a_0-a_1))z_2 
\end{split}
\]
with $S_V(x_0,x_1,x_2):=\Omega_V(x_0,x_1)+\Omega_V(x_1,x_2)+\Omega_V(x_2,x_0)$ is the phase for the Weyl product on $C^\infty_c(V)$ and $\bigoplus_{0,1,2}$ stands for cyclic summation, one has: 

\begin{enumerate}

\item\label{tBMi} 
	 $u\star_\theta v$ is smooth and the map $ C^\infty_c(\SUR_0)\times C^\infty_c(\SUR_0) \to C^\infty(\SUR_0)$ extends to an associative product on $\swE_\theta$. The pair $(\swE_\theta,\star_\theta)$ is a (pre-$C^\star$) Fréchet algebra.

\item\label{tBMii}
	 In coordinates $(a,x,z)$ the group multiplication law reads
\[ 
	L_{(a,x,z)}(a',x',z')=\left(a+a',e^{-a'}x+x',e^{-2a'}z+z'+\frac{1}{2}\Omega_V(x,x')e^{-a'}
\right).
\]
The phase and amplitude occurring in formula \eqref{PRODUCT} are both invariant under the left action $L:\SUR_0\times \SUR_0\to \SUR_0$.

\item\label{tBMiii}
	 Formula \eqref{PRODUCT} admits a formal asymptotic expansion of the form:
 \begin{equation*}
	u\star_\theta v\sim \,uv\,+\,\frac{\theta}{2i}\{u,v\}\,+O(\theta^2); 
\end{equation*}
where $\{\,,\,\}$ denotes the symplectic Poisson bracket on $C^\infty(\SUR_0)$ associated with $\omega$.  The full series yields an associative formal star product on $(\SUR_0,\omega)$ denoted by $\tilde{\star}_\theta$. 
 \end{enumerate}
\label{ThoDefHeizsansB}
\end{theorem}

The setting and \ref{tBMi} and \ref{tBMii} may be found in \cite{Biel-Massar}, while \ref{tBMiii} is a straightforward adaptation  to $\SUR_0$ of \cite{lcBBM}.

This theorem among with the isomorphisms given in \ref{SubsecIsomsdX} only provide a product on extensions of type $(d\mtu,0,2d)$. But we saw that the extensions $(\BX ,0,2d)$ with $\BX \neq\mtu$ are different. Hence the generalization of this result to other extensions is not straightforward. We address now this question.

\subsection{Extensions with non trivial \texorpdfstring{$\protect\BX $}{X}}	\label{subsecTrick}

The group $F(\mtu,0,2)$ is provided with a kernel $K\colon F\times F\times F\to \eC$ by theorem \ref{ThoDefHeizsansB}.  The symplectic group $\SP(V,\Omega)$ acts on $F$ by
\begin{equation}
\begin{aligned}
 \Phi\colon \SP(V,\Omega)\times F&\to F \\ 
(M,I(a,v,z))&\mapsto \Phi_{M}(I(a,v,z)):=I(a,M(v),z)
\end{aligned}
\end{equation}
where 
\begin{equation}
\begin{aligned}
 I\colon \mF&\to F \\ 
(a,n)&\mapsto  e^{aA} e^{n} 
\end{aligned}
\end{equation}
is the Iwasawa coordinate on $F$.

\begin{proposition}  
The kernel $K$ is invariant under this action: $\Phi^*_{M}K=K$.
 \label{PropkernelinvarSp}
\end{proposition}

\begin{proof}
We are looking on the kernel in expression \eqref{PRODUCT}. The amplitude of $K$, i.e. all what lies outside the exponential, and the cyclic sum in the phase only depend on the $a_{i}$'s. So $\Phi_{M}$ does not act on them. As far as $S_{0}$ is concerned, up to coefficients which only depend on the $a_{i}$'s, it is a sum of elements of the form $\Omega(Mv_{i},Mv_{j})=\Omega(v_{i},v_{j})$.
\end{proof}

Let $\BX $ be a matrix such that $\bar{\BX }=\BX -\mtu\in\gsp(V,\Omega)$ and $\mF'=\mF'(\BX,0,2)$. We consider $\mS$, the one dimensional subalgebra of $\gsp(V,\Omega)$ generated by $\bar{\BX}$ and we define
\begin{equation}
  \mG=\mS\oplus_{\rho}\mF
\end{equation} 
with 
\[ 
  \rho(\bar{\BX })(a,v,z)=(0,\bar{\BX }v,0).
\]
We denote by $G$ and $S$ the corresponding groups. We have in particular $F\simeq G/S$. An element of $\mG$ has the form
\begin{equation}	\label{EqDefkavzG}
  (k\bar\BX,a,v,z)=k\bar\BX+aA+v+zE.
\end{equation}

\begin{proposition}
The group $F'$ is a subgroup of $G$.
\end{proposition}
\begin{proof}
We will prove that $\mF'$ is isomorphic to a subalgebra of $\mG$, namely, the subalgebra $\mL\subset \mS\oplus_{\rho}\mF$,
\[ 
  \mL=\eR(A+\bar\BX)\oplus_{\sigma}(V+\eR E)
\]
where $\sigma$ is the splitting homomorphism \eqref{EqrhoBmudz} of $\mF$, which in the present case reads $\sigma(A+\bar\BX)(0,v,z)=(0,Xv,2z)$. In other words, the algebra $\mL$ is made of elements of the form \eqref{EqDefkavzG} with $k=a$.  The isomorphism is 
\begin{equation}
\begin{aligned}
\phi  \colon \mL&\to \mF'(\BX,0,2) \\ 
a(A+\bar\BX)+v+zE&\mapsto aA+v+zE. 
\end{aligned}
\end{equation}
Indeed, using formula \eqref{EqCommGeneF} with $d=1$ and $\mu=0$, we find
\[ 
\begin{split}
\Big[ \phi\big( a(A+\bar{\BX })+v+zE \big)&,\phi\big( a'(A+\bar{\BX })+v'+z'E \big) \Big]_{(\BX,0,2)}\\
		&=\big[ aA+v+zE,a'A,v'+z'E \big]\\
	 	&=\BX (av'-a'v)+\big(2(az'-a'z)+\Omega(v,v')\big)E\\
	&=\phi\big( X(av'-a'v),2(az'-a'z+\Omega(v,v'))  \big)\\
&=\phi\big[ a(A+\bar{\BX })+v+zE,a'(A+\bar{\BX })+v'+z'E \big].
\end{split}  
\]
\end{proof}

From now on, we identify $\mF'$ with $\mL$ by the isomorphism $\phi$ which will no longer be explicitly written. Image of $F'$ in $G$ by the isomorphism are elements of the form
\[ 
  g'= e^{a(A+\bar{\BX })} e^{v+zE}.
\]
Since the elements $ e^{A}$ and $ e^{\BX }$ commute in $G$, we can decompose an element $\phi^{-1}(g')$ as
\[ 
  \underbrace{e^{a\bar{\BX }}}_{\in S}\underbrace{e^{aA} e^{v+zE}}_{\in F}.
\]
The element $a(A+\bar{\BX })+v+zE$ seen in $\mS\oplus_{\rho}\mF$ will be denoted by $(a,v,z)$ as well
\[ 
  (a,v,z)=\phi^{-1}(aA+v+zE).
\]
We consider the following coordinate on $F'$:
\begin{equation} 
\begin{aligned}
 J\colon \mF'&\to F' \\ 
(a,v,z)&\mapsto  e^{a(A+\bar{\BX })} e^{v+zE}.
\end{aligned}
\end{equation}

\begin{proposition}
The group $F'$ is diffeomorphic to the homogeneous space $F\simeq G/S$.
\end{proposition}

\begin{proof}
We will prove that $F'$ acts simply transitively on $G/S$. Let us look at 
\begin{equation}		\label{EqgpJSF}
  g'=J(a,v,z)= \underbrace{e^{a(A+\bar{\BX })}}_{g'_{S}} \underbrace{e^{v+zE}}_{g'_{F}}.
\end{equation}
Noticing that $ e^{a\bar{\BX }}[e]=[ e^{a\bar{\BX }}]=[e]$ we find
\[ 
\begin{split}
g'[e]= e^{a\bar{\BX }} e^{aA} e^{v+zE} e^{-a\bar{\BX }} e^{a\bar{\BX }}[e]
		=\AD( e^{a\bar{\BX }})\big(  e^{aA} e^{v+zE} \big)[e].
\end{split}  
\]
In $\mG=\mS\oplus_{\rho}\mF$, by definition of $\rho$, we 
have 
$\AD( e^{a\bar{\BX }})\big(  e^{aA} e^{v+zE} \big)= e^{aA} e^{ e^{a\bar{\BX }}v+zE}$,
thus $g'[e]= e^{aA} e^{ e^{a\bar{\BX }}v+zE}[e]=[I(a, e^{a\bar{\BX }}v,z)]$. So, in order to get the element $[I(a,v,z)]\in G/S$, we have to act on $[e]$ with the element $g'=J(a, e^{-a\bar{\BX }}v,z)$. All that proves that the map
\begin{equation}
\begin{aligned}
 H\colon F'&\to G/S \\ 
(a,v,z)&\mapsto \big[ I(a, e^{a\bar{\BX }}v,z) \big]  
\end{aligned}
\end{equation}
is a diffeomorphism.
\end{proof}

The work done up to now provides a diffeomorphism
\begin{equation}
\begin{aligned}
 \varphi\colon F'&\to F \\ 
\varphi\big( J(a,v,z) \big)&=I(a, e^{a\bar{\BX }}v,z)
\end{aligned}
\end{equation}
which has suitable properties listed in the proposition below.

\begin{proposition}
This map $\varphi\colon F'\to F$ has the following properties:
\begin{enumerate}
\item if $g'=J(a,v,z)=g'_{S}g'_{F}$ in the sense of decomposition \eqref{EqgpJSF},
\begin{equation}
\varphi\circ L_{g'}=\AD(g'_{S})\circ L_{g'_{F}}\circ\varphi=\Phi_{ e^{a\bar{\BX }}}\circ L_{g'_{F}}\circ\varphi=\Phi_{g'_{S}}\circ L_{g'_{F}}\circ\varphi,
\end{equation}
\item the differential fulfils
\begin{equation}
d(\varphi\circ J)_{(0,0,0)}=dI_{(0,0,0)},
\end{equation}
\item if $\omega$ is the left-invariant symplectic form on $F$ and $\omega'$ the one on $F'$, we have
\[ 
  \varphi^*\omega=\omega',
\]
in other words, $\varphi$ is a symplectomorphism.

\end{enumerate}

\end{proposition}
\begin{proof}

The first point is a computation:
\[ 
\begin{split}
\varphi\big( L_{g'_Sg'_F}(g_sg_E) \big)&=\varphi\big( g'_Sg_S\AD(g_S^{-1})(g'_F)g_F \big)\\
		&=\AD(g'_Sg_S)\big( \AD(g_S^{-1})(g'_F)g_F \big)\\
		&=\AD(g'_S)\big( g'_F\AD(g_S)(g_F) \big)\\
		&=\big( \AD(g'_S)\circ L_{g'_F} \big)\big( \varphi(g_Sg_F) \big).
\end{split}  
\]
When $g'=g'_Sg'_F=J(a,v,z)$, we have $g'_S=\exp(a\BX)$ and $g'_F=I(a,v,z)$, so the result is given by
\[ 
  \AD(g'_S)(g'_F)= e^{ \Ad(a\BX) }I(a,v,z)=I(a, e^{a\BX}v,z)=\Phi_{ e^{a\BX}}I(a,v,z).
\]
That concludes the proof of the first point.  For the second statement, we have $(\varphi\circ J)(a,v,z)=\Phi_{ e^{a\BX}}I(a,v,z)$, so
\begin{equation}
\begin{split}
	d(\varphi \circ J)_{(0,0,0)}(Y_a,Y_v,Y_z)&=\Dsdd{  \Phi_{ e^{tY_a}}I(tY_a,tY_v,tY_z)   }{t}{0}\\
		&=\Dsdd{ I\big( tY_a, e^{tY_a\BX}tY_v,tY_z \big) }{t}{0}\\
		&=dI_{(0,0,0)}(Y_a,Y_v,Y_z).
\end{split}
\end{equation}
For the third point, we denote by $e$ and $e'$ the neutral of $F$ and $F'$. On the one hand,
\[ 
  (\varphi^*\omega)_{g'}=\omega_{\varphi(g')}\circ d\varphi_{g'}=\omega_e\circ d\big( L_{\varphi(g')^{-1}}\circ\varphi \big)_{g'};
\]
on the other hand, $\omega'_{g'}=\omega'_{e'}\circ d\big(L_{(g')^{-1}})_{g'}$. Hence, in order to have $\varphi^*\omega=\omega'$, it is necessary that
\[ 
  \omega'_{e'}\circ dJ_{(0,0,0)}=\omega_e\circ d\big(  L_{\varphi(g')^{-1}}\circ\varphi\circ L_{g'}  \big)_{e'}\circ dJ_{(0,0,0)}.
\]
But, for $g'=g'_Sg'_F$, we have
\[ 
\begin{split}
L_{\varphi(g')^{-1}}\circ\varphi\circ L_{g'}(g)&=\varphi(g')^{-1}\varphi(g'g)\\
		&=\varphi(g')^{-1}\AD(g'_S)\big( g'_F\varphi(g) \big)\\
		&=\AD\big( (g'_F)^{-1}\big) \AD(g'_S)\big( g'_F\varphi(g) \big)\\
		&=\big( \AD(g'_S)\circ\varphi \big)(g).
\end{split}  
\]
The first property yields
\[ 
  d\big( L_{\varphi(g')^{-1}}\circ\varphi\circ J \big)_{(0,0,0)}=\Ad(g'_S)\circ dI_{(0,0,0)}=d(\Phi_{ e^{a\BX}})_{e}\circ dI_{(0,0,0)}.
\]
Since $\omega_e$ is invariant under $\Phi_{ e^{a\BX}}$, it remains to be proved that $\omega'_{e'}\circ dJ_{(0,0,0)}=\omega_e\circ dI_{(0,0,0)}$. This is true because, in these coordinates, both sides applied on vectors $(Y_a,Y_v,Y_z)$ and $(Z_a,Z_v,Z_z)$ give
\[ 
  2(Y_aZ_z-Z_aY_z)+\Omega(Y_v,Z_v),
\]
so $\varphi$ is a symplectomorphism.

\end{proof}

Now, if $K$ is the kernel on $F$, we define the kernel on $F'$ by
\begin{equation}\label{EqKerRprime}
\begin{aligned}
 K'\colon F'\times F'\times F'&\to \eC\\
K'&=\varphi^*K. 
\end{aligned}
\end{equation}

\begin{theorem} 
The kernel $K'$ is\index{Kernel!for extension $(\BX ,0,2)$ of Heisenberg}
\begin{itemize}
\item left-invariant under $F'$,
\item associative on $F'$.
\end{itemize}
 \label{ThoDefoHeizAvecB}
\end{theorem}
\begin{proof}
For left-invariance, let $g'=J(a,v,z)$. We have
\[ 
  L_{g'}^*K'=\big( \varphi\circ L_{J(a,v,z)} \big)^*K=\big( \Phi_{ e^{a\BX}}\circ L_{g'_F}\circ\varphi \big)^*K=\varphi^*L_{g'_F}^*\Phi_{ e^{a\BX}}^*K=K',
\]
because of left-invariance of $K$ under $F$ and its invariance under $\Phi$. Associativity can be checked in much the same way as in lemma \ref{LemKerINvarIsom}.
\end{proof}

Let $F=F(\mtu,0,2)$ and $F'=F'(\BX,0,2)$. By proposition \ref{PropkernelinvarSp},  the kernel $K$ on $F$ is invariant under $\SP(V,\Omega)$, i.e. $\Phi^*_{M}K=K$ for all $M\in\SP(V,\Omega)$. The action of $\SP(V,\Omega)$ on $F$ is given by
\[ 
  \Phi_{M}\big( I(a,v,z) \big)=I(a,Mv,z).
\]
Define the map $\Phi'_M\colon F'\to F'$,
\begin{equation}   \label{EqDefPhiprimeM}
\Phi_{M}'\big( J(a,v,z) \big)=J\big( a, e^{-a\bar{\BX }}M e^{a\bar{\BX }}v,z \big)
\end{equation}
which fulfils 
\[ 
  \Phi_{M}\circ\varphi=\varphi\circ\Phi_{M}'.
\]
Thus, using the $\SP(V,\Omega)$-invariance of $K$, we have
\[ 
\phi_{M}'{}^*K'=(\varphi\circ\phi_{M}')^*K=(\phi_{M}\circ\phi)^*K=\varphi^*K=K'.
\]
This proves that $K'$ is also invariant under $\SP(V,\Omega)$ too.

\subsection{Jump from one kernel to another}

We have a kernel for the extensions $F_{\delta E^*}(d\mtu,0,2d)$ and $F_{\delta E^*}(\BX ,0,2)$. We can consider the isomorphism $L\colon F(\BX ,0,2)\to F(d\BX ,d\mu,2d)$ which is the lift of
\begin{equation}
\begin{aligned}
 dL\colon \mF(\BX ,0,2)&\to \mF(d\BX ,d\mu,2d) \\ 
  (a,v,z)&\mapsto (da,v+au,z).
\end{aligned}
\end{equation}
If $K'$ is a kernel on $F_{\delta E^*}(d\BX ,d\mu,2d)$, then
\[ 
  K_{0}=L^*K'
\]
is a kernel on $F_{\delta E^*}(\BX ,0,2)$.

An action $\Phi_{0}(M)\colon F(d\BX ,d\mu,2d)\to F(d\BX ,d\mu,2d)$ is given by
\begin{equation}
  \Phi_{0}(M)=L^{-1}\circ\Phi'(M)\circ L
\end{equation}
where $\Phi'(M)\colon F(\BX ,0,2)\to F(\BX ,0,2)$ is given by equation \eqref{EqDefPhiprimeM}. By lemma \ref{LemKerINvarIsom}, the kernel $K_{0}$ is left-invariant under the action of $F$ and invariant under the following action of $\SP(V,\Omega)$:
\[ 
  \Phi_{0}(M)^*K_{0}=K_{0}.
\]
\begin{lemma}

Let $\delta\eta^*$ and $\delta\xi^*$ be two exact forms on $\mF$ such that $\xi^*$ and $\eta^*$ belong to the same coadjoint orbit\index{Coadjoint!orbit}: there exists a $g\in F$ such that
\begin{equation}
\xi^*\circ\Ad(g)=\eta^*.
\end{equation}
A solution of the problem to find an automorphism $\sigma\colon F\to F$ such that
\begin{equation} \label{EqHypLemSigmaAD}
\delta\eta^*_{\sigma(h)}\big( d\sigma_{h}X_{h},d\sigma_{h}Y_{h} \big)=\delta\xi^*_{h}(X_{h},Y_{h})
\end{equation}
for all $h\in F$ and $X_{h}$, $Y_{h}\in T_{h}F$ is given by $\sigma=\AD(g^{-1})$.
\label{LemJumpCoadOrb} 
\end{lemma}

\begin{proof}
Transported to the identity, the condition \eqref{EqHypLemSigmaAD} becomes:
\[ 
\begin{split}
\delta\eta^*\big( dL_{\sigma(h)^{-1}}d\sigma_{h}X_{h},&dL_{\sigma(h)^{-1}}d\sigma_{h}Y_{h}  \big)\stackrel{!}{=}\delta\xi^*\big( dL_{h^{-1}}X_{h},dL_{h^{-1}}Y_{h} \big)\\
		&=\delta\eta^*\big( \Ad(g^{-1})dL_{h^{-1}}X_{h},\Ad(g^{-1})dL_{h^{-1}}Y_{h}  \big).
  \end{split}
\]
 If $X_{h}=\dsdd{ X_{h}(t)}{t}{0}$, we are searching for a $\sigma$ such that
\[ 
  \Dsdd{ \sigma(h)^{-1}\sigma\big( X_{h}(t) \big) }{t}{0}=\Dsdd{ \AD(g^{-1})\big( h^{-1}X_{h}(t) \big) }{t}{0}.
\]
Since $\sigma$ is a group isomorphism, $\sigma(h)^{-1}=\sigma(h^{-1})$ and the constraint on $\sigma$ becomes
\[ 
  \sigma\big(h^{-1}X_{h}(t)\big)=g^{-1}\big( h^{-1}X_{h}(t) \big)g.
\]
A solution is therefore
\begin{equation}
\sigma=\AD(g^{-1}).
\end{equation}

\end{proof}

\chapter{Toolbox}		\label{ChapTool}
\section{Connectedness of some usual groups}

\subsection{General results}

The following is a general result about Lie groups:
\begin{lemma}
If $G$ is a Lie group and $G_0$ is its identity component, the connected components of $G$ are lateral classes of $G_0$. More specifically, if $x\in G_1$, then $G_1=xG_0=G_0x$.
\label{LemConnSpecMo}
\end{lemma}

An other general result is lemma 2.4 of \cite{HelgasonSym} states that

\begin{lemma}
Connectedness of some usual groups:
\begin{itemize}
\item 
    The groups $\SU(p, q)$, $\SU^*(2n)$, $\SO^*(2n)$, $p(n, R)$, and $\SP(p, q)$ are
all connected.
\item 
    The group $\SO(p, q)$ ($0<p<p+q$) has exactly two connected components.
\end{itemize}
\label{LemConnSOpq}
\end{lemma}

\label{PgDisGeoConnSO}We are not going to prove this lemma here. Instead, we give some detail on the geometric nature of the two connected components of $\SO(p,q)$; a physical discussion in the case of $\SO(1,3)$ can be found in the reference \cite{Schomblond_em}. What is proved in \cite{HelgasonSym} is that $\SO(p,q)$ is homeomorphic to the topological product
\[ 
  \SO(p,q)=\SO(p,q)\cap\SU(p+q)\times \eR^{d}=\SO(p,q)\cap\SO(p+q)\times \eR^{d}
\]
for some $d\in\eN$. Hence an element of $\SO(p,q)$ reads
\[ 
  \begin{pmatrix}
A&0\\
0&B
\end{pmatrix}\times v
\]
where $v\in\eR^{d}$, $A\in \gO(p)$, $B\in\gO(q)$ are such that $\det A\det B=1$. The $v$ part corresponds to boost while $A$ and $B$ correspond to pure temporal and pure spatial rotations. An element of $\gO(n)$ has always determinant equals to $\pm 1$. Therefore one can decompose the rotation part as $(\det A=\det B=1)\otimes (\det A=\det B=-1)$. Both parts are connected.

Hence the first connected component contains $\mtu$ while the second one contains the element that simultaneously changes the sign of one spacial and one temporal direction.

\subsection{The quotient for anti de Sitter}

Homogeneous space considerations (see section \ref{SecSymeStructAdS}) will naturally lead us to define the anti de Sitter space as the quotient $G/H=\SO(2,l-1)/\SO(1,l-1)$ while the black hole definition (section \ref{SecCausal}) needs to consider Iwasawa decompositions of $G$. So we face the problem that the Iwasawa theorem \ref{ThoIwasawaVrai} only works with connected groups. In order to prevent any problems of this type, we prove now that, if $G_0$ and $H_0$ denote the identity component of $\SO(2,l-1)$ and $\SO(1,l-1)$ respectively, then $G/H=G_0/H_0$.

The groups that are considered here have only two connected components $G_0$ and $G_1$. We can chose $i_1\in G_1\cap H$ such that $i_1^2=\mtu$. Using lemma \ref{LemConnSpecMo}, it easy to prove that 
\begin{itemize}
\item $G_0G_0=G_0$,
\item $G_0G_1=G_1$,
\item $G_1G_1=G_0$.
\end{itemize}
For the last one, take $g$ and $g'$ in $G_1$. Then consider $g_0$ and $g'_0$ in $G_0$ such that $g=g_0i_1$ and $g'=g_0'i_1$. If $g_0(t)$ and $g'_0(t)$ are path from $\mtu$ to $g_0$ and $g_0'$, then the path $g_0(t)i_1g'_0(t)i_1$ is a path from $\mtu$ to $gg'$.

\begin{proposition}
The map
\begin{equation}
\begin{aligned}
 \psi\colon G/H&\to G_0/H_0 \\ 
[g]&\mapsto \overline{ g_0 } 
\end{aligned}
\end{equation}
where we define $g_0=g$ when $g\in G_0$ or $g_0=gi_1$ when $g\in G_1$ is a diffeomorphism.  The classes are $[g]=\{ gh\tq h\in H \}$ and $\overline{ g }=\{ gh_0\tq h_0\in H_0 \}$.
\label{PropGHconn}
\end{proposition}

\begin{proof}
First we prove that $\psi$ is well defined. For that we suppose that $[g]=[g']$. There are three cases:
\begin{enumerate}
\item The elements $g$ and $g'$ both belong to $G_0$. In this case, $g'=gh_0$ with $h_0\in H_0$ and $\overline{ gh }=\overline{ g }$.
\item The element $g$ belongs to $G_0$ while $g'$ belongs to $G_1$. In this case, $g'=gh$ with $h=h_0i_1$ and $h_0\in H_0$. Then $\psi[g]=\overline{ g }$ and $\psi[g']= \overline{ (gh_0i_1)_0 }=\overline{ gh_0i_1i_1 }=\overline{ gh_0 }=\overline{ g } $.
\item The case with $g$ and $g'$ in $G_1$ is similar.
\end{enumerate}

The fact that the map $\psi$ is surjective is clear. For injectivity, let $\psi[g]=\psi[g']$, i.e. there exists a $h_0$ in $H_0$ such that $g'_0=g_0h_0$. Thus we have $g'i_1^k=gi_1^lh_0$ with $k,l=0,1$ following the cases. Then $g'=gi_1^lh_0i_1^k$ in which $i_1^lh_0i_1^k$ belongs to $H$, so that $[g']=[g]$.

\end{proof}

\section{Iwasawa decomposition of Lie groups}

In this section, we show the main steps of the Iwasawa decomposition for a semisimple Lie group. We will by the way fix certain notations. For proofs, the reader will see \cite{Knapp} VI.4 and \cite{Helgason} III,\S\ 3,4 and VI,\S\ 3. In the whole section, $G$ denotes a semisimple group, and $\lG$ its real Lie algebra. The two main examples that are widely used during the thesis are $\SL(2,\eR)$ and $\SO(2,n)$.

\subsection{Cartan decomposition}

\begin{definition}
An involutive automorphism $\theta$ on a \emph{real} semi simple Lie algebra $\lG $ for which the form $B_{\theta}$,
\begin{equation}
          B_{\theta}(X,Y):=-B(X,\theta Y)
\end{equation}
($B$ is the Killing form on $\lG$) is positive definite is a \defe{Cartan involution}{Cartan!involution}.
\end{definition}

\begin{proposition}
There exists a Cartan involution for every real semisimple Lie algebra.
\end{proposition}

See \cite{Helgason}, theorem 4.1.  Since $\theta^2=id$, the eigenvalues of a Cartan involution are $\pm 1$, and we can define the \defe{Cartan decomposition}{Cartan!decomposition}\index{Decomposition!Cartan} $\lG$
\begin{equation}
                    \lG=\lK\oplus\lP
\end{equation}
into $\pm1$-eigenspaces of $\theta$ in such a way that $\theta=(-\id)|_{\lP}\oplus \id|_{\lK}$. These eigenspaces are subject to the following commutation relations:
\begin{equation} \label{Ieq:comm_KP}
[\lK,\lK]\subseteq\lK,\quad[\lK,\lP]\subseteq\lP,\quad [\lP,\lP]\subseteq\lK.
\end{equation}
The dimension of of maximal abelian subalgebra of $\lP$ is the \defe{rank}{Rank of a Lie algebra} of $\lG$. One can prove that it does not depend on the choices (Cartan involution and maximal abelian subalgebra). Let $\lA$ be one of such maximal abelian subalgebras.

\begin{lemma}
	The set of operators $\ad(\lA)$ is an abelian algebra and the elements are self-adjoint.
\end{lemma}

\subsection{Root space decomposition}

From the lemma, the operators $\ad(H)$ with $H\in\lA$ are simultaneously diagonalisable. There exists a basis $\{ X_i \}$ of $\lG$ and linear maps $\lambda_i\colon \lA\to \eR$ such that 
\[
	\ad(H)X_i=\lambda_i(H)X_i.
\]
 For any $\lambda\in\lA^*$, we define
\begin{equation}
		\lG_{\lambda}=\{X\in\lG|(\ad H)X=\lambda(H)X,\forall H\in\lA\}.
\end{equation}
Elements $0\neq\lambda\in\lA^*$ such that $\lG_{\lambda}\neq 0$ are called \defe{restricted roots}{Root!restricted} of $\lG$. The set of restricted roots is denoted by $\Sigma$, and have the important property to span (among with $\lA$ itself) the whole space:
\begin{equation}   \label{Ieq:somme_de_G}
             \lG=\lG_0 \oplus_{\lambda\in\Sigma}\lG_{\lambda},
\end{equation}
see \cite{Helgason} theorem 4.2 for a proof. This decomposition is called the \defe{restricted root space decomposition}{Root!space!decomposition}\index{Decomposition!root space}. Other properties of the root spaces are listed in the following proposition.

\begin{proposition}
The spaces $\lG_{\lambda_i}$ satisfy also:
\begin{enumerate}
\item $[\lG_{\lambda},\lG_{\mu}]\subseteq\lG_{\lambda+\mu}$,
\item $\theta\lG_{\lambda}=\lG_{-\lambda}$; in particular, when $\lambda$ belongs to $\Sigma$, $-\lambda$ belongs to $\Sigma$ too,
\item $\lG_0=\lA\oplus\lZ_{\lK}(\lA)$ orthogonally.
\end{enumerate}
\end{proposition}

\subsection{Iwasawa decomposition}

\begin{definition}
Let $V$ be a vector space. A \defe{positivity notion}{Positivity} (see \cite{Knapp}, page~154) is the data of a subset $V^+$ of $V$ such that
\begin{enumerate}
\item for any nonzero $v\in V$, $v\in V^+$ \emph{xor} $-v\in V^+$,
\item for any $v$, $w\in V^+$ and any $\mu\in\eR^+$, the elements $v+w$ and $\mu v$ are positive.
\end{enumerate}
\end{definition}
A \defe{positive}{} element of $V$ is an element of $V^+$. When $v$ is positive, we note $v>0$. 
Let us consider a notion of positivity on $\lA^*$ and denote by $\Sigma^+$ the set of positive roots. We define
\begin{equation}
      \lN:=\oplus_{\lambda\in\Sigma^+}\lG_{\lambda}.
\end{equation}
The \defe{Iwasawa decomposition}{Iwasawa decomposition}\index{Decomposition!Iwasawa} is given by the following theorem (\cite{Knapp}, theorem 5.12):

\begin{theorem}
Let $G$ be a linear connected semisimple group and $A=\exp\lA$, $N=\exp\lN$ where $\lA$ and $\lN$ are the previously defined algebras. Then $A$, $N$ and $AN$ are simply connected subgroups of $G$ and the multiplication map
\begin{equation}
\begin{aligned}
  \phi\colon A\times N\times K&\to G \\ 
 (a,n,k)&\mapsto ank 
\end{aligned}
\end{equation}
is a global diffeomorphism. In particular, the Lie algebra $\lG$ decomposes as vector space direct sum
\begin{equation}
            \lG=\lA\oplus\lN\oplus\lK.
\end{equation}
 The group $AN$ is a solvable subgroup of $G$ which is called the Iwasawa group, or Iwasawa component of $G$.
\label{ThoIwasawaVrai}
\end{theorem}

Notice that $A$, $N$ and $K$ are unique up to isomorphism. Their matricial representation of course depend on choices.

\section{Introduction to homogeneous spaces} \label{SecRoughomo}

Most of the material of this section can be found in a more general framework in the references \cite{Helgason, Loos, kobayashi, kobayashi2}. 

\subsection{Fundamental and invariant fields} \label{Subsec_Funda_conv}

Let $G$ be a Lie group with Lie algebra $\lG$. For each element of $\lG$, there are two distinguished vector fields on $G$, the \defe{left-invariant}{Left-invariant!vector field} and the \defe{right-invariant}{Right-invariant!vector field} one:
\begin{align}
\tilde X_g&=\Dsdd{  ge^{tX} }{t}{0}	&\underline X_g&=\Dsdd{ e^{tX}g }{t}{0}\\
dL_h\tilde X_g&=\tilde X_{hg}		&dR_h\underline X_g&=\underline X_{gh}.
\end{align}

When $G$ is a Lie group with an action on the manifold $M$ denoted by
\begin{equation}
\begin{aligned}
 \tau\colon G\times M&\to M \\ 
(g,x)&\mapsto \tau_g(x),
\end{aligned}
\end{equation}
we define the \defe{fundamental vector field}{Fundamental vector field} associated with $X\in\lG$ on the point $x\in M$ by
\begin{equation}			\label{EqDefChmpFonfOff}
X^*_x=\Dsdd{ \tau_{ e^{-tX}}(x) }{t}{0}.
\end{equation}
An usual case is the one of a Lie group acting on itself for which we have
\begin{equation}		\label{EqChmpFondGp}
  X^*_g=\Dsdd{ e^{-tX}g }{t}{0}.
\end{equation}

\subsection{Homogeneous spaces}		\label{SubSechoappahomsp}

An \defe{homogeneous space}{Homogeneous!space} is a differentiable manifold which posses a transitive diffeomorphism group. An important class of homogeneous spaces are quotients $M=G/H$ of a Lie group $G$ by a closed subgroup $H$. In this case, we use the classes at right:
\[ 
  [g]=\{ gh\tq h\in H \}
\]
and the action at left:
\[ 
  \tau_g[g']=[gg'].
\]
The canonic projection is $\pi\colon G\to M$ and we denote $\mfo=[e]$. We will only deal with this kind of homogeneous spaces. The Lie algebras of $G$ and $H$ are denoted by $\lG$ and $\lH$ respectively.

One know that (almost) every homogeneous space is of this kind in the following way.  Let $M$ be a homogeneous space and $\mfo$, a point of $M$. We consider $G$, a group which acts transitively on $M$ (in particular, $G\mfo=M$) and $H$, the subgroup of $G$ which fixes $\mfo$. Then, one proves that the map $[g]\mapsto g\mfo$ is a homogeneous space isomorphism between $M$ and $G/H$. 

One can prove that $\ker(d\pi_e)=\lH$, and from the very definition of the objects, one has
\begin{equation}  \label{IEqdpigdtaudpi}
  d\pi_g \circ dL_g=d\tau_g\circ d\pi_e.
\end{equation}
For sake of simplicity, we will use the notation $\mu_g=\tau_g\circ \pi$.

The homogeneous space $G/H$ is endowed with its \defe{natural topology}{Natural topology} which is defined by the requirement that the projection $\pi$ is continuous and open. We refer to \cite{Helgason} for the properties of that topology.

\begin{definition}
The homogeneous space $M=G/H$ is \defe{reductive}{Homogeneous!space reductive} is there exists a subspace $\lQ$ of $\lG$ such that 
\begin{align*}
\lG&=\lQ\oplus\lH&[\lH,\lQ]&\subset\lQ.
\end{align*}
\end{definition}

\begin{proposition} 
In an reductive homogeneous space, the restriction of the projection $d\pi_e\colon \lQ\to T_{\mfo}M$ is an isomorphism.
\label{IPropdpiisomMTM}
\end{proposition}

\begin{proof}
The map $\dpt{d\pi_e}{\lG=\lQ\oplus\lH}{T_{\mfo}M}$ is of course surjective; then, since $\lH$ is the kernel, $\dpt{d\pi_e}{\lQ}{T_{\mfo}M}$ must be surjective too. Now, if we have $d\pi_eX=d\pi_eY$ for $Y$, $X\in\lQ$, the difference $(X-Y)$ must belongs to the kernel of $d\pi_e$ which is nothing but $\lH$. This situation is impossible because $\lG=\lH\oplus\lQ$ is a direct sum.
\end{proof}

We can generalize this proposition by considering the space $\lQ_g=dL_g\lQ$. Using equality \eqref{IEqdpigdtaudpi}, the map $d\pi_g\circ dL_g\colon \lQ\to T_{[g]}M$ is an isomorphism. Since, by definition, the map $dL_{g^{-1}}\colon \lQ_g\to \lQ$ is an isomorphism, we conclude that

\begin{corollary} 
The restriction $d\pi_g\colon \lQ_g\to T_{[g]}M$ is a vector space isomorphism.
\label{ICordpiietwii}
\end{corollary}

\subsection{Killing induced product}		\label{SubsecKillHomo}

Since the Killing form $B$ is an $\Ad_H$-invariant product on $\lQ$, we can define
\begin{equation}
B_g(X,Y)=B_e(dL_{g^{-1}}X,dL_{g^{-1}}Y)
\end{equation}
which descent (see \cite{Kerin} for properties) to a homogeneous metric on $T_{[g]}M$:
\begin{equation}  \label{EqDefMetrHomo}
B_{[g]}(d\pi X,d\pi Y)=B_g(\pr X,\pr Y)
\end{equation}
where $\dpt{\pr}{T_gG}{dL_g\lQ}$ is the canonical projection. An useful property of that projection is $\pr(dL_gX)=dL_gX_Q$ when $X=X_Q+X_H$. Using that property, we can write the product under the more manageable form
\[ 
  B_{[g]}(d\mu_gX,d\mu_gY)=B_e(\pr X,\pr Y)
\]
for all $X$, $Y\in\lG$.

Although equation \eqref{EqChmpFondGp} looks like \eqref{EqDefChmpFonfOff}, we find a major difference here: the norm of $q_i^*[g]$ is not a constant. One should expect that it was a constant because \eqref{EqDefChmpFonfOff} expresses a left translation while the Killing form is invariant under left translations. But the metric \eqref{EqDefMetrHomo} is a composition of the Killing form with a projection. Let us study this case in details in computing the product of two vectors of the form
\[ 
  X^*_{[g]}=d\pi\Dsdd{  e^{-tX}g }{t}{0},
\]
with $X\in\lQ$:
\[ 
\begin{split}
  B_{[g]}(X^*,Y^*)&=B_g\big( \pr\Dsdd{  e^{-tX}g }{t}{0},\pr\Dsdd{  e^{-tY}g }{t}{0} \big)\\
		&=B_g\big( dL_g\pr\Ad(g^{-1})X,dL_g\pr\Ad(g^{-1})Y \big)\\
		&=B_e\Big(   \big( \Ad(g^{-1})X \big)_{\lQ},\big( \Ad(g^{-1})Y \big)_{\lQ}  \Big)\\
		&\neq B_e\Big(   \Ad(g^{-1})X_{\lQ},\Ad(g^{-1})Y_{\lQ}  \Big)\\
		&=B_e(X,Y)
\end{split}
\]
where the symbol $\neq$ has to be understood as ``not equal in general'' because equality holds of course for certain particular vectors such as zero.
\section{Toolbox for \texorpdfstring{$SL(2,\eR)$}{SL2R}}	\label{SecToolSL}

\subsection{Iwasawa decomposition}
\index{Iwasawa decomposition!of $SL(2,\eR)$}

Let $G=\SL(2,\eR)$ the group of $2\times 2$ matrices with unit determinant. The Lie algebra $\lG=\gsl(2,\eR)$ is the algebra of matrices with vanishing trace:
\begin{equation}
\begin{split}
 \lG &=  \{ X\in\End(\eR^2)\tq \tr(X) = 0\} \\
&=\left\{ \begin{pmatrix}
x & y \\
z & -x
\end{pmatrix}\textrm{ with }x,y,z\in\eR  \right\}. 
\end{split}
\end{equation}
The following elements will be intensively used:
\[
H=\begin{pmatrix}
1 & 0 \\
0 & -1
\end{pmatrix}
,\quad
  E=\begin{pmatrix}
0 & 1 \\
0 & 0
\end{pmatrix}
,\quad
 F=\begin{pmatrix}
0 & 0 \\
1 & 0
\end{pmatrix},
\quad
T=\begin{pmatrix}
0&1\\
-1&0
\end{pmatrix}
\]
where $T=E-F$ has been introduced for later convenience. The commutators are
\begin{subequations}\label{EqTableSLdR}
\begin{align}  
  [H,E]&=2E	&[T,H]&=-2T  \\
  [H,F]&=-2F	&[T,E]&=H   \\
  [E,F]&=H	&[T,F]&=H.
\end{align}
\end{subequations}
Notice that the sets $\{ H,E,F \}$, $\{ H,E,F \}$ and $\{ H,E+F,T \}$ are basis. A Cartan involution is given by $\theta(X)=-X^t$, and the corresponding Cartan decomposition is
\[
\begin{split}
   \lK&=\Span\{ T \}\\
\lP&=\Span\{ H,E+F \}
\end{split}
\]

Up to some choices, the Iwasawa decomposition\label{pg_iwasldr} of the group $\SL(2,\eR)$ is given by the exponentiation of $\lA$, $\lN$ and~$\lK$
\begin{equation}
\begin{aligned}
  \lA&=\Span\{ H \}
&\lN&=\Span\{ E \}
&\lK&=\Span\{T\},
\end{aligned}
\end{equation}
so that
\begin{equation}\label{eq:expo_ANK}
A=\begin{pmatrix}
e^a & 0 \\
0 & e^{-a}
\end{pmatrix}\quad
N=\begin{pmatrix}
1 & l \\
0 & 1
\end{pmatrix}\quad
K=\begin{pmatrix}
\cos k & \sin k \\
-\sin k & \cos k
\end{pmatrix}.
\end{equation}

A common parametrization of $AN$ by $\eR^2$ is provided by
\begin{equation}   \label{EqParmalSL} 
(a,l)=
\begin{pmatrix}
  e^a&le^a\\
  0  &e^{-a}
\end{pmatrix}.
\end{equation}
One immediately has the following formula for the left action of $AN$ on itself:
\[
  L_{(a,l)}(a',l')=\begin{pmatrix}
e^{a+a'} & e^{a+a'}l'+e^{a-a'}l \\
0 & e^{-a-a'}
\end{pmatrix}=(a+a',l'+e^{-2a'}l).
\]
In this setting, the inverse is given by $(a,l)^{-1}=(-a,-l e^{2a})$.  

\subsection{Killing form}

The Killing form $B(X,Y)=\tr(\ad X\circ\ad Y)$ takes the following values:
\begin{subequations}
\begin{align}
B(T,H)&=0  & B(H,H)&=8\\
B(T,E)&=-4 & B(E,E)&=0\\
B(H,E)&=0  &  B(T,T)&=-4.
\end{align}
\end{subequations}
Expressed in the basis $\{H,E,F\}$, the matrix of the Killing form reads
\begin{equation}
B=
\begin{pmatrix}
8&&\\
&&4\\
&4&
\end{pmatrix}
\end{equation} 
while, in the basis  $\{H,E+F,T\}$, we find
\begin{equation}   \label{EqBHEFTsldR}
B=
\begin{pmatrix}
8\\
&8\\
&&-8
\end{pmatrix}.
\end{equation}
The latter is the reason of the name of the vector $T$: the sign of its norm is different, so that $T$ is candidate to be a time-like direction.

\subsection{Abstract root space setting}

Looking on the table \eqref{EqTableSLdR} from an abstract point of view, we see that $E$ and $F$ are eigenvectors of $\ad(H)$ with eigenvalues $2$ and $-2$. So $\lA=\lG_0=\eR H$; $\lG_2=\eR E$; and $\lG_{-2}=\eR F$. Using a more abstract notation, the table of $\SL(2,\eR)$ becomes
\begin{subequations}  \label{subeq_rootSLR}
\begin{align}
  [A_{0},A_{2}]&=2A_{2}\\
	[A_{0},A_{-2}]&=-2A_{-2}\\
	[A_{2},A_{-2}]&=A_{0}.
\end{align}
\end{subequations}

\subsection{Isomorphism}

As pointed out in the chapter II, \S6 of \cite{Knapp_reprez}, the map (seen as a conjugation in $\SL(2,\eC)$)
\begin{equation}
\begin{aligned}
 \psi\colon \SU(1,1)&\to \SL(2,\eR) \\ 
  U&\mapsto AUA^{-1} 
\end{aligned}
\end{equation}
with $A=\begin{pmatrix}
1&i\\i&1
\end{pmatrix}$ is an isomorphism between $\SL(2,\eR)$ and $\SU(1,1)$.

\section{Root spaces for \texorpdfstring{$\so(2,1)$}{so(2,1)}}

The algebra $\so(2,1)$ is made up from $3\times 3$ matrices such that $X^t\eta+\eta X=0$ with vanishing trace. If we choice $\eta=diag(-,-,+)$, we find matrices of the form
\[ 
  \so(2,1)\leadsto\begin{pmatrix}
a	&u^t\\
u	&0
\end{pmatrix}
\]
where $a$ is an antisymmetric $2\times 2$ matrix and $u$ is any $1\times 2$ matrix. We find the Cartan decomposition
\[ 
  \iK\leadsto\begin{pmatrix}
\so(2)\\&0
\end{pmatrix},\qquad
\iP\leadsto\begin{pmatrix}
0&0&u_1\\
0&0&u_2\\
u_1&u_2&0
\end{pmatrix}.
\]
If one chooses
\[ 
  J=\begin{pmatrix}
0&0&1\\0\\1
\end{pmatrix}
\]
as generator for the abelian subalgebra of $\iP$, one finds
\[ 
  V_1=\begin{pmatrix}
0&1&0\\
-1&0&1\\
0&1&0
\end{pmatrix},\qquad
V_{-1}=\begin{pmatrix}
0&1&0\\
-1&0&-1\\
0&-1&0
\end{pmatrix}
\]
as eigenvectors for $\ad(J)$ with eigenvalues $1$ and $-1$. The root space decomposition commutator table of $\so(2,1)$ is thus given by
\begin{subequations}
\begin{align}
 [V_0,V_1]&=V_1\\
[V_0,V_{-1}]&=-V_{-1}\\
[V_1,V_{-1}]&=-2V_0.
\end{align}
\end{subequations}
Notice that the map $\phi(A_0)=2V_0$, $\phi(A_2)=V_1$, $\phi(A_{-2})=-V_{-1}$ provides an isomorphism between this table and the one of $\gsl(2,\eR)$, equations \eqref{subeq_rootSLR}. This fact assures a Lie algebra isomorphism $\gsl(2,\eR)\simeq\so(2,1)$. We actually have a stronger result: 
\begin{proposition}
The group $\SL(2,\eR)$ is a double-covering of the identity component $\SO_0(1,2)$.
\end{proposition}

\section{Iwasawa decomposition for \texorpdfstring{$\SO(1,n)$}{SO1n}}
\index{Iwasawa decomposition!of $\SO(1,n)$}

We saw in the section \ref{SecSymeStructAdS} that the quotient $\SO(2,n)/\SO(1,n)$ has a particular importance. Hence, we will work out the Iwasawa decompositions of these groups imposing certain compatibility conditions. We already build the Cartan involution $\theta$ in such a way that $[\sigma,\theta]=0$.

Now we show that $\theta$ descent to a Cartan involution on $H$. It is clear that the restriction of $\theta$ is an involutive automorphism of $H$. Lemma \ref{lem:Killing_ss_descent} assures us that the restriction of the Killing form of $G$ to $H$ is the Killing form of $H$, so that the condition of positivity of $B_{\theta}$ holds on $H$ as well as on $G$. Last, $\theta$ leaves $\sH$ invariant. Indeed suppose that  $\theta X_{\sH}=X'_{\sH}+X_{\sQ}\in\sH\oplus\sQ$. Then $\sigma\theta X_{\sH}=h'-q$ and $\theta\sigma h=h'+q$; since $[\theta,\sigma]=0$, we have $q=0$.

All that proves that we can use the same Cartan involution on $G$ as well as on $H$. Since $\theta=id|_{\iK}\oplus(-id)|_{\iP}$, it is clear that 
\begin{align}
   \iKH&=\iK\cap\sH,
   &\iPH&=\iP\cap\sH
\end{align}
is the Cartan decomposition of $\sH$. We can write explicit matrices as
\begin{equation}
\iKH=so(n)\leadsto
\begin{pmatrix}
  0&0&\cdots  \\
  0&0&\cdots  \\
  \vdots&\vdots& B
\end{pmatrix},
\end{equation}
where $B$ is skew-symmetric, and
\begin{equation}
\iPH\leadsto
\begin{pmatrix}
  0&0&\cdots  \\
  0&0& u^t  \\
  \vdots&u& 0
\end{pmatrix}.
\end{equation}
One remark that there are no two-dimensional subalgebra of $\iP_{\sH}$. So $\iA_{\sH}$ reduces to the choice of any element $J_1\in\iP_{\sH}$.  A positivity notion is easy to find: the form $\omega\in\iAH^*$ such that $\omega(J_1)=1$ is positive while $-\omega$ is negative. We choose
\[
  J_1=\begin{pmatrix}&0\\0&0&0&1\\&0\\&1\end{pmatrix}.
\]
The computation of $\iNH=\{X\in\sH\tq (\ad J_1)X=X\}$ yields the following :
\begin{equation}\label{eq:re_N_H}
\iNH\leadsto
\begin{pmatrix}
     \begin{pmatrix} &&0&0\\&&a&0\\0&a&0&-a\\0&0&a&0\end{pmatrix}
     & \begin{pmatrix}\cdots&0&\cdots\\\leftarrow&\overline{v}&\rightarrow\\
        \cdots&0&\cdots\\\leftarrow&\overline{v}&\rightarrow\end{pmatrix}\\
     \begin{pmatrix}\vdots&\uparrow&\vdots&\uparrow\\
                    0& \overline{v}&0&-\overline{v}\\
		    \vdots&\downarrow& \vdots&\downarrow \end{pmatrix}
     &0		    
\end{pmatrix}.
\end{equation}
Finally, we consider the algebra $\iRH=\iAH\oplus\iNH$.
\section{Iwasawa decomposition for \texorpdfstring{$\SOdn$}{SO2n}} \label{subsecIwasawa_un}
\index{Iwasawa decomposition!of $\SO(2,n)$}

As seen in the general construction and in previous examples, the Iwasawa decomposition of a group or an algebra depends on several choices. We will study two out of them in the case of $\SO(2,n)$ and see that some ``compatibility conditions'' with the decomposition of $\SO(1,n)$ and the symmetric space structure of $AdS$ (see section \ref{SecSymeStructAdS}) fix most of choices. 

The Lie algebra $\sG=\mathfrak{so}(2,n)$ is the set 
\begin{equation}\label{def_sodn}
\big\{X\in M_{(2+n)\times (2+n)}\textrm{ such that }X^t\eta+\eta X=0\textrm{ and } \tr X=0\big\}
\end{equation}
where $\eta$ is the diagonal metric $\eta=diag(-,-,+,\ldots,+)$. An element of $\mathfrak{so}(2,n)$ can be written as
$
    X=\begin{pmatrix}
a & u ^t \\
v & B
\end{pmatrix}
$
with the matrices $a\in M_{2\times 2}$, $u\in M_{n\times 2}$, $v\in M_{n\times 2}$, and $B\in M_{n\times n}$. The conditions in \eqref{def_sodn} give: $a=-a^t$, $u=v$, and $B=-B^t$. Hence, a general matrix of $\sodn$ is given by
\begin{equation}	\label{eq:gene_sodn}
X=\begin{pmatrix}
a & u^t \\
u & B
\end{pmatrix}
\end{equation}
where $a,B$ are skew-symmetric.

\subsection{Cartan decomposition}		\label{SubSecCartandeuxN}

The Cartan decomposition of $\so(2,n)$ associated with the Cartan involution $\theta(X)=-X^t$ is
\begin{align}\label{K_et_P}
   \iK&\leadsto
\begin{pmatrix}
\mathfrak{so}(2) \\
 & \mathfrak{so}(n)
\end{pmatrix},
&\iP&\leadsto\begin{pmatrix}
0 & u^t \\
u & 0
\end{pmatrix}.
\end{align}
 Elements of $\SO(2)$ are represented by
\[ 
  \begin{pmatrix}
\cos\mu&\sin\mu\\
-\sin\mu&\cos\mu
\end{pmatrix}.
\]
A common abuse of notation in the text will be to identify the angle $\mu$ with the element of $\SO(2)$ itself. In the same spirit, when we speak about a matrix of $A\in \SO(2)$, we mean a matrix whose upper left corner is $A$ and the rest is the unit matrix. For example, for $AdS_3$, the matrix $-\mtu\in \SO(2)$ is
\[
\begin{pmatrix}
\begin{matrix}
-1&0\\
0&-1
\end{matrix}\\
&1\\
&&1
\end{pmatrix}.
\]

Remark that $\iK$, the compact part of $\mG$ is made up from ``true'' rotations while $\iP$ contains boosts.  This remark allows us to guess a right choice of maximal abelian subalgebra in $\iP$. Indeed elements of $\iA$ must be boosts and the fact that there are only two time-like directions restricts $\sA$ to a two dimensional algebra. Up to reparametrization, it is thus generated by $t\partial_x+x\partial_t$ and $u\partial_y+y\partial_t$. Hence the following choice seems to be logical:
\begin{align}   \label{EqDevineA}
   J_1&=
\begin{pmatrix}
&0\\
0&0&0&1\\
&0\\
&1
\end{pmatrix}\in\sH,
&J_2&=q_1=
\begin{pmatrix}
0&0&1&0\\
0\\
1\\
0
\end{pmatrix}\in\sQ.
\end{align}

\subsection{Maximal abelian subalgebra}

The generators of $\iA$ that we choose are the following linear combination of $J_1$ and $J_2$:
\begin{equation}\label{mtr_A}
       H_p=E_{13}+E_{31}+\mue{p}(E_{24}+E_{42}).
\end{equation}

\subsection{Nilpotent part}

We search the eigenvectors and eigenvalues of $\ad(H_p)$ under the form $E=\begin{pmatrix}A&B\\C&D\end{pmatrix}$ with $A\in M_{4\times 4}$, $B\in M_{(n-2)\times 4}$, $C\in M_{4\times (n-2)}$, $D\in M_{(n-2)\times (n-2)}$.  Remark that, thanks to \eqref{K_et_P}, the matrix $C$ is completely determined by $B$: 
\begin{equation}\label{C_de_B}
\left\{\begin{aligned}C^{i1}&=B^{1i},&C^{i2}&=B^{2i},\\
                    C^{i3}&=-B^{3i},&C^{i4}&=-B^{4i}.\end{aligned}\right.
\end{equation}
The equation to be solved is
\begin{equation}                                                                        \label{eq_gene}
(\ad H_p)E=\begin{pmatrix}[N_p,A]&N_pB\\-CN_p&0\end{pmatrix}\stackrel{!}{=}\lambda_pE,
\end{equation}
with the notation $\lambda_p=\lambda(H_p)$.

\subsubsection{Search for two non zero eigenvalues}

By ``two non zero eigenvalues'', we mean a $\lambda\in\iA^*$ such that $\lambda_{1}\neq 0$ and $\lambda_{2}\neq 0$. We immediately find $D=0$.

The next step is to determine $B$ by the condition $N_pB=\lambda_p B$. We change the range of the indices. Now, $a$, $b:1\rightarrow 4$, and $i$, $j:5\rightarrow n+2$ and a few computation give $\sum_{a}^{}N_p^{ca}B^{ai}=\lambda_pB^{ci}$ (with sum over $a$). Taking successively $c=1,2,3,4$ and taking into account $\lambda_p\neq 0$, we find:
\begin{subequations}\label{pour_B}
\begin{align}
B^{3i}&=\lambda_pB^{1i} \\
\mue{p}B^{4i}&=\lambda_pB^{2i}\\
\lambda_p&=\pm 1.
\end{align}
\end{subequations}
We can check that the equations obtained by $-CN_p=\lambda_pC$ are exactly the one that we can find directly using \eqref{C_de_B} and \eqref{pour_B}.

 Now, we determine $A$ by the condition $[N_p,A]=\lambda_pA$. We know that $A$ and $N_p$ are $4\times 4$ matrices. Again, we redefine the range of the indices: $a=1,2$ and $i=3,4$. Symmetry properties of $A$ are $A^{ai}=A^{ia}$, $A^{ij}=-A^{ji}$ and $A^{ab}=-A^{ba}$, so that
\[
 A=A^{12}(E_{12}-E_{21})+A^{ai}(E_{ai}+E_{ia})+A^{34}(E_{34}-E_{43})
\]
Using equation \eqref{mtr_A}, a quite tedious (but direct) computation give the following for $[N_p,A]$:
\begin{equation}\label{grosse_A}
\begin{pmatrix}
0&A^{23}-(-1)^p A^{14}&0&A^{34}-(-1)^pA^{12}\\
-A^{23}+ (-1)^pA^{14}&0&A^{12}- (-1)^pA^{34}&0\\
0&A^{12}- (-1)^pA^{34}&0&A^{14}- (-1)^pA^{23}\\
A^{34}- (-1)^pA^{12}&0&A^{14}- (-1)^pA^{23}&0
\end{pmatrix}
\end{equation}
which has to be equated to $\lambda_pA$. We immediately have $A^{13}=A^{24}=A^{31}=A^{42}=0$. The others conditions are:
\begin{subequations}\label{pour_A}
\begin{align}
\lambda_pA^{12}=&A^{23}- (-1)^pA^{14}\\
\lambda_pA^{14}=&A^{34}- (-1)^pA^{12}
\end{align}
\end{subequations}
Since we are in the case $\lambda_p\neq 0$, using the fact that $\lambda_p=\pm 1$, we easily find $A=0$. Now, we define $\lambda_{\pm\pm}\in\iA^*$ by
\begin{eqnarray}
\begin{aligned}
\lambda_{++}(H_1)&=1 &\lambda_{++}(H_2)&=1,\\
\lambda_{+-}(H_1)&=1 &\lambda_{+-}(H_2)&=-1,\\
\lambda_{-+}(H_1)&=-1 &\lambda_{-+}(H_2)&=1,\\
\lambda_{--}(H_1)&=-1 &\lambda_{--}(H_2)&=-1.
\end{aligned}
\end{eqnarray}
Root spaces are (with $i:5\rightarrow n+2$):
\begin{eqnarray}
\begin{aligned}
\lambda_{++}&\leadsto V_i=E_{3i}+E_{1i}+E_{i1}-E_{i3},\\
\lambda_{-+}&\leadsto W_i=E_{4i}+E_{2i}-E_{i4}+E_{i2},\\
\lambda_{--}&\leadsto X_i=E_{3i}-E_{1i}-E_{i1}-E_{i3},\\
\lambda_{+-}&\leadsto Y_i=E_{4i}-E_{2i}-E_{i4}-E_{i2}.
\end{aligned}
\end{eqnarray}

\subsubsection{Search for eigenvalues with one zero}

We denote by $\lambda(a,b)$ the element of $\iA^*$ defined by $\lambda(x_1,x_2)(H_i)=x_i$. Equations \eqref{grosse_A} and \eqref{pour_A} with for example $a=0$ and $b\neq 0$ give $\lambda_2=\pm 2$. Serious computation give:
\begin{subequations}
\begin{equation}
\lambda(0,-2)\leadsto F=
\begin{pmatrix}
                0&1&0&1\\
		-1&0&-1&0\\
		0&-1&0&-1\\
		1&0&1&0
              \end{pmatrix},\quad
\lambda(2,0)\leadsto N
=\begin{pmatrix}
                0&1&0&1\\
		-1&0&1&0\\
		0&1&0&1\\
		1&0&-1&0
              \end{pmatrix}
\end{equation}
\begin{equation}
\lambda(0,2)\leadsto M
=\begin{pmatrix}
                0&1&0&-1\\
		-1&0&1&0\\
		0&1&0&-1\\
		-1&0&1&0
              \end{pmatrix},\quad
\lambda(-2,0)\leadsto L=
\begin{pmatrix}
                0&1&0&-1\\
		-1&0&-1&0\\
		0&-1&0&1\\
		-1&0&-1&0
              \end{pmatrix}.
\end{equation}
\end{subequations}
In these expressions, we only wrote the upper left part of the matrices which are zero everywhere else.

\subsubsection{Two vanishing eigenvalues}

We now search for a matrix $E$ of $\so(2,n)$ such that $(\ad H_p)E=0$ for $p=1,2$. Taking a look at \eqref{eq_gene}, we see that $D$ has no more constraints (apart the usual symmetries). The equation \eqref{grosse_A} gives us $A^{12}=A^{23}=A^{14}=A^{34}=0$, but $A^{13}$, $A^{24}$,$A^{31}$,$A^{42}$ are free. Therefore we can write:
\begin{equation}
  \mG_{\lambda(0,0)}=\{x(E_{13}+E_{31})+y(E_{42}+E_{24})+
    \left(\begin{array}{c|c}
    0&0\\
    \hline
    0&D
    \end{array}\right)
  \}.
\end{equation}
Remark that $\mG_{\lambda(0,0)}\cap\iP$ is spanned by matrices of the form
\[
\begin{pmatrix}
	0&0&x&0\\
	0&0&0&y\\
	x&0&0&0\\
	0&y&0&0
\end{pmatrix},
\]
so that $\mG_{\lambda(0,0)}\cap\iP=\iA$. This is a simple consequence of the very definition of $\iA$ as maximal abelian subalgebra of $\iP$.

\subsubsection{Choice of positivity}

We are now able to write down the component $\iN$. We just have to \oge select\fge\ some $\mG_{\lambda_{(a,b)}}$ with a notion of positivity. Our choice is:
 \begin{equation}
     \iN=\{V_i,W_i,M,N\},\\
 \end{equation}
with $i$, $j:5\to n+2$. A basis of $\iK$ is given by $K_{rs}=E_{rs}-E_{sr}$, and $K_a=E_{12}-E_{21}$ with $\dpt{r,s}{3}{n+2}$.  

Let us summarize the result obtained.  The Iwasawa decomposition of $\SOdn$ is given by
\begin{subequations}\label{eq:Iwasawa_explicite}
\begin{align}
  \mA&=\{H_1,H_2\}\\
  \mN&=\{V_i,W_j,M,N\}\\
  \mK&=\{\begin{pmatrix}a&0\\0&B\end{pmatrix}\}.
\end{align}
\end{subequations}
with $i$, $j:5\to n+2$, $a\in M_{2\times 2}$, $B\in M_{n\times n}$ skew-symmetric, and
$H_p=E_{13}+E_{31}+(-1)^p(E_{24}+E_{42})$.

The non-zero commutators in $\mA\oplus\mN$ are :
\begin{subequations}		\label{TabelPrem}
\begin{align}
[V_i,W_j]&=\delta_{ij}M &[W_i,N]&=-2V_i\\
[H_1,V_i]&=V_i          &[H_2,V_i]&=V_i\\
[H_1,N]&=2N          &[H_2,M]&=2M\\
[H_1,W_i]&=-W_i   &[H_2,W_i]&=W_i.
\end{align}
\end{subequations}
\subsection{Second Iwasawa decomposition for \texorpdfstring{$\SO(2,n)$}{SO(2,n)}}\label{subsecIwasawa_deux}

We know the $\pm 1$ eigenspaces decompositions $\sG\stackrel{\sigma}{=}\sH\oplus\sQ\stackrel{\theta}{=}\sK\oplus\sP$ with $[\sigma,\theta]=0$, and the Iwasawa decomposition $\sA_{\sH}\oplus\sN_{\sH}\oplus\sK_{\sH}$ of $\sH$. 

For compatibility and simplicity purposes, we want the Iwasawa decompositions $\sG=\sA\oplus\sN\oplus\sK$ of $\sG$ in such a way that $\sA_{\sH}\subset\sA$ and $\sN_{\sH}\subset\sN$.
We denote by $A$, $N$, $K$, $A_H$, $N_H$ and $K_H$ the analytic connected subgroups of $G$ whose Lie algebras are $\sA$, $\sN$, $\sK$, $\sA_{\sH}$, $\sN_{\sH}$, and $\sK_{\sH}$ respectively.

 For the $\sA$-part, we just perform a change of basis
\begin{subequations}
\begin{align}
J_{1}&=\frac{ 1 }{2}(H_{2}-H_{1})\\
J_{2}&=\frac{ 1 }{2}(H_{1}+H_{2}) 
\end{align}
\end{subequations}
in order to have $J_{1}\in\sP\cap\sH$ and $J_{2}\in\sP\cap\sQ$. The involution $\sigma$ has a simpler expression in this basis.
If $\alpha\in\sA^*$, the notation $\sigma^*\alpha$ means $f\circ\alpha$, and $\sG_{\alpha}$ denotes the root space associated with the form $\alpha$.

\begin{lemma} 
The involutions $\sigma$ and $\theta$ act on the root spaces by 
\begin{align}\label{eq:sig_G_alpha}
\sigma\sG_{\varphi}&=\sG_{\sigma^*\varphi},
  &\theta\sG_{\varphi}&=\sG_{-\varphi}
\end{align}
for all $\varphi\in\sA^*$,
  \label{LemSigmaRootSpaces}
\end{lemma}

\begin{proof}
Let $H\in \sA$ and $X\in\sG_{\varphi}$; by definition: $[H,X]=\varphi(H)X$. We have
\[ 
  \varphi(H)\sigma X=\sigma[H,X]=[\sigma H,\sigma X]=\varphi(\sigma H)\sigma X.
\]
We conclude that $\sigma X\in\sG_{\sigma^*\varphi}$. For the second equality, we take $X\in\sG_{\varphi}$ and 
\[ 
  [H,\theta X]=\theta[\theta H,X]=\theta\big( \varphi(\theta H)X \big)=-\varphi(H)\theta(X).
\]
\end{proof}

\begin{proposition} 
The involution $\sigma$ changes the sign of the $J_{2}^*$-part of the root spaces:
\[
	\sigma\sG_{(x,y)}=\sG_{(x,-y)}
\]
where $(x,y)$ denote the coordinates of a root in the basis $\{ J_1^*,J_2^* \}$ of $\sA$.
\end{proposition}

\begin{proof}
Since $\sigma$ is an involutive automorphism, it satisfies $[\sigma X,Y]=\sigma[X,\sigma Y]$. So when $X\in\sG_{(x,y)}$, we have $(\ad(J_{1}))(\sigma X)=x\sigma X$ and $\ad(J_{2})(\sigma X)=-y\sigma X$.

\end{proof}

It is also clear that $\sigma^*\alpha=0$ implies $\alpha(J_2)=0$. Indeed, $(\sigma^*\alpha)(J_2)=\alpha(-J_2) =-\alpha(J_2)$, because $J_{2}\in\sQ$.

 We turn now our attention to the Iwasawa decomposition. As before we are searching for matrices under the form $E=\begin{pmatrix}A&B\\C&D\end{pmatrix}$ with dimensions $4+(n-2)$. The condition which determines the root spaces reads
\begin{equation}
  (\ad J_i)E=
\begin{pmatrix}
  [j_i,A]&j_iB\\
  -Cj_i&0
\end{pmatrix}
\stackrel{!}{=}\lambda_i
\begin{pmatrix}A&B\\C&D\end{pmatrix}.
\end{equation}

The computations are rather the same as the ones of the first time. The result is\label{pg:root_n}
\begin{equation}
\sG_{(0,0)}\leadsto
\begin{pmatrix}
&&x&0\\
&&0&y\\
x&0\\
0&y\\
&&&& D
\end{pmatrix},
\end{equation}
where $D\in M_{(n-2)\times(n-2)}$ is skew-symmetric. Notice that all but the $D$-part of this space is spanned by $q_1$ and $J_1$, so the $\sQ$-component of that matrix is a multiple of $q_1$. Other root spaces are given by
\begin{subequations}
\begin{align}
\sG_{(1,0)}&\leadsto W_i=E_{2i}+E_{4i}+E_{i2}-E_{i4}\in\sH,\\
\sG_{(-1,0)}&\leadsto Y_i=-E_{2i}+E_{4i}-E_{i2}-E_{i4},\\
\sG_{(0,1)}&\leadsto V_i=E_{1i}+E_{3i}+E_{i1}-E_{i3},\\
\sG_{(0,-1)}&\leadsto X_i=-E_{1i}+E_{3i}-E_{i1}-E_{i3}
\end{align}
\end{subequations}
with\footnote{Let us remember that we are dealing with $\SO(2,n)$ and that $AdS_{l}$ is a quotient of $\SO(2,l-1)$, so in the case of $AdS_{l}$ the index $j$ runs from $5$ to $l+1$. The first anti de Sitter space which contains such root spaces is $AdS_{4}$. More generally, remark that the table \eqref{EqTableSOIwa} of $\mathfrak{so}(2,n)$ gives the feeling that if something works with $AdS_4$, it will work for $AdS_{l\geq 4}$.} $\dpt{i}{5}{n+2}$. For example,
\begin{equation}
V_{5}=\begin{pmatrix}
&&&&1\\
&&&&0\\
&&&&1\\
&&&&0\\
1&0&-1&0&0
\end{pmatrix},
\qquad
W_5=
\begin{pmatrix}
&&&&0\\
&&&&1\\
&&&&0\\
&&&&1\\
0&1&0&-1&0
\end{pmatrix}
\end{equation}

\begin{equation}
\sG_{(1,1)}\leadsto M=
\begin{pmatrix}
   0&1&0&-1\\
   -1&0&1&0\\
   0&1&0&-1\\
   -1&0&1&0
\end{pmatrix},
\quad
\sG_{(1,-1)}\leadsto L=
\begin{pmatrix}
   0&1&0&-1\\
   -1&0&-1&0\\
   0&-1&0&1\\
   -1&0&-1&0
\end{pmatrix},
\end{equation}
\begin{equation}
\sG_{(-1,1)}\leadsto N=
\begin{pmatrix}
   0&1&0&1\\
   -1&0&1&0\\
   0&1&0&1\\
   1&0&-1&0
\end{pmatrix},
\quad
\sG_{(-1,-1)}\leadsto F=
\begin{pmatrix}
   0&1&0&1\\
   -1&0&-1&0\\
   0&-1&0&-1\\
   1&0&1&0
\end{pmatrix}.
\end{equation}
These are the same spaces as the previous ones. The subtlety is that we will choice an other notion of positivity, so that the space $\sN$ will be different. 

\begin{figure}[ht]
\begin{center}
\begin{pspicture}(-2,-2)(2,2)
  \psaxes[dotsep=1pt]{->}(0,0)(-1.9,-1.9)(1.9,1.9)
  \psdots[dotscale=2](1,0)(1,1)(1,-1)(0,1)
  \psdots[dotstyle=diamond,dotscale=2](-1,0)(-1,1)(-1,1)(0,-1)(-1,-1)
  \rput(2.3,-0.3){$\sA^*_{\sH}$}
\end{pspicture}
\end{center}
\caption{The root spaces}\label{fig:root}
\end{figure}
 
Let us recall the aim of our new decomposition: we want to have $\sR_{\sH}\subset\sR$. For this purpose, the equation \eqref{eq:re_N_H} gives us a constraint on the choice of the positivity notion on $\sA^*$. First we must have $\sG_{(0,1)}\subset\sN$. The the upper left $4\times 4$ corner of $\sN_{\sH}$ is spanned by $\sG_{(1,1)}-\sG_{(1,-1)}$. We complete our choice of $\sN$ with $\sG_{(1,0)}$. The underlying notion of positivity is that the element $\alpha(a,b)$ is positive\index{Positivity} in $\sA^*$ when $ (a>0) \LogOu(a=0 \LogEt b>0)$. 

The difference between decomposition and the previous one is the replacement of $N$ by $L\in\sG_{(1,-1)}$. Now\label{PgTablaIwa}
\begin{subequations}
\begin{align}
\sN&=\{W_i,V_j,M,L\}\\
\sA&=\{ J_1, J_2\},
\end{align}
\end{subequations}
with the commutator table 
\begin{subequations}  \label{EqTableSOIwa}
\begin{align}
[V_i,W_j]&=\delta_{ij}M &[V_j,L]&=2W_j\\
[ J_1,W_j]&=W_j       &[ J_2,V_i]&=V_i\\
[ J_1,L]&=L           &[ J_2,L]&=-L\\
[ J_1,M]&=M           &[ J_2,M]&=M.
\end{align}
\end{subequations}
It is important to note that $W_{i}$, $J_{1}\in\sH$ and $J_{2}\in\sQ$. The following change of basis in $\sA$ reveals to be useful in some circumstances:
\begin{subequations}  \label{EqChmHJ}
\begin{align}
  H_1&=J_1-J_2
	&H_2&=J_1+J_2
\end{align}
\end{subequations}
which leads to the table
\begin{subequations}		\label{TableSeconde}
\begin{align}
{}[V_i,W_j]&=\delta_{ij}M		&[V_j,L]&=2W_j\\
[H_1,V_i]&=-V_i				&[H_2,V_i]&=V_i\\
[H_1,W_i]&=W_i				&[H_2,W_i]&=W_i\\
[H_1,L]&=2L				&[H_2,M]&=2M.
\end{align}
\end{subequations}

\section[Low dimensional anti de Sitter spaces]{Group realization of low dimensional anti de Sitter spaces}

\subsection{Two dimensional anti de Sitter}		\label{SubsecTwoDimAdSAdGH}

Using notations and conventions of section \ref{SecToolSL}, the set $\Ad(G)H$ is a subset of $\gsl(2,\eR)$ made of elements of norm $8$. One can show by brute force computation or using the commutation relations that 
\[ 
\begin{split}
 \Ad( e^{x_KT} e^{x_NE})H&=\big( -\sin(2x_K)x_N+\cos(2x_K) \big)H\\
			&+\big( -\cos(2x_K)x_N-\sin(2x_K) \big)(E+F)\\
			&-x_NT,
\end{split}  
\]
which is a, following the Killing form \eqref{EqBHEFTsldR},  general element of norm $8$ in $\gsl(2,\eR)$. Thus as sets, $AdS_2$ is $\Ad(KN)H$.  Since $A$ is the stabilizer of $H$ for the adjoint action of $G$ on $H$ and $G=ANK=KNA$, we also have
\[ 
  AdS_2=G/A=\Ad(KN)H=\Ad(G)H.
\]
That provides isomorphisms
\begin{equation}
\begin{aligned}
 \phi\colon [0,\pi[\times\eR&\to AdS_2 \\ 
(x_K,x_B)&\mapsto \Ad( e^{xK T} e^{xNE})H,
\end{aligned}
 \end{equation}
or
\begin{equation}	 \label{EqCylAdSDeux}
\begin{aligned}
 \phi\colon Cyl&\to AdS_2 \\ 
(\theta,h)&\mapsto \Ad( e^{\frac{ \theta }{ 2 }T} e^{hE})H. 
\end{aligned}
\end{equation}
where $Cyl$ is the usual cylinder in $\eR^3$.

\subsection{Three dimensional anti de Sitter}		\label{SubsecGpAdsDeux}

The space $AdS_3$ is the hyperboloid 
\begin{equation} \label{hyperboloide}
 u^2 + t^2 - x^2 - y^2 = 1
 \end{equation}
embedded in $\eR^{2,2}$. There exists a bijection between $\eR^{2,2}$ and the two by two real matrices given by
\[ 
 v= \begin{pmatrix}
u\\t\\x\\y
\end{pmatrix}
\mapsto
g(v)=\begin{pmatrix}
u+x&y+t\\
y-t&u-x
\end{pmatrix}.
\]
As far as norm is concerned, we have $\| v \|=\det g(v)$. Among these matrices, the ones of $\SL(2,\eR)$ are given by the condition $\det g=1$, which is precisely the equation of the hyperboloid in $\eR^{4}$. That shows that $AdS_3=\SL(2,\eR)$.

\section{Symmetric space structure on anti de Sitter}\label{SecSymeStructAdS}

The $l$-dimensional anti de Sitter space $AdS_l$ can be described as set of points $(u,t,x_1,\ldots,x_{l-1})\in \eR^{2,l-1}$  such that $u^2+t^2-x_1^2-\ldots-x_{l-1}^2=1$. The next few pages are devoted to describe the homogeneous and symmetric space structures on $AdS_l$ induced by the transitive an isometric action of $\SO(2,l-1)$. We suppose that the groups $\SO(2,l-1)$ and $\SO(1,l-1)$ are parametrized in such a way that the second, seen as subgroup of the first one, leaves unchanged the vector $(1,0,\ldots,0)$. In this case, proposition 4.3 of chapter II in \cite{Helgason} provides the homogeneous space isomorphism
\begin{equation}
\begin{aligned}
  \SO(2,l-1)/\SO(1,l-1)&\to AdS_l \\ 
[g]&\mapsto  
 g\cdot
\begin{pmatrix}
1\\0\\\vdots
\end{pmatrix}
\end{aligned}
\end{equation}
where the dot denotes the usual ``matrix times vector'' action of the representative $g\in [g]$ in the defining representation of $\SO(2,l-1)$ on $\eR^{2,l-1}$. As far as notations are concerned, the classes are taken from the right:  $[g]=\{gh\tq h\in H\}$; in particular the class of the identity $e$ is denoted by $\mfo$; the groups $\SO(2,l-1)$ and $\SO(1,l-1)$ are denoted by $G$ and $H$ respectively and their Lie algebras by $\sG$ and $\sH$. Following proposition \ref{PropGHconn}, we can in fact only consider the identity components of $G$ and $H$.

\begin{proposition}
The homogeneous space $AdS_l$ is reductive\index{Reductive!$AdS_n$}.
\label{PropAdSreduct}
\end{proposition}

The proof relies on the following lemma and the fact that $\SO(2,n)$ is semisimple.
\begin{lemma} 
If $G$ is a semisimple Lie group and $H$ a semisimple subgroup of $G$, the restrictions on $H$ of the Killing form of $G$ is nondegenerate.
 \label{lem:Killing_ss_descent}
\end{lemma}

\begin{proof}[Proof of proposition \ref{PropAdSreduct}]
From the Killing form of $G$ , one defines
\[
   \sQ=\sH^{\perp}=\{X\in\sG:B(X,H)=0\,\forall H\in\sH\}.
\]
Let $H$, $H'\in\sH$ and $Y\in\sQ$. From $\ad$-invariance of the Killing form, we have $B([H,Y],H')=0$. Hence $(\ad(\sH)\sQ)\subset \sQ$ and the claim is proved.

\end{proof}

Matrices of $\SO(2,n)$ are $(2+n)\times(2+n)$ matrices while the $n$-dimensional anti de Sitter space is a quotient of $\SO(2,n-1)$. In order to avoid confusions, we will reserve the letter $n$ to the study of the group $\SO(2,n)$ and the letter $l$ will denote the dimension of the anti de Sitter space which will thus be $AdS_l$.

Let us provide a matrix representation now. The matrices of $\so(1,n)$ have to be seen as matrices of $\so(2,n)$ with the condition $Y^t\sigma+\sigma Y=0$ for the  ``metric''\ $\sigma=diag(0,-,+,\ldots,+)$. Hence,
\begin{equation}\label{eq:gene_H}
\sH=\soun\leadsto
  \begin{pmatrix}
     \begin{matrix}
       0&0\\
       0&0
     \end{matrix}
                       &  \begin{pmatrix}
		             \cdots 0\cdots\\
			    \leftarrow v^t\rightarrow
                          \end{pmatrix}\\
    \begin{pmatrix}	  
       \vdots & \uparrow\\
         0    & v \\
       \vdots & \downarrow
    \end{pmatrix} &  B
  \end{pmatrix}
\end{equation}
where  $v\in M_{n\times 1}$ and $B\in M_{n\times n}$ is skew-symmetric. Comparing this with the general form \eqref{eq:gene_sodn} of a matrix of $\sodn$ matrix, one immediately finds that, with the choice
\begin{equation}\label{eq:gene_M}
\sQ\leadsto
 \begin{pmatrix}
     \begin{matrix}
       0&a\\
       -a&0
     \end{matrix}
                       &  \begin{pmatrix}		             
			  \leftarrow w^t\rightarrow \\
			     \cdots 0\cdots\\
                          \end{pmatrix}\\
    \begin{pmatrix}	  
      \uparrow   & \vdots\\
          w      &  0\\
      \downarrow & \vdots 
    \end{pmatrix} & 0
  \end{pmatrix},
 \end{equation}
the decomposition $\sG=\sH\oplus\sQ$ is reductive:
\begin{align}\label{EqDefRedHQ}
  [\sH,\sQ]&\subseteq\sQ,
 &[\sQ,\sQ]&\subseteq\sH,
\end{align}
and $B(\sH,\sQ)=0$. In the sequel, we will use the basis of $\sQ$ defined by 
\begin{align}		\label{EqDefBaseqi}
  q_0&=E_{12}-E_{21}, &q_i&=E_{1i}+E_{i1}.
\end{align}
We define the involutive automorphism $\sigma=\id|_{\sH}\oplus(-\id)|_{\sQ}$.  The vector space $\sQ$ can be identified with the tangent space $T_{[e]}AdS_l$, and that identification can be extended by defining $\sQ_g=dL_g\sQ$. In this case $\dpt{d\pi}{\sQ_g}{T_{[g]}AdS_l}$ is a vector space isomorphism.\label{PgdpibaseQTgM} A homogeneous metric on $T_{[g]}AdS_l$ is defined as in subsection \ref{SubsecKillHomo}.

  Cartan decomposition of $\SO(2,l-1)$ are of crucial importance in chapter~\ref{ChapAdS}, so that we want to use a Cartan involution $\theta$ such that $[\sigma,\theta]=0$ (see \cite{Loos} page 153, theorem 2.1). One can show that $X\mapsto -X^t$ has that property. The corresponding Cartan decomposition is described in appendix \ref{SubSecCartandeuxN}.

As a consequence of relations \eqref{EqDefRedHQ}, 
\begin{equation}  \label{EqdpiAdpi}
d\pi\Ad(h)=\Ad(h) d\pi
\end{equation}
because, if $X\in\sQ$, $d\pi^{-1}(X)=\{ X+Y\tq Y\in\sH \}$, so $\Ad(h)Y\in\sH$ and $\Ad(h)X\in \sQ$.

\section{Iwasawa decomposition for \texorpdfstring{$\gsl(2,\eC)$}{sl2C}}
\index{Iwasawa decomposition!of $\SL(2,\eC)$}

Matrices of $\gsl(2,\eC)$ are acting on $\eC^2$ as 
\[ 
\begin{split}
  \begin{pmatrix}
\alpha&\beta\\\gamma&-\alpha
\end{pmatrix}&
\begin{pmatrix}
a+bi\\c+di
\end{pmatrix}\\
&=
\begin{pmatrix}
(\alpha_1a-\alpha_2b+\beta_1c-\beta_2d)+i(\alpha_2a+\alpha_1b+\beta_2c+\beta_1d)\\
(\gamma_1a-\gamma_2b-\alpha_1c+\alpha_2d)+i(\gamma_2a+\gamma_1b-\alpha_2c-\alpha_1d)
\end{pmatrix}
\end{split}
\]
if $\alpha=\alpha_1+i\alpha_2$.  Our aim is to embed $\SL(2,\eC)$ in $\SP(2,\eR)$ (see sections \ref{SecSympleGp} and \ref{SecDirADs}), so that we want a four dimensional realization of $\gsl(2,\eC)$. It is easy to rewrite the previous action under the form of $\begin{pmatrix}
\alpha&\beta\\\gamma&-\alpha
\end{pmatrix}$ acting of the vertical four component vector $(a,b,c,d)$. The result is that a general matrix of $\gsl(2,\eC)$ reads
\begin{equation}		\label{EqGenslMatr}
\gsl(2,\eC)\leadsto
\begin{pmatrix}
\boxed{
\begin{array}{cc}
\alpha_1&-\alpha_2\\
\alpha_2&\alpha_1
\end{array}
}&
\begin{array}{cc}
\beta_1&-\beta_2\\
\beta_2&\beta_1
\end{array}\\
\begin{array}{cc}
\gamma_1&-\gamma_2\\
\gamma_2&\gamma_1
\end{array}&
\boxed{
\begin{array}{cc}
-\alpha_1&\alpha_2\\
-\alpha_2&-\alpha_1
\end{array}
}
\end{pmatrix}.
\end{equation}
The boxes are drawn for visual convenience.  Using the Cartan involution $\theta(X)=-X^t$, we find the following Cartan decomposition:
\begin{equation}
\begin{split}
\iK_{\gsl(2,\eC)}&\leadsto
\begin{pmatrix}
\boxed{
\begin{array}{cc}
0&-\alpha_2\\
\alpha_2&0
\end{array}
}&
\begin{array}{cc}
\beta_1&-\beta_2\\
\beta_2&\beta_1
\end{array}\\
\begin{array}{cc}
-\beta_1&-\beta_2\\
\beta_2&-\beta_1
\end{array}&
\boxed{
\begin{array}{cc}
0&\alpha_2\\
-\alpha_2&0
\end{array}
}
\end{pmatrix},\\
\iP_{\gsl(2,\eC)}&\leadsto
\begin{pmatrix}
\boxed{
\begin{array}{cc}
\alpha_1&0\\
0&\alpha_1
\end{array}
}&
\begin{array}{cc}
\beta_1&-\beta_2\\
\beta_2&\beta_1
\end{array}\\
\begin{array}{cc}
-\beta_1&-\beta_2\\
\beta_2&-\beta_1
\end{array}&
\boxed{
\begin{array}{cc}
0&\alpha_2\\
-\alpha_2&0
\end{array}
}
\end{pmatrix}.
\end{split}
\end{equation}
We have $\dim\iP_{\gsl(2,\eC)}=3$ and $\dim\iP_{\gsl(2,\eC)}=3$. A maximal abelian subalgebra of $\iP_{\gsl(2,\eC)}$ is the one dimensional algebra generated by
\[ 
  A_1=
\begin{pmatrix}
1\\&1\\&&-1\\&&&-1
\end{pmatrix}.
\]
The corresponding root spaces are
\begin{itemize}
\item $\gsl(2,\eC)_0$:
\[ 
  I_1=
\begin{pmatrix}
1\\&1\\&&-1\\&&&-1
\end{pmatrix},\quad
I_2=
\begin{pmatrix}
0&-1\\
1&0\\
&&0&1\\
&&-1&0
\end{pmatrix}
\]
\item $\gsl(2,\eC)_2$:
\[ 
  D_1=\begin{pmatrix}
&&1&0\\
&&0&1\\
0&0\\
0&0
\end{pmatrix},\quad
D_2=
\begin{pmatrix}
&&0&-1\\
&&1&0\\
0&0&\\
0&0&
\end{pmatrix}
\]
\item $\gsl(2,\eC)_{-2}$
\[ 
  C_1=\begin{pmatrix}
&&0&0\\
&&0&0\\
1&0\\
0&1
\end{pmatrix},\quad
C_2=\begin{pmatrix}
&&0&0\\
&&0&0\\
0&-1\\
1&0
\end{pmatrix}.
\]
\end{itemize}
It is natural to choice $\gsl(2,\eC)_2$ as positive root space system. In this case, $\iN_{\gsl(2,\eC)}=\{ D_1,D_2 \}$, $\iA_{\gsl(2,\eC)}=\{ I_1 \}$ and the table of $\iA\oplus\iN$ is
\begin{align}
[I_1,D_1]&=2D_1&		[D_1,D_2]&=0\\
[I_1,D_2]&=2D_2&		
\end{align}

\section{Symplectic group}		\label{SecSympleGp}

\subsection{Iwasawa decomposition}
\index{Iwasawa decomposition!of $\SP(2,\eR)$}

A simple computation shows that $4\times 4$ matrices subject to $A^t\Omega+\Omega A=0$ are given by
\[ 
  \begin{pmatrix}
A&B\\
C&-A^t
\end{pmatrix}
\]
where $A$ is any $2\times 2$ matrix while $B$ and $C$ are symmetric matrices. Looking at general form \eqref{EqGenslMatr}, we see that the operation to invert the two last column and then to invert the two last lines provides a homomorphism $\phi\colon \gsl(2,\eC)\to \gsp(2,\eR)$. The aim is now to build an Iwasawa decomposition of $\gsp(2,\eR)$ which ``contains'' the one of $\gsl(2,\eC)$.

Using the Cartan involution $\theta(X)=-X^t$, we find the Cartan decomposition
\begin{align}
\iK_{\gsp(2,\eR)}&\leadsto
\begin{pmatrix}
A&S\\-S&A
\end{pmatrix},
&\iP_{\gsp(2,\eR)}&\leadsto
\begin{pmatrix}
S&S'\\S'&-S
\end{pmatrix}
\end{align}
where $S$ and $S'$ are any symmetric matrices while $A$ is a skew-symmetric one. We have $\dim\iK_{\gsp(2,\eR)}=4$ and $\dim\iP_{\gsp(2,\eR)}=6$. It turns out that $\phi(\iK_{\gsl(2,\eC)})\subset\iK_{\gsp(2,\eR)}$ and $\phi(\iP_{\gsl(2,\eC)})\subset \iP_{\gsp(2,\eR)}$. A maximal abelian subalgebra of $\iP_{\gsp(2,\eR)}$ is spanned by the matrices $A'_1$ and $A'_2$ listed below and the corresponding root spaces are:
\begin{itemize}
\item $\gsp(2,\eR)_{(0,0)}$:
\[ 
  A'_1=
\begin{pmatrix}
1&0\\
0&1\\
&&-1&0\\
&&0&-1
\end{pmatrix},
\quad
A'_2=
\begin{pmatrix}
0&1\\
1&0\\
&&0&-1\\
&&-1&0
\end{pmatrix}
\]
\item $\gsp(2,\eR)_{(0,2)}$:
\[ 
 X'= \begin{pmatrix}
1&-1&\\
1&-1&\\
&&-1&-1\\
&&1&-1
\end{pmatrix}
\]
\item $\gsp(2,\eR)_{(0,-2)}$:
\[ 
 V'= \begin{pmatrix}
1&1\\
-1&-1\\
&&-1&1\\
&&-1&1
\end{pmatrix}
\]
\item $\gsp(2,\eR)_{(2,0)}$:
\[ 
 W'= \begin{pmatrix}
&&1&0\\
&&0&-1\\
0&0\\0&0
\end{pmatrix}
\]
\item $\gsp(2,\eR)_{(2,2)}$:
\[ 
  L'=
\begin{pmatrix}
&&1&1\\
&&1&1\\
0&0\\0&0
\end{pmatrix}
\]
\item $\gsp(2,\eR)_{(2,-2)}$:
\[ 
  M'=
\begin{pmatrix}
&&1&-1\\
&&-1&1\\
0&0\\0&0
\end{pmatrix}
\]
\item $\gsp(2,\eR)_{(-2,0)}$
\[ 
Y'=
\begin{pmatrix}
&&0&0\\&&0&0\\
1&0\\0&-1
\end{pmatrix}
\]
\item $\gsp(2,\eR)_{(-2,2)}$:
\[ 
  N'=\begin{pmatrix}
&&0&0\\&&0&0\\
1&-1\\
-1&1
\end{pmatrix}
\]
\item $\gsp(2,\eR)_{(-2,-2)}$:
\[ 
  F'=\begin{pmatrix}
&&0&0\\
&&0&0\\
1&1\\
1&1
\end{pmatrix}
\]
\end{itemize}
It is important to notice how do the root spaces of $\gsl(2,\eC)$ embed:
\begin{align}
\phi(I_1)&=A'_1	&\phi(I_2)&=\frac{ V'-X' }{ 2 }\\
\phi(D_1)&=\frac{ L'-M' }{2}	&\phi(D_2)&=-W'\\
\phi(C_1)&=\frac{ F'-N' }{2}	&\phi(C_2)&=Y'. 
\end{align}
So $\iN_{\gsp(2,\eR)}$ must at least contain the elements $L'$, $M'$ and $W'$. We complete the notion of positivity by $V'$. The Iwasawa algebra reads
\[ 
\begin{split}
\iA_{\gsp(2,\eR)}&=\{ B_1,B_2 \}\\
\iN_{\gsp(2,\eR)}&=\{ L',M',W',V' \}
\end{split}  
\]
with
\begin{align*}
[L',V']&=-4W'	&[W',V']&=-2M'\\
[B_1',L']&=2L'	&[B'_2,M']&=2M'\\
[B_1',W']&=W'	&[B_2',W']&=W'\\
[B_1',V']&=-V'	&[B_2',V']&=V'
\end{align*}
where $B'_1=\frac{ 1 }{2}(A'_1+A'_2)$ and $B_2=\frac{ 1 }{2}(A_1'-A_2')$. Generators of $\iK_{\gsp(2,\eR)}$ by
\begin{align*}
K'_t&=
\begin{pmatrix}
&&1&0\\
&&0&1\\
-1&0\\
0&-1
\end{pmatrix}
	&K'_1&=
\begin{pmatrix}
0&1\\-1&0\\
&&0&1\\
&&-1&0
\end{pmatrix}\\
K'_2&=
\begin{pmatrix}
&&0&1\\&&1&0\\0&-1\\-1&0
\end{pmatrix}
	&K'_3&=
\begin{pmatrix}
&&1&0\\
&&0&-1\\
-1&0\\
0&1
\end{pmatrix}.
\end{align*}
Notice that $[K'_t,K'_i]=0$ for $i=1$, $2$, $3$.

\subsection{Isomorphism}		\label{SubSecIsosp}

The following provides an isomorphism $\psi\colon \so(2,3)\to \gsp(2,\eR)$:
\begin{align*}
\psi(H_i)&=B'_i		&\psi(u)&=K'_t\\
\psi(W)&=W'		&\psi(R_1)&=\frac{ 1 }{2}K'_1\\
\psi(M)&=M'		&\psi(R_2)&=\frac{ 1 }{2}K'_2\\
\psi(L)&=L'		&\psi(R_3)&=\frac{ 1 }{2}K'_3\\
\psi(V)&=\frac{ 1 }{2}V'
\end{align*}
where the $R_i$'s are the generators of the $\so(3)$ part of $\sK_{\so(2,3)}$ satisfying the relations $[R_i,R_j]=\epsilon_{ijk}R_k$. It is now easy to check that the image of the embedding $\phi\colon \gsl(2,\eC) \to \gsp(2,\eR)$ is exactly $\so(1,3)$, so that
\begin{equation}
\psi^{-1}\circ\phi\colon \gsl(2,\eC)\to \sH
\end{equation}
is an isomorphism which realises $\sH$ as subalgebra of $\gsp(2,\eR)$. This circumstance will be useful in defining a spin structure on $AdS_4$.

One can prove that the kernel of the adjoint representation of $\SP(2,\eR)$ on its Lie algebra is $\pm\mtu$, in other words, $\Ad(a)=\id$ if and only if $a=\pm\mtu$. We define a bijective map $h\colon \SO(2,3)\to \SP(2,\eR)/\eZ_2$ by the requirement that
\begin{equation}		\label{Eqdefhspsl}
  \psi\big( \Ad(g)X \big)=\Ad\big( h(g) \big)\psi(X)
\end{equation}
for every $X\in\so(2,3)$. The following is true for all $\psi(X)$:
\[ 
\begin{split}
\Ad\big(h(gg'\big)) \psi(X)&=\psi\Big( \Ad(g)\big( \Ad(g')X \big) \Big)\\
			&=\Ad\big( h(g) \big)\psi\big( \Ad(g')X \big)\\
			&=\Ad\big( h(g)h(g') \big)\psi(X),
\end{split}
\]
 the map $h$ is therefore a homomorphism. If an element $a\in \SP(2,\eR)$ reads $a= e^{X_A} e^{X_N} e^{X_K}$ in the Iwasawa decomposition, the property $\Ad(a)\psi(X)=\psi\big( \Ad(g)X \big)$ holds for the element\label{PgSolhpsiSP} $g= e^{\psi^{-1}X_A} e^{\psi^{-1}X_N} e^{\psi^{-1}X_K}$ of $\SO(2,3)$. This shows that $h$ is surjective.

\subsection{Reductive structure on the symplectic group}		\label{SubSecRedspT}

A lot of structure of $\so(2,3)$, such as the reductive homogeneous space decomposition as $\sQ\oplus\sH$, can be immediately transported from $\so(2,3)$ to $\gsp(2,\eR)$. Indeed, let $\mT=\psi(\sQ)$ and $\mI=\phi\big( \gsl(2,\eC) \big)$. We have the direct sum decomposition 
\[ 
\gsp(2,\eR)=\mT\oplus\mI.
\]
 Let $X\in\mT\cap\mI$, then $\psi^{-1}X$ belongs to $\sQ\cap\sH$ which only contains $0$. The fact that $\psi$ is an isomorphism yields that $X=0$. Since $\psi$ preserves linear independence, a simple dimension counting shows that the sum actually spans the whole space.

Putting $g=h^{-1}(a)$ in the definition \eqref{Eqdefhspsl} of $h$, we find
\[ 
  \psi\left( \Ad\big( h^{-1}(a) \big)X \right)=\Ad(a)\psi(X).
\]
Considering a path $a(t)$ with $a(0)=e$, we differentiate this expression with respect to $t$ at $t=0$ we find
\[ 
  \ad(dh^{-1}\dot a)X=d\psi^{-1}\big( \ad(\dot a)\psi(X) \big)=\ad(d\psi^{-1}\dot a)(d\psi^{-1}\psi X),
\]
but $d\psi=\psi$ because $\psi$ is linear, hence $[dh^{-1}\dot a,X]=[\psi^{-1}\dot a,X]$ for all $X\in \so(2,3)$ and $\dot a\in \gsp(2,\eR)$. We deduce that $(dh^{-1})_e=\psi^{-1}$. We define 
\begin{align*}
	\theta_{\gsp}&=\id|_{\iK_{\gsp}}\oplus(-\id)|_{\iP_{\gsp}}\\
	\sigma_{\gsp}&=\id|_{\mT}\oplus(-\id)|_{\mI}.
\end{align*}
We can check that $\psi^{-1}\circ\theta_{\gsp}\circ\psi=\theta$ and $\psi^{-1}\circ\theta_{\gsp}\circ\psi=\theta$. Then it is clear that
\[ 
  [\sigma_{\gsp},\theta_{\gsp}]=0
\]
using the corresponding vanishing commutator in $\so(2,3)$. We denote $\mT_a=dL_a\mT$ and the fact that $dp= d\pi\circ dh^{-1}= d\pi\circ \psi^{-1}$ shows that $dp(\mT_a)$ is a basis of $T_{p(a)}(G/H)$. So we consider the basis $t_i=\psi(q_i)$ of $\mT$ and the corresponding left-invariant vector fields $\tilde t_i(a)=dL_at_i$.
\section{Some symplectic and Poisson geometry}\label{sec:symple}

\subsection{Symplectic manifold}

A \defe{symplectic structure}{Symplectic!vector space} on a vector space $V$ is a skew-symmetric, nondegenerate bilinear $2$-form $\dpt{\Omega}{V\times V}{\eR}$. We define the \defe{symplectic group}{Symplectic!group} $\SP(\Omega)$\nomenclature{$\SP(V,\Omega)$}{Symplectic group} as the group of linear operators $\dpt{A}{V}{V}$ such that $\Omega(Au,Av)=\Omega(u,v)$ for every $u$, $v\in V$. It is easy to see that elements of $\SP(V)$ satisfy
\begin{equation} \label{EqPropodefMtrsym}
   A^t\Omega A=\Omega.
\end{equation}
The Lie algebra of $\SP(\Omega)$ is denoted by $\gsp(\Omega)$\nomenclature{$\mfsp$}{Symplectic algebra}. Taking the derivative of equation \eqref{EqPropodefMtrsym} with respect to $A$, one finds the following condition for $B\in\gsp(\Omega)$:
\begin{equation}
 \Omega B+B^{t}\Omega=0.
\end{equation}

A \defe{symplectic manifold}{Symplectic!manifold} is the data of a smooth manifold $M$ and a symplectic structure $\omega_{x}$ on each tangent space $T_{x}M$. The map $x\mapsto \omega_{x}$ is required to be a smooth section of the $2$-tensor bundle.

\begin{definition} 
A \defe{symplectic Lie algebra}{Symplectic!Lie algebra} is a Lie algebra $\mfs$ endowed with a symplectic structure $\omega$  such that $\forall x$, $y$, $z\in\mfs$,
\begin{equation}   \label{EqDefAlgSymple}
\omega([x,y],z)+\omega([y,z],x)+\omega([z,x],y)=0.
\end{equation}
\label{DefSympleAlg}
\end{definition}
\subsection{Poisson manifold}

Let $M$ be a smooth manifold. A \defe{Poisson bracket}{Poisson bracket}, or a \emph{Poisson structure} on $M$ is a map $\{ .,. \}\colon  C^{\infty}(M)\times C^{\infty}(M)\to  C^{\infty}(M)$ such that
\begin{enumerate}
\item $\big(  C^{\infty}(M), \{ .,. \} \big)$ is a Lie algebra,
\item for each $f\in C^{\infty}(M)$, the map $\{ f,. \}$ is a derivation of the algebra $ C^{\infty}(M)$:
\[ 
  \{ f,gh \}=\{ f,g \}h+g\{ f,h \}.
\]
\end{enumerate}

\subsection{Hamiltonian action}\label{app:ham_act}

Let $(M_1,\omega_1)$ and $(M_2,\omega_2)$ be symplectic manifolds. A \defe{symplectomorphism}{Symplectomorphism} from $M_1$ to $M_2$ is a diffeomorphism $\varphi\colon M_1\to M_2$ such that $\varphi^*\omega_2=\omega_1$.

For any function $f\in C^{\infty}(M)$, we define the \defe{Hamiltonian field}{Hamiltonian!field} $X_f\in\cvec(M)$ associated with $f$ by
\begin{equation}   \label{EqDefHamVect}
  i(X_f)\omega=df.
\end{equation}
 Existence is assured because $\omega$ is nondegenerate.  A symplectic structure induces a Poisson bracket by defining\index{Poisson bracket!on symplectic manifold}
\begin{equation}\label{eq:def_Poisson}
   \{f,g\}=-\omega(X_f,X_g)=-X_g(f)=X_f(g).
\end{equation}
In local coordinates, one can write $\omega=\frac{1}{2}\omega_{ij}dx^i\wedge dx^j$ and $X_f=\omega^{ij}\partial_if\partial_j$, where $(\omega^{ij})$ is the inverse matrix of $(\omega_{ij})$. The Poisson tensor defined by
\begin{equation}  \label{EqPoisson} 
   \{f,g\}=P^{kl}\partial_kf\partial_lg,
\end{equation} 
is nothing else than $P=\omega^{-1}$.

\begin{theorem}

If $\varphi\colon M\to M'$ is a diffeomorphism between two Poisson manifolds $(M,P)$ and $(M',P')$, then the following are equivalent:
\begin{enumerate}
\item $\varphi_*(X_{f\circ\varphi})=X'_f$,
\item\label{ite_equivii} $\{u\circ\varphi,v\circ\varphi\}=\{u,v\}'\circ\varphi$,
\item $\varphi_*P=P'$,
\end{enumerate}

\setcounter{bidon}{\value{enumi}}
\noindent If moreover the Poisson structures $P$ and $P'$ come from symplectic forms $\omega$ and $\omega'$,
\begin{enumerate}
\setcounter{enumi}{\value{bidon}}
\item $\varphi^*\omega'=\omega$.
\end{enumerate}
\label{tho:equiv_Poisson}
\end{theorem}

Now, we consider a symplectic action $\dpt{\tau}{G\times M}{M}$ of a Lie group $G$ on $M$ (i.e. $\dpt{\tau_g}{M}{M}$ is a symplectic transformation of $M$ for each $g\in G$).  The action is \defe{Hamiltonian}{Hamiltonian!action} if, for every $X\in\mG$, there exists a map $\lambda_X\in C^{\infty}(M,\eC)$ such that
\begin{subequations}\label{eq:act_ham}
\begin{align} 
  i(X^*)\omega&=d\lambda_X,\label{1s17d}\\
  \{\lambda_X,\lambda_Y\}&=\lambda_{[X,Y]}.\label{eq:hamil}
 \end{align}
\end{subequations}

\begin{definition}
The map $\dpt{\lambda}{\mG}{\Cinf(M)}$ which satisfies \eqref{eq:act_ham} is the \defe{dual}{Momentum map!dual} momentum map while the \defe{momentum map}{Momentum map} is $\dpt{J}{M}{\mG^*}$ defined by
   \begin{equation} \label{eq:defmomm ap}
     \lambda_X(x)=\langle J(x),X\rangle
   \end{equation}
for all $X\in\mG$.
\label{def:app_mom_mom_duale}
\end{definition}

\subsection{Coadjoint orbits}		\label{sub:coadjoint}

Let $G$ be a Lie group and $\mG$ its Lie algebra.  We know that $G$ acts on the dual $\mG^*$ by
\begin{equation}    \label{EqDefActCoadj}
  g\cdot \xi=\xi\circ \Ad(g^{-1})=\Ad(g)^*\xi
\end{equation}
for $g\in G$ and $\xi\in\mG^*$. The second equality defines the \defe{coadjoint action}{Coadjoint!action} $\Ad^*\colon G\times\mG^*\to \mG^*$. In other words, for all $X\in\mG$,
\[
(g\cdot \xi)(X)=\langle \xi,\Ad(g^{-1})X\rangle=\langle \Ad(g)^*(\xi),X\rangle.
\]
In this context, the notion of fundamental fields is given bu $X^*_{\xi}=\xi\circ\ad(X)$. Let $\theta_{\xi}=\{g\cdot \xi|g\in G\}$, the orbit of $\xi$ in $\mG^*$. It can be shown that 
\begin{equation}\label{eq:omega_Gs}
  \tilde\omega_x(X^*_x,Y^*_x)=\langle x,[X,Y]\rangle
\end{equation}
 defines a symplectic form on $\theta_{\xi}$\nomenclature{$\theta_{\xi}$}{Coadjoint orbit of $\xi$}, the coadjoint orbit of $\xi$.

\begin{proposition}
The coadjoint action is Hamiltonian.
\end{proposition}

\subsection{Central extension}

Let $\mG$ be a Lie algebra. A \defe{Chevalley coboundary}{Chevalley coboundary} is a $2$-form which reads $\delta\xi$ for a certain $\xi\in\mG^*$ with $\delta$ defined by
\begin{equation}  \label{EqDefChevCoycl}
   (\delta\xi)(A,B)=-\xi([A,B]).  
\end{equation}
Let $\Omega$ be a $2$-cocycle. If it is not a coboundary, we add an element $C$ in $\mG$ and we consider $\mG'=\mG\oplus\eR C$ with the Lie algebra structure
\begin{equation}
  [A+s,B+t]_{\mG'}=[A,B]_{\mG}+\Omega(A,B)C.
\end{equation}
This is the \defe{central extension}{Central extension} of $\mG$ with respect to the $2$-cocycle $\Omega$. The terminology comes from the fact that the extension $\eR C$ belongs to the center of $\mG'$. The point is that $\Omega$ is a coboundary in $\mG'$ because 
\begin{equation}
\begin{split}
(\delta C^*)(A,B)=C^*[A,B]_{\mG'}
		=C^*\big( [A,B]_{\mG}+\Omega(A,B)C \big)
		=\Omega(A,B),
\end{split}
\end{equation}
so that $\Omega=\delta C^*$.

Now we suppose that the group $G$ acts on a manifold $M$. We define the action of the extended group $G'=G\otimes e^{\eR C}$ by saying that the ``new'' part does not act: $(g,s)\cdot x=g\cdot x$. Fundamental fields remains unchanged:
\begin{equation}  \label{eq:XseqXss}
(X,s)^*=X^*.
\end{equation}
 If the action of $G$ on $M$ is weakly Hamiltonian, we have functions $\dpt{ \mu_x}{M}{\eC}$ such that $i(X^*)\omega=d\mu_X$. These functions fulfil $X^*=\{ \mu_X,\,. \}$. We define 
\begin{equation} \label{eq:lamXsmuXs}
  \lambda_{X,s}=\mu_X+s.
\end{equation}

\begin{proposition}
The action of $G'$ is (strongly) Hamiltonian for these functions. 
\end{proposition}

\begin{proof}

From equation \eqref{eq:XseqXss}, we have $\{ \mu_X,. \}=\{ \mu_{X,s},. \}$ hence
 \begin{equation}
  \{ \lambda_{(X,s)},\lambda_{(Y,t)} \}=\{ \mu_X,\mu_Y \}
		=\mu_{[X,Y]}+C_{X,Y}
\end{equation}
for certain constants $C_{XY}$ which satisfy the property $d\big( \{ \mu_X,\mu_Y \}-\mu_{[X,Y]} \big)=0$.
Therefore
\begin{equation}
  \{ \lambda_{(X,s)},\lambda_{(Y,t)} \}=\lambda_{[X,Y],C_{X,Y}}
			=\lambda_{[(X,s),(Y,t)]}.
\end{equation}

\end{proof}

The sense of the whole construction is the following. When the action $G$ is weakly Hamiltonian on $M$, we have functions $\mu_X$ which define $\Omega$ by
\[ 
  \{ \mu_X,\mu_Y \}=\mu_{[X,Y]}+\Omega(X,Y).
\]
In this case, the corresponding group extension has a strongly Hamiltonian action with momentum maps given by \eqref{eq:lamXsmuXs}.

\pagestyle{empty}

\chapter*{Conclusion}
\addcontentsline{toc}{chapter}{Conclusion}

In a first time we defined a black hole in anti de Sitter space. This construction is not related to any metric divergence but is a dimensional generalization of a causal black hole whose singularity is dictated by causal issues. The originality of our approach lies in the fact that our method uses essentially group theoretical and symmetric spaces techniques. That result should be generalisable to any semisimple symmetric space.

Then we proved that the physical domain of the black hole (the non singular part) is equivalent to a group in the sense that there exists a group which acts freely and transitively by diffeomorphisms. So we identify the group with the manifold and it is easy to prove that the latter group is a split extension of an Heisenberg group which happens to be quantizable by a twisted pull-back of a previously known quantization of $\SU(1,n)$. 

We also proved two somewhat out of subject small results. The first one is the fact that a deformation of the half-plane by Unterberger can be transported to a deformation of the Iwasawa subgroup of $\SL(2,\eR)$ which can in turn deform (by the group action method) the dual of its Lie algebra. We  showed however that that deformation is not universal; indeed we pointed out two different actions of the Iwasawa subgroup of $\SL(2,\eR)$ on $AdS_2$ for which the deformation by group action method reveals to be unable to even multiply two compactly supported functions. An interesting question is to know the precise point in the construction of Unterberger which makes his product non universal.

The second small result is a proof of concept for quantization of the Iwasawa subgroup of $\SO(2,n)$ by the method of the extension lemma. We wrote $\SO(2,n)$ as a symplectic split extension of $\SU(1,n)$ by $\SU(1,1)$. The extension lemma then provided a kernel on $\SO(2,n)$ because kernels were known on $\SU(1,1)$ and $\SU(1,n)$. Is that quantization equivalent in some sense to the one that we performed in the main line of the black hole deformation ? That question still has to be solved.

As a final remark, I want to point out that the major challenge of this century is not quantization, but global warming.

\bibliographystyle{unsrt}			
\def\polhk#1{\setbox0=\hbox{#1}{\ooalign{\hidewidth
  \lower1.5ex\hbox{`}\hidewidth\crcr\unhbox0}}}

\begin{theindex}

  \item Associativity
    \subitem of a WKB quantization, \hyperpage{61}
  \item Automorphism
    \subitem of affine symplectic manifold, \hyperpage{59}

  \indexspace

  \item Black hole, \hyperpage{15}, \hyperpage{18}
  \item BTZ black hole, \hyperpage{14}

  \indexspace

  \item Cartan
    \subitem decomposition, \hyperpage{93}
    \subitem involution, \hyperpage{93}
  \item Causal solvable symmetric black hole, \hyperpage{14}
  \item Causal structure, \hyperpage{16}
  \item Central extension, \hyperpage{116}
  \item Chevalley coboundary, \hyperpage{115}
  \item Coadjoint
    \subitem action, \hyperpage{115}
    \subitem orbit, \hyperpage{88}
  \item Covariant star product, \hyperpage{63}

  \indexspace

  \item Decomposition
    \subitem Cartan, \hyperpage{93}
    \subitem Iwasawa, \hyperpage{94}
    \subitem root space, \hyperpage{93}
  \item Deformation, \hyperpage{59}
  \item Differentiable vector, \hyperpage{62}

  \indexspace

  \item Event horizon, \hyperpage{18}
  \item Extension
    \subitem lemma, \hyperpage{78}
    \subitem of Heisenberg algebra, \hyperpage{81}
  \item Exterior, \hyperpage{18}

  \indexspace

  \item Fundamental
    \subitem vector field
      \subsubitem on $\eR^2$, \hyperpage{68}
  \item Fundamental vector field, \hyperpage{95}
  \item Future directed vector, \hyperpage{18}

  \indexspace

  \item Globally group type, \hyperpage{30}

  \indexspace

  \item Hamiltonian
    \subitem action, \hyperpage{115}
    \subitem field, \hyperpage{114}
  \item Heisenberg algebra, \hyperpage{81}
  \item Homogeneous
    \subitem space, \hyperpage{95}
    \subitem space reductive, \hyperpage{95}

  \indexspace

  \item Interior region, \hyperpage{18}
  \item Iwasawa decomposition, \hyperpage{94}
    \subitem of $\SL(2,\eC)$, \hyperpage{109}
    \subitem of $\SO(1,n)$, \hyperpage{99}
    \subitem of $\SO(2,n)$, \hyperpage{100}
    \subitem of $\SP(2,\eR)$, \hyperpage{110}
    \subitem of $SL(2,\eR)$, \hyperpage{96}

  \indexspace

  \item Kernel
    \subitem for extension $(\BX ,0,2)$ of Heisenberg, \hyperpage{87}
  \item Kernel of a WKB quantization, \hyperpage{61}

  \indexspace

  \item Left-invariant
    \subitem kernel, \hyperpage{61}
    \subitem vector field, \hyperpage{94}
  \item Light
    \subitem cone, \hyperpage{17}
    \subitem ray, \hyperpage{16}
  \item Locally group type, \hyperpage{30}

  \indexspace

  \item Momentum map, \hyperpage{115}
    \subitem dual, \hyperpage{115}

  \indexspace

  \item Natural topology, \hyperpage{95}

  \indexspace

  \item Partial Fourier transform, \hyperpage{74}
  \item Poisson bracket, \hyperpage{114}
    \subitem and Moyal product, \hyperpage{71}
    \subitem on symplectic manifold, \hyperpage{114}
  \item Positivity, \hyperpage{94}, \hyperpage{105}

  \indexspace

  \item Quantization, \hyperpage{59}

  \indexspace

  \item Rank of a Lie algebra, \hyperpage{93}
  \item Reductive
    \subitem $AdS_n$, \hyperpage{107}
  \item Regular
    \subitem representation, \hyperpage{59}
  \item Representation
    \subitem regular left, \hyperpage{59}
  \item Right-invariant
    \subitem vector field, \hyperpage{94}
  \item Root
    \subitem restricted, \hyperpage{93}
    \subitem space
      \subsubitem decomposition, \hyperpage{93}

  \indexspace

  \item Singular point, \hyperpage{15}
  \item Singularity, \hyperpage{15}
  \item Symplectic
    \subitem group, \hyperpage{113}
    \subitem Lie algebra, \hyperpage{114}
    \subitem manifold, \hyperpage{114}
    \subitem vector space, \hyperpage{113}
  \item Symplectomorphism, \hyperpage{114}

  \indexspace

  \item Time orientation, \hyperpage{18}

  \indexspace

  \item WKB quantization, \hyperpage{60}

\end{theindex}

\end{document}